\documentclass[a4paper,10pt]{scrartcl}
\pdfoutput=1
\usepackage{ifpdf}
\usepackage[english]{babel}
\usepackage[utf8]{inputenc}
\usepackage[T1]{fontenc}
\usepackage{lmodern}
\usepackage{enumerate}
\usepackage{geometry}
\usepackage{amsfonts,amsmath,amssymb,amsopn}
\usepackage{amsthm}
\usepackage{mathtools}
\usepackage{stmaryrd}
\usepackage{nicefrac}
\usepackage{exscale}
\usepackage{empheq}
\usepackage{mathrsfs}

\usepackage[numbers,square]{natbib}
\usepackage{url} 
\usepackage[colorlinks=true,pdfpagelabels,unicode]{hyperref}
\hypersetup{linkcolor=blue, urlcolor=blue, citecolor=blue}
\usepackage[dvipsnames,svgnames,table]{xcolor}
\usepackage{orcidlink} 

\allowdisplaybreaks

\theoremstyle{plain}

\newtheoremstyle{theo}
	{3pt} 
	{3pt} 
	{\itshape} 
	{} 
		{\bfseries} 
	{\\} 
	{ } 
	{\thmname{#1}\thmnumber{ #2.}\thmnote{ - #3}} 
\theoremstyle{theo}

\newtheorem{definition}{Definition}[section]
\newtheorem{lemma}[definition]{Lemma}
\newtheorem{theorem}[definition]{Theorem}
\newtheorem{corollary}[definition]{Corollary}

\newenvironment{bew}{\begin{proof}[\bfseries Proof:]}{\end{proof}}

\newtheoremstyle{remark}
	{3pt} 
	{3pt} 
	{} 
	{} 
		{\bfseries} 
	{} 
	{ } 
	{\thmname{#1}\thmnumber{ #2.}\thmnote{ - #3}} 
\theoremstyle{remark}

\DeclareMathOperator{\bomega}{\overline{\Omega}}
\DeclareMathOperator{\romega}{\partial\Omega}

\DeclareMathOperator{\intd}{d\!}
\DeclareMathOperator{\wto}{\rightharpoonup}
\DeclareMathOperator{\wsto}{\stackrel{\star}{\wto}}
\newcommand{\fracpow}{\varrho}
\newcommand{\epsi}{\varepsilon}
\newcommand{\eps}{\varepsilon}

\newcommand{\nep}{n_\epsi}
\newcommand{\cep}{c_\epsi}
\newcommand{\uep}{u_\epsi}
\newcommand{\wep}{w_\epsi}
\newcommand{\Dep}{D_\epsi}
\newcommand{\Sep}{S_\epsi}
\newcommand{\fep}{f_\epsi}
\newcommand{\gep}{g_\epsi}

\newcommand{\GNI}{Gagliardo--Nirenberg inequality}
\newcommand{\intttau}{\int_t^{t+\tau}\!}
\newcommand{\inttpone}{\int_t^{t+1}\!\!} 

\newcommand{\into}[1]{\int_0^{#1}\!}

\newcommand{\intoT}{\into{T}}

\newcommand{\intomega}{\int_{\Omega}\!} 
\newcommand{\intoTomega}{\intoT\!\intomega}

\newcommand{\intinfomega}{\int_0^\infty\!\!\intomega}

\newcommand{\intromega}{\int_{\romega}\!} 

\newcommand{\Lo}[1][1]{L^{#1}(\Omega)} 
\newcommand{\W}[1][1,2]{W^{#1}(\Omega)}

\newcommand{\LSp}[2]{L^{#1\;\!}\!\left(#2\right)} 
\newcommand{\LSpn}[2]{L^{#1\;\!}\!(#2)}

\newcommand{\LSploc}[2]{L_{loc}^{#1}\!\left(#2\right)} 

\newcommand{\CSp}[2]{C^{#1}\!\left(#2\right)}

\newcommand{\CSpnl}[2]{C^{#1}\!\,(#2)} 
\newcommand{\CSploc}[2]{C_{loc}^{#1}\!\left(#2\right)}

\newcommand{\R}{\mathbb{R}}
\newcommand{\N}{\mathbb{N}}
\newcommand{\F}{\mathcal{F}}
\newcommand{\G}{\mathcal{G}}
\newcommand{\Po}{\mathcal{P}}
\newcommand{\hb}{\hat{b}}
\newcommand{\Fep}{\F_{\epsi}}
\newcommand{\Gep}{\G_{\epsi}}
\newcommand{\nfrac}[2]{{\nicefrac{#1}{#2}}}

\newcommand{\M}{\mathcal{M}}
\newcommand{\dimN}{N}

\DeclareMathOperator*{\esssup}{ess\, sup}

\makeatletter
 \def\@lbibitem[#1]#2#3{%
  \if\relax\@extra@b@citeb\relax\else
    \@ifundefined{br@#2\@extra@b@citeb}{}{%
     \@namedef{br@#2}{\@nameuse{br@#2\@extra@b@citeb}}%
    }%
  \fi
  \@ifundefined{b@#2\@extra@b@citeb}{%
   \def\NAT@num{}%
  }{%
   \NAT@parse{#2}%
  }%
  \def\NAT@tmp{#1}%
  \expandafter\let\expandafter\bibitemOpen\csname NAT@b@open@#2\endcsname
  \expandafter\let\expandafter\bibitemShut\csname NAT@b@shut@#2\endcsname
  \@ifnum{\NAT@merge>\@ne}{%
   \NAT@bibitem@first@sw{%
    \@firstoftwo
   }{%
    \@ifundefined{NAT@b*@#2}{%
     \@firstoftwo
    }{%
     \expandafter\def\expandafter\NAT@num\expandafter{\the\c@NAT@ctr}%
     \@secondoftwo
    }%
   }%
  }{%
   \@firstoftwo
  }%
  {%
   \global\advance\c@NAT@ctr\@ne
   \@ifx{\NAT@tmp\@empty}{\@firstoftwo}{%
    \@secondoftwo
   }%
   {%
    \expandafter\def\expandafter\NAT@num\expandafter{\the\c@NAT@ctr}%
    \global\NAT@stdbsttrue
   }{}%
   \bibitem@fin
   \item[\href{#3}{\hfil\NAT@anchor{#2}{\NAT@num}}]
   \global\let\NAT@bibitem@first@sw\@secondoftwo
   \NAT@bibitem@init
  }%
  {%
   \NAT@anchor{#2}{}%
   \NAT@bibitem@cont
   \bibitem@fin
  }%
  \@ifx{\NAT@tmp\@empty}{%
    \NAT@wrout{\the\c@NAT@ctr}{}{}{}{#2}%
  }{%
    \expandafter\NAT@ifcmd\NAT@tmp(@)(@)\@nil{#2}%
  }%
}
\makeatother

\author{
Tobias Black\footnote{tblack@math.upb.de}\ \orcidlink{0000-0001-9963-0800}\\
{\small Institute of Mathematics, Paderborn University,}\\[-5pt]
{\small 33098 Paderborn, Germany}
}
\title{Very mild diffusion enhancement\\ and singular sensitivity:\\Existence of bounded weak solutions in a two-dimensional chemotaxis-Navier--Stokes system}
\date{}

\begin{document}
\maketitle
\begin{abstract}
\noindent
{\textbf{Abstract.} We consider an initial-boundary value problem for the chemotaxis-Navier--Stokes system
\begin{align*}
\left\{
\begin{array}{r@{\,}l@{\quad}l@{\quad}l@{\,}c}
n_{t}+u\cdot\!\nabla n&=\nabla\cdot\big(D(n)\nabla n-nS(x,n,c)\cdot\nabla c\big),\ &x\in\Omega,& t>0,\\
c_{t}+u\cdot\!\nabla c&=\Delta c-cn,\ &x\in\Omega,& t>0,\\
u_{t}+(u\cdot\nabla)u&=\Delta u+\nabla P+n\nabla\Phi,\quad \nabla\cdot u=0,\ &x\in\Omega,& t>0,\\
\big(D(n)\nabla n-n&S(x,n,c)\cdot\nabla c)\cdot\nu=\nabla c\cdot\nu=0,\ u=0,\ &x\in\romega,& t>0,\\
n(\cdot,0)=n_0,\ &\ c(\cdot,0)=c_0,\ u(\cdot,0)=u_0,\ &x\in\Omega.
\end{array}\right.
\end{align*}
in a smoothly bounded domain $\Omega\subset\R^2$. Assuming $S:\overline{\Omega}\times[0,\infty)\times(0,\infty)\rightarrow \R^{2\times 2}$ to be sufficiently regular and such that with $\gamma\in[0,\tfrac56]$ and some non-decreasing $S_0:(0,\infty)\to(0,\infty)$, we have
\begin{align*}
\big|S(x,n,c)\big|\leq \frac{S_0(c)}{c^\gamma}\quad\text{for all }(x,n,c)\in\overline{\Omega}\times[0,\infty)\times(0,\infty),
\end{align*} 
we show that if $D:[0,\infty)\to[0,\infty)$ is suitably regular and positive throughout $(0,\infty)$, then for all $M>0$ one can find $L(M)>0$ such that whenever
$$\liminf_{n\to\infty} D(n)>L\quad\text{and}\quad \liminf_{n\searrow0}\frac{D(n)}{n}>0$$
are satisfied and the initial data $(n_0,c_0,u_0)$ are suitably regular and satisfy $\|c_0\|_{\Lo[\infty]}\leq M$ there is a global and bounded weak solution for the initial-boundary value problem above. Under the additional assumption of $D(0)>0$ this solution is moreover a classical solution of the same problem.

}\medskip



\end{abstract}

\newpage
\section{Introduction}\label{sec1:intro}
Colonies of the aerobic bacteria \emph{Bacillus subtilis} suspended in a drop of water have been observed to form plume-like areas with higher concentration of cells near the air-water contact line (\cite{tuval2005bacterial}). In addition to their experimental findings on this pattern formation, the authors also proposed a coupled chemotaxis-fluid model of the form
\begin{align}\label{problem}
\left\{
\begin{array}{r@{\,}l@{\quad}l@{\quad}l@{\,}c}
n_{t}+u\cdot\!\nabla n&=\nabla\cdot\big(D(n)\nabla n-nS(x,n,c)\cdot\nabla c\big),\ &x\in\Omega,& t>0,\\
c_{t}+u\cdot\!\nabla c&=\Delta c-cn,\ &x\in\Omega,& t>0,\\
u_{t}+(u\cdot\nabla)u&=\Delta u+\nabla P+n\nabla\Phi,\quad \nabla\cdot u=0,\ &x\in\Omega,& t>0,\\
\big(D(n)\nabla n-n&S(x,n,c)\cdot\nabla c)\cdot\nu=\nabla c\cdot\nu=0,\ u=0,\ &x\in\romega,& t>0,\\
n(\cdot,0)=n_0,\ &\ c(\cdot,0)=c_0,\ u(\cdot,0)=u_0,\ &x\in\Omega.
\end{array}\right.
\end{align}
to mathematically describe the spatio-temporal evolution of the population. Herein, $n,c,u,P$ denote the density of the bacteria, the concentration of the oxygen, the fluid velocity and the associated pressure, respectively, and $\Omega\subset\R^\dimN$, $\dimN\in\{2,3\}$ is a smoothly bounded domain. In its original  formulation in \cite{tuval2005bacterial} it was assumed that $D$  is constant and $S$ a real-valued function solely depending on $c$. Since its first introduction different nuances on the influence of the bacterial migration have been considered. Two popular extensions consist of tensor-valued chemotactic-sensitivity functions incorporating a rotational component to the cross-diffusive flux (\cite{XueOthmer-Multiscale-SIAM09,xueMacroscopicEquationsBacterial2015,Winkler18_stokesrot,CaoWang-GlobClass-DCDSB15,caolan16_smalldatasol3dnavstokes}) and nonlinear diffusion operators of porous medium type (\cite{FrancescoLorzMarkowich-DCDS10,TaoWin-GlobExAndBdd-DCDS12,TaoWin-LocBddGlobSol-AnnInstHP13,duanNoteGlobalExistence2014}). See also the survey \cite{BBWT15} for additional references and results on the fluid-free counterparts to \eqref{problem}. 

A major facet of the analysis of chemotaxis systems consists in ruling out the possibility blow-up phenomena, which famously occur in the classical fluid-free Keller--Segel systems obtained by replacing the second equation with $c_t+u\cdot\nabla c=\Delta c-c+n$ and setting $u\equiv0$ (\cite{KS70,NagSen_JAMA98,Win13pure,CieslakStinner-JDE12}). In contrast, the setting with signal consumption is generally assumed to be less prone to emit generalized solutions blowing up in finite time (see e.g. \cite{Tao-consumption_JMAA11,win_fluid_CPDE12,win-stab2d-ArchRatMechAna12,win15_chemorot,caoGlobalintimeBoundedWeak2014}) as long as the chemotactic sensitivity is suitably well-behaved, because of the inherent boundedness of the signal concentration entailed by the consumption process. Nevertheless, global classical solutions have only been constructed under the assumption of suitably small initial data in this setting. If, however, perception of the signal is modeled after the Weber--Fechner law (\cite{HP09,ROSEN1978}) the sensitivity function is of a singular type, i.e. scalar valued $S(x,n,c)=\frac{1}{c}$, which together with the depletion of signal leads to difficult mathematical obstacles in the analysis. (See also the survey \cite{lankeitDepletingSignalAnalysis2023} for more details.)

\textbf{Counteracting singularities in chemoattraction.} Considering sensitivities which are possibly less severe than the singularity imposed by the Weber--Fechner law, i.e. letting $S(x,n,c)=\frac{1}{c^\gamma}$ with $\gamma\in[0,1]$, a variety of structural approaches to deal with the interplay between consumption of signal and singular sensitivity can be found in recent literature. For instance, in two-dimensional domains global classical solutions for sufficiently small initial data were obtained in \cite{liuGlobalClassicalSolvability2023} for a chemotaxis-fluid system with linear diffusion and $\gamma\in[0,1)$. In a setting with diffusion enhanced movement, i.e. taking porous medium-type diffusion $D(n)=n^{m-1}$ and suitably large $m>1$, global and bounded weak solutions have been established for two- and three-dimensional domains in \cite{wangGlobalBoundedSolution2023} under the assumption of $m>\frac{3\dimN-2}{2\dimN}$ without the need for a smallness condition on the initial data.

A different approach (\cite{fuestChemotaxisFluidSystems2023}) consists of limiting the speed of consumption at low levels of chemical concentration by replacing $-nc$ in the second equation by a term acting superlinear with respect to $c$ near zero, e.g. $-c^2 n$. In this slow consumption case, \cite{fuestChemotaxisFluidSystems2023} established global generalized solutions to a two-dimensional chemotaxis-fluid system with $D(n)\equiv 1$ and $S(x,n,c)=\frac{1}{c}$ for arbitrarily large initial data.

Limiting the range of $\gamma$ while still keeping the singular form also seems to act in favor for global solutions. For example the existence of global weak solutions in a chemotaxis-fluid model with linear diffusion and $\gamma\in\big(0,\frac{1}{2}\sqrt{\frac{\dimN+8}{\dimN}}-\frac{1}{2}\big)$ can be found in \cite[Remark 1.2]{kimGlobalSolutionsChemotaxisfluid2023}. For sufficiently small values of $\gamma$ even weaker assumptions on the diffusion enhancement of the bacteria lead to global solutions, as witnessed by \cite{winklerChemotaxisStokesInteractionVery2022} for the case $\gamma=\frac{1}{2}$. There, it was shown in a 3D-Stokes-variant of \eqref{problem} that for all $M>0$ one can find a constant $L=L(M)>0$ such that if $D$ is sufficiently smooth and satisfies
$$\liminf_{n\to\infty} D(n)>L\quad\text{and}\quad\liminf_{n\searrow 0}\frac{D(n)}{n}>0$$
and the initial concentration fulfills $\|c_0\|_{\Lo[\infty]}\leq M$, then there exists a global bounded weak solution. In particular, arbitrarily slow diffusion processes $D(n)$ can be accommodated as long as they have divergent behavior near $n=\infty$.

Having the results above in mind, we dedicate our investigations to the following. On the one hand, we want to extend \cite{winklerChemotaxisStokesInteractionVery2022} to the two-dimensional Navier--Stokes setting despite the more limited fluid-regularity in this setting, and, on the other hand, we want to determine if power-type singularities in the chemotactic sensitivity beyond the prototypical case $\gamma=\frac12$ featured in \cite{winklerChemotaxisStokesInteractionVery2022} can be counteracted by the mild diffusion enhancement expressed above.

\textbf{Main results.} 
Before stating our main theorem, let us briefly fix some notations. Here and below, $\Omega\subset\R^2$ is a smoothly bounded domain and we let $A_r$ denote the Stokes Operator $A_r:=-\mathcal{P}\Delta$ on $L^r(\Omega;\R^2)$ with its domain $D(A_r):=W^{2,r}(\Omega;\R^2)\cap W_{0}^{1,r}(\Omega)\cap L_{\sigma}^{r}(\Omega)$, where $L_{\sigma}^{r}(\Omega):=\{\varphi\in L^r(\Omega;\R^2)\,\vert\,\nabla\cdot\varphi=0\}$ represents the solenoidal vector fields in $L^r(\Omega;\R^2)$ and $\mathcal{P}$ stands for the Helmholtz projection from $L^r(\Omega;\R^2)$ onto $L_\sigma^r(\Omega)$. Moreover, for $\fracpow>0$ we denote by $A_r^\fracpow$ the fractional power of $A_r$. (See also \cite{gig85,sohr}.)
\begin{theorem}\label{theo1}
Let $M>0$, $\gamma\in[0,\frac56]$, $\Phi\in\W[2,\infty]$ and assume $\Omega\subset\R^2$ to be a bounded domain with smooth boundary. Suppose the tensor-valued function $S\in\CSp{2}{\bomega\times[0,\infty)\times(0,\infty);\R^{2\times2}}$ fulfills
\begin{align}\label{eq:prop-S}
|S(x,n,c)|\leq\frac{S_0(c)}{c^\gamma}\quad\text{for all }(x,n,c)\in\bomega\times[0,\infty)\times(0,\infty),
\end{align}
with some non-decreasing $S_0:(0,\infty)\to(0,\infty)$. Then, there exists $L=L(M)>0$ such that whenever
$D\in\bigcup_{\theta\in(0,1)}\CSp{\theta}{[0,\infty)}\cap\CSp{2}{(0,\infty)}$ is positive on $(0,\infty)$ and satisfies
\begin{align}\label{eq:liminf}
\liminf_{n\to\infty} D(n)>L\quad\text{and}\quad\liminf_{n\searrow0}\frac{D(n)}{n}>0,
\end{align}
and the initial data $(n_0,c_0,u_0)$ fulfill
\begin{align}\label{IR}
\begin{cases}
n_0\in\W[1,\infty]\ \text{with } n_0\geq 0\text{ and }n_0\not\equiv0\text{ in }\bomega,\\
c_0\in\W[1,\infty]\ \text{with } c_0>0\text{ in }\bomega,\\
u_0\in D(A^\fracpow_r)\ \text{for some }\fracpow\in(\tfrac{1}{2},1)\text{ and all }r\in(1,\infty), 
\end{cases}
\end{align}
with $\|c_0\|_{\Lo[\infty]}\leq M$ and, in the case of $\gamma\in(\tfrac12,\tfrac56]$ also $\|n_0\|_{\Lo[1]}\leq M$, there is a global weak solution $(n,c,u)$ which solves \eqref{problem} in the sense of Definition~\ref{def:weak_sol}.
This weak solution satisfies the boundedness properties
\begin{align}\label{eq:bounded}
\esssup_{t>0}\Big(\big\|n(\cdot,t)\big\|_{\Lo[\infty]}+\big\|c(\cdot,t)\big\|_{\W[1,\infty]}+\big\|u(\cdot,t)\big\|_{\W[1,\infty]}\Big)<\infty.
\end{align}
If, additionally, $D(0)>0$, then the solution features improved regularity properties, i.e.
\begin{align}\label{classical}
\begin{cases}
n\in\CSp{0}{\bomega\times[0,\infty)}\cap\CSp{2,1}{\bomega\times(0,\infty)},\\
c\in\CSp{0}{\bomega\times[0,\infty)}\cap\CSp{2,1}{\bomega\times(0,\infty)},\\
u\in\CSp{0}{\bomega\times[0,\infty);\R^2}\cap\CSp{2,1}{\bomega\times(0,\infty);\R^2},
\end{cases}
\end{align}
and there is $P\in\CSp{1,0}{\bomega\times(0,\infty)}$ such that $(n,c,u,P)$ forms a classical solution of \eqref{problem} in $\Omega\times(0,\infty)$.
\end{theorem}

\begin{corollary}\label{cor:1}
Let $\gamma\in[0,\tfrac56]$ and assume $\Omega\subset\R^2$ to be a bounded domain with smooth boundary. Suppose $S\in\CSp{2}{\bomega\times[0,\infty)\times(0,\infty);\R^{2\times2}}$ fulfills \eqref{eq:prop-S} with some non-decreasing $S_0:(0,\infty)\to(0,\infty)$ and $D\in\bigcup_{\theta\in(0,1)}\CSp{\theta}{[0,\infty)}\cap\CSp{2}{(0,\infty)}$ is positive on $(0,\infty)$, with $\liminf_{n\searrow 0}\frac{D(n)}{n}>0$ and
$$D(n)\to\infty\quad\text{as }n\to\infty.$$ Then, for any choice of $(n_0,c_0,u_0)$ satisfying \eqref{IR} there is a global weak solution $(n,c,u)$ which solves \eqref{problem} in the sense of Definition~\ref{def:weak_sol} and remains bounded in the sense of \eqref{eq:bounded}. If, additionally, $D(0)>0$, then $n,c,u$ are of the classes specified in \eqref{classical} and there is some $P\in\CSp{1,0}{\bomega\times(0,\infty)}$ such that $(n,c,u,P)$ forms a classical solution of \eqref{problem} in $\Omega\times(0,\infty)$.
\end{corollary}

In particular, Corollary~\ref{cor:1} applies to porous medium type diffusion of the form $D(n)=n^{m-1}$ with $m\in(1,2]$. Whereas, for linear diffusion, Theorem~\ref{theo1} on the one hand partially recovers the small-data global bounded weak solutions established in \cite{liuGlobalClassicalSolvability2023} (for $\gamma\in[0,1)$), and on the other hand assuming that $D(n)\equiv D_0\geq L$ for some sufficiently large $L$, which in the case $\gamma\leq\frac12$ solely depends on the initial chemical concentration, Theorem~\ref{theo1} asserts the existence of a global bounded classical solution. Moreover, in case of porous medium type diffusion $D(n)=n^{m-1}$, $m>1$, the global weak solution is actually continuous.
\begin{corollary}\label{cor:2}
Let $\gamma\in[0,\tfrac56]$ and assume that $\Omega\subset\R^2$ is a domain with smooth boundary. Assume that $D(n)=n^{m-1}$ for some $m\in(1,2]$ and that $S\in\CSp{2}{\bomega\times[0,\infty)\times(0,\infty);\R^{2\times 2}}$ fulfills \eqref{eq:prop-S} with some non-decreasing $S_0:(0,\infty)\to(0,\infty)$. Then, for any choice of $(n_0,c_0,u_0)$ satisfying \eqref{IR} the global weak solution $(n,c,u)$ to \eqref{problem} obtained in Theorem~\ref{theo1} is continuous.
\end{corollary}
\smallskip

\textbf{Assembling an appropriate energy-like functional.} After making use of a result related to the Trudinger--Moser inequality (cf. Lemma~\ref{lem:mos-trud-2}) to obtain spatio-temporal bounds on $n\ln(n)$ and time-uniform bounds on $\intomega u^2$, we can adapt the approach of \cite{winklerChemotaxisStokesInteractionVery2022} to our Navier--Stokes setting. For this, following the philosophy of \cite{winklerChemotaxisStokesInteractionVery2022}, we set $D_{2}(n):=\int_0^n\!\int_0^s D(\sigma)\intd\sigma\intd s$ and make the time-evolution of $\intomega D_{2}(n)$ the focal point of our analysis (see Section~\ref{sec3:energy}). Indeed, for small values of $\gamma$ some control on the ill-signed taxis term $\intomega\frac{n^2|\nabla c|^2}{c^{2\gamma}}$ can be established by considering a combination of functionals of the form $$\intomega D_{2}(n)+\intomega\frac{n|\nabla c|^2}{c}+\intomega\frac{|\nabla c|^4}{c^3}.$$ For $\gamma>\frac{1}{2}$, however, the well-signed term $\intomega\frac{n^2|\nabla c|^2}{c}$ is insufficient to control this taxis-term and adjusting to $\intomega\frac{n|\nabla c|^2}{c^{\rho}}+\intomega\frac{|\nabla c|^4}{c^{\rho+2}}$ for some $\rho>1$ destroys the subtle cancellations utilized in Lemma~\ref{lem:ener-part3}. Accordingly, for $\gamma>\frac{1}{2}$ a slightly different approach is needed (Section~\ref{sec4:large-g}), which necessitates the additional condition that also $\|n_0\|_{\Lo[1]}\leq M$ is satisfied. Under this assumption splitting the taxis-related term by means of Young's inequality control is provided solely by the differential inequality for $\intomega\frac{|\nabla c|^4}{c^3}$. Afterwards, having control on $\intomega|\nabla c|^4$ from our functional inequality at hand, we can obtain the boundedness $n$ by an iterated testing procedure for $(n-1)_+$ (Lemma~\ref{lem:moser-it}). These bounds at hand, an Aubin--Lions type argument will establish the existence of the weak solution (Lemma~\ref{lem:aubin} as commonly witnessed in previous chemotaxis precedents).
\smallskip

For the remainder of the work, we will always assume that $S$, $D$ and the initial data to be of the regularity classes specified in Theorem~\ref{theo1} and \eqref{IR}, respectively.

\setcounter{equation}{0} 
\section{Approximate problems and basic estimates}\label{sec2:}
In order to make \eqref{problem} accessible to our arguments, we approximate both the diffusion coefficient $D$ and the sensitivity function $S$ in the following way. We introduce the perturbation parameter $\epsi\in(0,1)$, and pick $\Dep\in\CSp{2}{[0,\infty)}$ such that 
\begin{align}\label{eq:Dep}
\epsi\leq \Dep(n) \quad\text{and}\quad D(n)\leq\Dep(n)\leq D(n)+2\epsilon\quad\text{for all }n\geq 0.
\end{align} 
Moreover, we chose $\rho_\epsi\in C_0^\infty(\Omega)$ and $\chi_\epsi\in C_0^\infty([0,\infty))$ satisfying  
\begin{align*}
0\leq \rho_\epsi\leq 1\text{ in }\Omega,\qquad 0\leq\chi_\epsi\leq 1\text{ in }[0,\infty)
\end{align*}
and
\begin{align*}
\rho_\epsi\nearrow 1\text{ in }\Omega,\quad\text{and}\quad \chi_\epsi\nearrow 1\text{ in }[0,\infty)\quad\text{as }\epsi\searrow0.
\end{align*}
Then, with $\Sep\in\CSp{2}{\bomega\times[0,\infty)^2;\R^{2\times2}}$ defined via 
\begin{align}\label{eq:Sep}
\Sep(x,n,c)=\rho_\epsi(x)\chi_\epsi(n)S(x,n,c+\epsi),\quad\text{for }x\in\bomega,\ n\geq0,\ c\geq0,
\end{align}
we consider the non-degenerate approximate problems of the form
\begin{align}\label{approxprob}
\left\{
\begin{array}{r@{\,}l@{\quad}l@{\quad}l@{\,}c}
n_{\epsi t}+\uep\cdot\!\nabla \nep&=\nabla\cdot\big(\Dep(\nep)\nabla \nep-\nep\Sep(x,\nep,\cep)\cdot\nabla \cep\big),\ &x\in\Omega,& t>0,\\
c_{\epsi t}+\uep\cdot\!\nabla \cep&=\Delta \cep-\cep\nep,\ &x\in\Omega,& t>0,\\
u_{\epsi t}+(\uep\cdot\nabla)\uep&=\Delta \uep+\nabla P_\epsi+\nep\nabla\Phi,\quad \nabla\cdot \uep=0,\ &x\in\Omega,& t>0,\\
\frac{\partial\nep}{\partial\nu}=\frac{\partial\cep}{\partial\nu}&=0,\ \uep=0,\ &x\in\romega,& t>0,\\
\nep(\cdot,0)=n_0,\ &\ \cep(\cdot,0)=c_0,\ \uep(\cdot,0)=u_0,\ &x\in\Omega.
\end{array}\right.
\end{align}
Moreover, we infer -- for later purposes -- from \eqref{eq:liminf} and \eqref{eq:Dep} that our choices above entail
\begin{align}\label{eq:s0}
\exists\, s_0=s_0(L)>0\text{ such that }\Dep(n)\geq L\quad\text{for all }n\geq s_0\text{ and }\epsi\in(0,1).
\end{align}

Adapting well-established fixed points arguments, we find that for each $\epsi\in(0,1)$ the problem \eqref{approxprob} has a global classical solution (e.g. \cite[Lemma 2.1]{Win-ct_fluid_3d-CPDE15}, \cite[Lemma 2.1]{tao_winkler_chemohapto11-siam11} and \cite[Lemma 2.2]{Lan17-LocBddGlobSolNonlinDiff-JDE}).

\begin{lemma}\label{lem:locex}
Let $\Omega\subset\R^2$ be a bounded domain with smooth boundary and $\Phi\in\W[2,\infty]$. Suppose that $n_0,c_0$ and $u_0$ comply with \eqref{IR}. Then for any $\epsi\in(0,1)$ there exist functions
\begin{align*}
\nep&\in\CSp{0}{\bomega\times[0,\infty)}\cap\CSp{2,1}{\bomega\times(0,\infty)},\\
\cep&\in\CSp{0}{\bomega\times[0,\infty)}\cap\CSp{2,1}{\bomega\times(0,\infty)},\\
\uep&\in\CSp{0}{\bomega\times[0,\infty);\R^2}\cap\CSp{2,1}{\bomega\times(0,\infty);\R^2},\\
P_\epsi&\in\CSp{1,0}{\bomega\times(0,\infty)},
\end{align*}
which solve \eqref{approxprob} classically in $\Omega\times(0,\infty)$. Moreover, $\nep$ and $\cep$ are positive in $\bomega\times(0,\infty)$.
\end{lemma}

The solutions of \eqref{approxprob} satisfy basic mass-conservation for the first and an $L^\infty$-boundedness property for the second solution components, respectively.

\begin{lemma}\label{lem:mass}
For all $\epsi\in(0,1)$ the solution $(\nep,\cep,\uep)$ of \eqref{approxprob} satisfies
\begin{align*}
\intomega \nep(\cdot,t)=\intomega n_0\quad\text{for all }t>0,
\end{align*}
and
\begin{align*}
\big\|\cep(\cdot,t)\big\|_{\Lo[\infty]}\leq \|c_0\|_{\Lo[\infty]}\quad\text{for all }t>0.
\end{align*}
\end{lemma}

\begin{bew}
The mass-conservation property of $\nep$ follows from direct integration of the first equation in \eqref{approxprob}, whereas the bound for $\|\cep(\cdot,t)\|_{\Lo[\infty]}$ is a consequence of the maximum principle employed to the second equation of \eqref{approxprob}.
\end{bew}

Straightforward testing of the second equation against $\cep^{q-1}$ with $q\in(0,1)$ additionally provides a first -- albeit very mild -- bound on a time-space integral containing the gradient of $\cep$.

\begin{lemma}\label{lem:nab_cep_2-q}
For all $q\in(0,1)$ one can find $C=C(q,n_0,c_0)>0$ such that for all $\epsi\in(0,1)$ the solution $(\nep,\cep,\uep)$ of \eqref{approxprob} satisfies
\begin{align*}
\inttpone\intomega\frac{|\nabla\cep|^2}{\cep^{2-q}}\leq C\quad\text{for all }t>0.
\end{align*}
\end{lemma}

\begin{bew}
Using the second equation of \eqref{approxprob}, the fact that $\frac{\partial\cep}{\partial\nu}=0$ on $\romega\times(0,\infty)$ and the divergence-free property of $\uep$, we calculate
\begin{align*}
-\frac{1}{q}\frac{\intd}{\intd t}\intomega\cep^q=(q-1)\intomega\cep^{q-2}|\nabla\cep|^2+\intomega\nep\cep^q\quad\text{on }(0,\infty)
\end{align*}
for all $\epsi\in(0,1)$. In view of Lemma~\ref{lem:mass}, we can further estimate the last term on the right, to conclude
\begin{align*}
-\frac{\intd}{\intd t}\intomega\cep^q+q(1-q)\intomega\frac{|\nabla\cep|^2}{\cep^{2-q}}\leq q\|c_0\|_{\Lo[\infty]}^q\intomega n_0\quad\text{on }(0,\infty)
\end{align*}
for all $\epsi\in(0,1)$. Integrating over $[t,t+1]$ therefore yields
\begin{align*}
\intomega\cep^q(\cdot,t)+q(1-q)\inttpone\intomega\frac{|\nabla\cep|^2}{\cep^{2-q}}&\leq\intomega\cep^q(\cdot,t+1)+q\|c_0\|_{\Lo[\infty]}^q\intomega n_0\\&\leq \|c_0\|_{\Lo[\infty]}^q\Big(|\Omega|+ q\intomega n_0\Big)
\end{align*}
for all $t>0$ and $\epsi\in(0,1)$, from which the claim follows in light of $q\in(0,1)$ and the positivity of $\cep$.
\end{bew}

Making use of the second condition in \eqref{eq:liminf} and Lemma~\ref{lem:nab_cep_2-q} we can also establish some feeble time-space information linked to the diffusion coefficient and the gradient of $\nep$. Here, we introduce $\gamma_0\in(\tfrac12,1)$ in order to keep the dependencies of the constant clear of $\gamma$. 

\begin{lemma}\label{lem:st-bound-Dep_nab_nep}
Let $\gamma_0\in(\tfrac12,1)$ and assume that \eqref{eq:prop-S} is satisfied for some $\gamma\in[0,\gamma_0]$. For all $L>0$ there is $C=C(L,n_0,c_0,\gamma_0)>0$ such that if $D$ fulfills \eqref{eq:liminf}, then the solution $(\nep,\cep,\uep)$ of \eqref{approxprob} fulfills
\begin{align*}
\inttpone\intomega\frac{\Dep(\nep)}{(\nep+1)^2}|\nabla\nep|^2\leq C
\end{align*}
for all $t>0$ and all $\epsi\in(0,1)$.
\end{lemma}

\begin{bew}
Relying on the first equation of \eqref{approxprob}, we find from integrating by parts that
\begin{align*}
\frac{\intd}{\intd t}\intomega\ln(\nep+1)&=\intomega\frac{\nabla\cdot\big(\Dep(\nep)\nabla\nep\big)}{\nep+1}-\intomega\frac{\nabla\cdot\big(\nep\Sep(x,\nep,\cep)\cdot\nabla\cep\big)}{\nep+1}-\intomega\frac{\uep\cdot\nabla\nep}{\nep+1}\\
&=\intomega\frac{\Dep(\nep)|\nabla\nep|^2}{(\nep+1)^2}-\intomega\frac{\nep\Sep(x,\nep\cep)\cdot\nabla\cep\cdot\nabla\nep}{(\nep+1)^2}\quad\text{on }(0,\infty)
\end{align*}
for all $\epsi\in(0,1)$, where the integral featuring the fluid component disappeared due to the divergence-free property and homogeneous Dirichlet boundary condition for $\uep$. Accordingly, reordering and employing Young's inequality, we find that
\begin{align*}
\frac{1}{2}\intomega\frac{\Dep(\nep)|\nabla\nep|^2}{(\nep+1)^2}\leq \frac{\intd}{\intd t}\intomega\ln(\nep+1)+\frac{1}{2}\intomega\frac{\nep^2|\Sep(x,\nep,\cep)|^2 |\nabla\cep|^2}{(\nep+1)^2\Dep(\nep)}\quad\text{on }(0,\infty)
\end{align*}
for all $\epsi\in(0,1)$. Integrating from $t$ to $t+1$, we draw on \eqref{eq:prop-S} and the non-decreasing property of $S_0$ to find that
\begin{align}\label{eq:st-bound-Dep_nab_nep-eq1}
\inttpone\intomega\frac{\Dep(\nep)|\nabla\nep|^2}{(\nep+1)^2}&\leq 2\intomega\ln\big(\nep(\cdot,t+1)+1\big)-2\intomega\ln\big(\nep(\cdot,t)+1\big)\nonumber\\
&\qquad+S_0^2\big(\|c_0\|_{\Lo[\infty]}+1\big)\inttpone\intomega\frac{\nep^2|\nabla\cep|^2}{(\nep+1)^2\Dep(\nep)\cep^{2\gamma}}\quad\text{for all }t>0
\end{align}
and all $\epsi\in(0,1)$. Next, with $s_0=s_0(L)>0$ as in \eqref{eq:s0} we have $\Dep(n)\geq L$ for all $n\geq s_0$ and $\epsi\in(0,1)$ and introducing $$\kappa=\kappa(s_0):=\inf_{n\in(0,2s_0)}\frac{D(n)}{n}>0,$$ the positivity of which is due to \eqref{eq:Dep} and the positivity of $D$ in $(0,\infty)$, we obtain from \eqref{eq:Dep} that
\begin{align*}
\frac{n}{\Dep(n)}\leq\frac{n}{D(n)}\leq\frac{1}{\kappa}\quad\text{for all }n\in(0,s_0)\text{ and all }\epsi\in(0,1).
\end{align*}
Thus, splitting the domain in two parts, we can estimate
\begin{align}\label{eq:st-bound-Dep_nab_nep-eq2}
\inttpone\intomega\frac{\nep^2|\nabla\cep|^2}{(\nep+1)^2\Dep(\nep)\cep^{2\gamma}}&=\inttpone\int_{\{\nep\geq s_0\}}\frac{\nep^2|\nabla\cep|^2}{(\nep+1)^2\Dep(\nep)\cep^{2\gamma}}+\inttpone\int_{\{\nep<s_0\}}\frac{\nep^2|\nabla\cep|^2}{(\nep+1)^2\Dep(\nep)\cep^{2\gamma}}\nonumber\\
&\leq \frac{1}{L}\inttpone\int_{\{\nep\geq s_0\}}\frac{\nep^2|\nabla\cep|^2}{(\nep+1)^2\cep^{2\gamma}}+\frac{1}{\kappa}\inttpone\int_{\{\nep<s_0\}}\frac{\nep|\nabla\cep|^2}{(\nep+1)^2\cep^{2\gamma}}\nonumber\\
&\leq \Big(\frac{1}{L}+\frac{1}{\kappa}\Big)\inttpone\intomega\frac{|\nabla\cep|^2}{\cep^{2\gamma}}\quad\text{for all }t>0\text{ and }\epsi\in(0,1).
\end{align}

The assumption $\gamma\leq\gamma_0$ ensures $2\gamma_0-2\gamma\geq 0$, so that by Young's inequality and $\gamma_0<1$
$$\|c_0\|_{\Lo[\infty]}^{2(\gamma_0-\gamma)}\leq \|c_0\|_{\Lo[\infty]}^2+1=:C_1(c_0).$$
Hence, making use of Lemma~\ref{lem:mass} and Lemma~\ref{lem:nab_cep_2-q} for $q=2-2\gamma_0\in(0,1)$, we find $C_2=C_2(n_0,c_0,\gamma_0)>0$ such that
\begin{align}\label{eq:st-bound-Dep_nab_nep-eq3}
\inttpone\intomega\frac{|\nabla\cep|^2}{\cep^{2\gamma}}=\inttpone\intomega\frac{\cep^{2\gamma_0-2\gamma}|\nabla\cep|^2}{\cep^{2\gamma_0}}\leq \|c_0\|_{\Lo[\infty]}^{2(\gamma_0-\gamma)}\inttpone\intomega\frac{|\nabla\cep|^2}{\cep^{2-q}}\leq C_1C_2
\end{align}
for all $t>0$ and $\epsi\in(0,1)$. Aggregating \eqref{eq:st-bound-Dep_nab_nep-eq1}-\eqref{eq:st-bound-Dep_nab_nep-eq3} and utilizing that $\ln(s+1)\leq s$ for $s\geq0$, then entails the existence of $C_3=C_3(L,n_0,c_0,\gamma_0)>0$ satisfying
\begin{align*}
\inttpone\intomega\frac{\Dep(\nep)|\nabla\nep|^2}{(\nep+1)^2}\leq 2\intomega n_0+S_0^2\big(\|c_0\|_{\Lo[\infty]}+1\big)\Big(\frac{1}{L}+\frac{1}{\kappa}\Big)C_1C_2=: C_3\quad\!\text{for all }t>0\text{ and }\epsi\in(0,1),
\end{align*}
which concludes the proof.
\end{bew}

In order to exploit the above information, we will rely on the fact that we are considering a two-dimensional domain and the well-known Trudinger inequality (see e.g. \cite{NSY97}). A direct outcome of said inequality, which is almost in the precise form we will make use of, has been formulated in \cite[Lemma 2.2]{winklerSmallMassSolutionsTwoDimensional2020} and reads as follows.
\begin{lemma}\label{lem:mos-trud}
Let $\Omega\subset\R^2$ be a bounded domain with smooth boundary. Then for any choice of $\eta>0$ there is $K=K(\eta,\Omega)>0$ such that if $0\not\equiv\phi\in\CSp{0}{\bomega}$ is nonnegative and $\psi\in\W[1,2]$, then for each $a>0$,
\begin{align*}
\intomega\phi|\psi|\leq\frac{1}{a}\intomega\phi\ln\Big(\frac{\phi}{\overline{\phi}}\Big)+\frac{(1+\eta)a}{8\pi}\Big(\intomega\phi\Big)\intomega|\nabla\psi|^2+Ka\Big(\intomega\phi\Big)\Big(\intomega|\psi|\Big)^2+\frac{K}{a}\intomega\phi,
\end{align*}
where $\overline{\phi}:=\frac{1}{|\Omega|}\intomega\phi$.
\end{lemma}

Let us reformulate the lemma above to be more suitable for an application to $\nep$ in our setting. In particular, making use of suitably chosen $\psi$ and some properties of $D$, we want to feature a term of the form $\intomega\frac{D(\phi)|\nabla\phi|^2}{(\phi+1)^2}$ on the right-hand side in order to be able to rely on Lemma~\ref{lem:st-bound-Dep_nab_nep} in a more direct fashion.

\begin{lemma}\label{lem:mos-trud-2}
Let $\Omega\subset\R^2$ be a bounded domain with smooth boundary and suppose that $0\not\equiv \phi\in\CSp{0}{\bomega}\cap\W[1,2]$ is nonnegative. Assume that there are $\tilde{s}_0>0$ and $\widetilde{D}\in\CSp{0}{[0,\infty)}$ fulfilling $\widetilde{D}(s)>L$ for all $s\geq \tilde{s}_0$. Then for any choice of $\eta>0$ there is $K=K(\eta,\tilde{s}_0,\Omega)>0$ such that
\begin{align*}
\int\limits_{\{\phi>\tilde{s}_0+1\}}\hspace*{-0.4cm}\phi\ln(\phi+1)\leq\frac{(1+\eta)K}{L}\Big(\intomega\phi\Big)\intomega\frac{\widetilde{D}(\phi)|\nabla\phi|^2}{(\phi+1)^2}+K\Big(\intomega\phi\Big)^3+\bigg(K-\ln\Big(\frac{1}{|\Omega|}\intomega\phi\Big)\bigg)\intomega\phi+K.
\end{align*}
\end{lemma}

\begin{bew}
We introduce a smooth cut-off function $g\in\CSp{\infty}{[0,\infty)}$ satisfying $0\leq g\leq 1$, $g(s)=0$ for $s<\tilde{s}_0$ and $g(s)=1$ for $s>\tilde{s}_0+1$ and $g'(s)\leq 2$ on $(0,\infty)$ and employ Lemma~\ref{lem:mos-trud} for $\psi=g(\phi)\ln(\phi+1)$ and $a=2$ to find $K_1=K_1(\eta,\Omega)>0$ such that
\begin{align}\label{eq:mos-trud-2-eq1}
\intomega\phi g(\phi)\ln(\phi+1)\leq\frac{1}{2}\intomega\phi\ln\Big(\frac{\phi}{\overline\phi}\Big)&+\frac{1+\eta}{4\pi}\Big(\intomega\phi\Big)\Big(\intomega \big(g'(\phi)\big)^2(\phi+1)^2|\nabla\phi|^2+\intomega \frac{g(\phi)|\nabla\phi|^2}{(\phi+1)^2}\Big)\nonumber\\
&+2K_1\Big(\intomega\phi\Big)^3+\frac{K_1}{2}\intomega\phi,
\end{align}
where we made use of $\ln(s+1)\leq s$ for $s\geq0$. Herein, we rely on the properties of $g$ and the fact that $\widetilde{D}(s)>L$ for $s\geq \tilde{s}_0$ to estimate
\begin{align*}
\intomega \big(g'(\phi)\big)^2(\phi+1)^2|\nabla\phi|^2&=\int\limits_{\{\phi\in[\tilde{s}_0,\tilde{s}_0+1]\}}\hspace*{-0.4cm}\big(g'(\phi)\big)^2(\phi+1)^2|\nabla\phi|^2\\
&\leq 4(\tilde{s}_0+2)^2\int\limits_{\{\phi\in[\tilde{s}_0,\tilde{s}_0+1]\}}\hspace*{-0.28cm}\frac{(\phi+1)^2 \widetilde{D}(\phi)}{(\phi+1)^2\widetilde{D}(\phi)}|\nabla\phi|^2\leq \frac{4(\tilde{s}_0+2)^4}{L}\intomega\frac{\widetilde{D}(\phi)|\nabla\phi|^2}{(\phi+1)^2}
\end{align*}
and
\begin{align*}
\intomega\frac{g(\phi)|\nabla\phi|^2}{(\phi+1)^2}\leq\int\limits_{\{\phi\geq \tilde{s}_0\}}\frac{|\nabla\phi|^2}{(\phi+1)^2}\leq \frac{1}{L}\intomega\frac{\widetilde{D}(\phi)|\nabla\phi|^2}{(\phi+1)^2}.
\end{align*}
Hence, in view of $g(s)=1$ for $s\geq \tilde{s}_0+1$, we obtain from \eqref{eq:mos-trud-2-eq1} that
\begin{align}\label{eq:mos-trud-2-eq2}
\int\limits_{\{\phi>\tilde{s}_0+1\}}\hspace*{-0.4cm}\phi \ln(\phi+1)\leq\frac{1}{2}\intomega\phi\ln\Big(\frac{\phi}{\overline\phi}\Big)&+\frac{1+\eta}{4\pi}\Big(\intomega\phi\Big)\frac{1+4(\tilde{s}_0+2)^4}{L}\intomega\frac{\widetilde{D}(\phi)|\nabla\phi|^2}{(\phi+1)^2}\nonumber\\
&+2K_1\Big(\intomega\phi\Big)^3+\frac{K_1}{2}\intomega\phi.
\end{align}
Further estimating
\begin{align*}
\frac{1}{2}\intomega\phi\ln\Big(\frac{\phi}{\overline{\phi}}\Big)&=\frac{1}{2}\intomega\phi\ln\phi-\frac{1}{2}\ln\Big(\frac{1}{|\Omega|}\intomega\phi\Big)\intomega\phi\\
&\leq \frac{1}{2}\int\limits_{\{\phi>\tilde{s}_0+1\}}\hspace*{-0.3cm}\phi\ln(\phi+1)+\frac{1}{2}\int\limits_{\{\phi\leq \tilde{s}_0+1\}}\hspace*{-0.3cm}\phi\ln\phi-\frac{1}{2}\ln\Big(\frac{1}{|\Omega|}\intomega\phi\Big)\intomega\phi\\
&\leq \frac{1}{2}\int\limits_{\{\phi>\tilde{s}_0+1\}}\hspace*{-0.3cm}\phi\ln(\phi+1)+\frac12|\Omega|(\tilde{s}_0+1)^2-\frac{1}{2}\ln\Big(\frac{1}{|\Omega|}\intomega\phi\Big)\intomega\phi,
\end{align*}
plugging this back into \eqref{eq:mos-trud-2-eq2} and rearranging finally yields the claim upon choosing \[K:=\max\Big\{\frac{1+4(\tilde{s}_0+1)^4}{2\pi},4K_1,|\Omega|(\tilde{s}_0+1)^2\Big\}.\qedhere\]
\end{bew}

These preparations at hand, we can directly refine the spatio-temporal bound of Lemma~\ref{lem:st-bound-Dep_nab_nep} together with the mass-conservation of $\nep$ into boundedness of the quantity $\int_t^{t+1}\!\intomega\nep\ln(\nep+1)$.

\begin{lemma}\label{lem:st-bound-neplnnep}
Let $\gamma_0\in(\tfrac12,1)$ and assume that \eqref{eq:prop-S} is satisfied for some $\gamma\in[0,\gamma_0]$. For all $L>0$ there is $C=C(L,n_0,c_0,\gamma_0)>0$  such that if $D$ fulfills \eqref{eq:liminf}, then for all $\epsi\in(0,1)$
\begin{align*}
\inttpone\intomega\nep\ln(\nep+1)\leq C\quad\text{for all }t>0.
\end{align*}
\end{lemma}

\begin{bew}
With $s_0=s_0(L)>0$ as considered in \eqref{eq:s0}, the prerequisites of Lemma~\ref{lem:mos-trud-2} are satisfied for $\phi=\nep$, $D=\Dep$ and $\eta=1$, so that an employment of Lemma~\ref{lem:mos-trud-2} in combination with the mass conservation established in Lemma~\ref{lem:mass} yields $C_1=C_1(L,n_0)>0$ such that
\begin{align*}
\int\limits_{\{\nep>s_0+1\}}\hspace*{-0.3cm}\nep\ln(\nep+1)\leq C_1\intomega\frac{\Dep(\nep)|\nabla\nep|^2}{(\nep+1)^2}+C_1\quad\text{for all }t>0\text{ and }\epsi\in(0,1).
\end{align*}
Moreover, we have
\begin{align*}
\int\limits_{\{\nep\leq s_0+1\}}\hspace*{-0.3cm}\nep\ln(\nep+1)\leq \int\limits_{\{\nep\leq s_0+1\}}\hspace*{-0.3cm}\nep^2\leq (s_0+1)^2|\Omega|\quad\text{for all }t>0\text{ and }\epsi\in(0,1),
\end{align*}
and accordingly find $C_2=C_2(L,n_0)>0$ satisfying
\begin{align*}
\inttpone\intomega\nep\ln(\nep+1)\leq C_2\inttpone\intomega\frac{\Dep(\nep)|\nabla\nep|^2}{(\nep+1)^2}+C_2\quad\text{for all }t>0\text{ and }\epsi\in(0,1).
\end{align*}
The claim now evidently follows from Lemma~\ref{lem:st-bound-Dep_nab_nep}.
\end{bew}

Another useful tool for deriving additional boundedness properties will be the following lemma, which in this form is a verbatim copy of \cite[Lemma 3.4]{winklerThreedimensionalKellerSegel2019}.

\begin{lemma}\label{lem:ode}
Let $t_0\in\R$, $T\in(t_0,\infty]$ and suppose that the nonnegative function $h\in\LSploc{1}{\R}$ has the property that there exist $\tau>0$ and $b>0$ such that
\begin{align*}
\frac{1}{\tau}\intttau h(s)\intd s\leq b\quad\text{for all }t\in(t_0,T).
\end{align*}
Then for any choice of $a>0$ we have
\begin{align*}
\int_{t_0}^{t}e^{-a(t-s)}h(s)\intd s\leq\frac{b\tau}{1-e^{-a\tau}}\quad\text{for all }t\in[t_0,T).
\end{align*}
Consequently, if $y\in\CSp{0}{[t_0,T)}\cap\CSp{1}{(t_0,T)}$ has the property that
\begin{align*}
y'(t)+ay(t)\leq h(t)\quad\text{for all }t\in(t_0,T),
\end{align*}
then
\begin{align*}
y(t)\leq e^{-a(t-t_0)}y(t_0)+\frac{b\tau}{1-e^{-a\tau}}\quad\text{for all }t\in[t_0,T),
\end{align*}
and in particular
\begin{align*}
y(t)\leq y(t_0)+\frac{b\tau}{1-e^{-a\tau}}\quad\text{for all }t\in[t_0,T).
\end{align*}
\end{lemma}

Using Lemma~\ref{lem:st-bound-neplnnep} and Lemma~\ref{lem:mos-trud} in synergy with the ODE-lemma presented above, we can now draw on reasoning as undertaken in e.g. \cite[Lemma 2.5]{winklerLpboundLocalSensing} and \cite[Lemma 3.4]{heihoffTwoNewFunctional2023} to conclude this section with establishing uniform bounds for the fluid-component.

\begin{lemma}\label{lem:uep-bound}
Let $\gamma_0\in(\tfrac12,1)$ and assume \eqref{eq:prop-S} is fulfilled for some $\gamma\in[0,\gamma_0]$. For all $L>0$ there is $C=C(L,n_0,c_0,u_0,\gamma_0)>0$ such that if $D$ satisfies \eqref{eq:liminf}, then
\begin{align*}
\intomega|\uep(\cdot,t)|^2+\inttpone\intomega|\nabla\uep|^2\leq C\quad\text{for all }t>0\text{ and all }\epsi\in(0,1).
\end{align*}
\end{lemma}

\begin{bew}
Direct calculations drawing on the third equation of \eqref{approxprob} yield that
\begin{align}\label{eq:uep-bound-eq1}
\frac{1}{2}\frac{\intd}{\intd t}\intomega|\uep|^2\leq -\intomega|\nabla\uep|^2+\|\nabla \Phi\|_{\Lo[\infty]}\intomega|\uep|\nep
\end{align}
holds for all $\epsi\in(0,1)$ and $t>0$. In view of a Poincaré-inequality, we can pick $C_1=C_1>0$ satisfying 
\begin{align*}
\intomega|\uep|^2\leq C_1\intomega|\nabla\uep|^2\quad\text{on }(0,\infty)
\end{align*}
and by the Cauchy-Schwarz-inequality we hence also obtain
\begin{align}\label{eq:uep-bound-eq2}
\Big(\intomega|\uep|\Big)^2\leq |\Omega|\intomega|\uep|^2\leq C_1|\Omega|\intomega|\nabla\uep|^2\quad\text{on }(0,\infty).
\end{align}
Denoting by $K=K(\Omega)>0$ the constant given in Lemma~\ref{lem:mos-trud} for $\eta=1$, we also introduce $$C_2=C_2(\Omega)=\frac{1}{\frac1{4\pi}+KC_1|\Omega|}$$ and employ Lemma~\ref{lem:mos-trud} with $\phi=\nep$, $\psi=|\uep|$, and $a=C_2(2\|\nabla\Phi\|_{\Lo[\infty]}\intomega n_0)^{-1}$
to find that by \eqref{eq:uep-bound-eq2} we have
\pagebreak
\begin{align*}
&\intomega\nep|\uep|\\
\leq\ &\frac{1}{a}\intomega\nep\ln\Big(\frac{\nep}{\overline{n}_0}\Big)+\frac{2a}{8\pi}\Big(\intomega n_0\Big)\intomega|\nabla\uep|^2+Ka\Big(\intomega n_0\Big)\Big(\intomega|\uep|\Big)^2+\frac{K}{a}\intomega n_0\\
\leq\ &\frac{1}{a}\intomega\nep\ln\Big(\frac{\nep}{\overline{n}_0}\Big)+\frac{a}{C_2}\Big(\intomega n_0\Big)\intomega|\nabla\uep|^2+\frac{K}{a}\intomega n_0\\
=\ &\frac{1}{a}\intomega\nep\ln\Big(\frac{\nep}{\overline{n}_0}\Big)+\frac{1}{2\|\nabla\Phi\|_{\Lo[\infty]}}\intomega|\nabla\uep|^2+\frac{K}{a}\intomega n_0
\end{align*}
on $(0,\infty)$ for all $\epsi\in(0,1)$. Plugging this back into \eqref{eq:uep-bound-eq1} shows that for all $\epsi\in(0,1)$
\begin{align*}
&\frac{1}{2}\frac{\intd}{\intd t}\intomega|\uep|^2+\frac{1}{2}\intomega|\nabla\uep|^2\\
\leq\ &\frac{\|\nabla\Phi\|_{\Lo[\infty]}}{a}\intomega\nep\ln\Big(\frac{\nep}{\overline{n}_0}\Big)+\frac{\|\nabla\Phi\|_{\Lo[\infty]} K}{a}\intomega n_0\\
\leq &\frac{2\|\nabla\Phi\|_{\Lo[\infty]}^2}{C_2}\Big(\intomega n_0\Big)\Big(\intomega\nep\ln(\nep+1)+|\ln(\overline{n}_0)|\intomega n_0\Big)+\frac{2\|\nabla\Phi\|_{\Lo[\infty]}^2 K}{C_2}\Big(\intomega n_0\Big)^2
\end{align*}
holds on $(0,\infty)$ and accordingly, we can find $C_3=C_3(n_0)>0$ such that
\begin{align}\label{eq:uep-bound-eq3}
\frac{\intd}{\intd t}\intomega|\uep|^2+\intomega|\nabla\uep|^2\leq C_3\intomega\nep\ln(\nep+1)+C_3\quad\text{on }(0,\infty)\text{ for all }\epsi\in(0,1).
\end{align}
In view of Lemma~\ref{lem:st-bound-neplnnep} and Lemma~\ref{lem:ode} this firstly provides $C_4=C_4(L,n_0,c_0,u_0,\gamma_0)>0$ satsifying 
\begin{align*}
\big\|\uep(\cdot,t)\big\|_{\Lo[2]}\leq C_4\quad\text{for all }\epsi\in(0,1)\text{ and }t>0
\end{align*}
and then, in a second step, by integrating \eqref{eq:uep-bound-eq3} also the bound on $\inttpone\intomega|\nabla\uep|^2$.
\end{bew}

\setcounter{equation}{0} 
\section{Assembling components of an energy-like functional}\label{sec3:energy}
In order to derive further $\epsi$-uniform bounds, we will follow arguments of \cite[Section 3]{winklerChemotaxisStokesInteractionVery2022}, where -- after introducing the primitive functions $D_{1,\epsi}(n)=\int_0^{n}\Dep(s)\intd s$, $D_{2,\epsi}(n)=\int_0^{n} D_{1,\epsi}(s)\intd s$ for $n\geq 0$ -- the time-evolution of 
$$y_\epsi(t)=\intomega D_{2,\epsi}(\nep)+b_1\intomega\frac{\nep|\nabla\cep|^2}{\cep}+b_2\intomega\frac{|\nabla\cep|^4}{\cep^3}+b_3\intomega\Psi_{2,\epsi}^{(s_0)}(\nep)$$
was considered for suitable $b_i\geq 0$, $i\in\{1,2,3\}$ and a certain $\Psi_{2,\epsi}^{(s_0)}$ (see \eqref{eq:Psi}).

Let us start by considering a differential inequality of $\intomega D_{2,\epsi}(\nep)$. The treatment of the right-hand side appearing therein will depend on the size of $\gamma$. As we will see later, restricting  $\gamma$ to $[0,\tfrac12]$ will provide us with the freedom to make use of the favorably-signed terms featured in the inequality for $\frac{\intd}{\intd t}\intomega\frac{\nep|\nabla\cep|^2}{\cep}$, to gain control on the term of the right-hand side of this lemma.

\begin{lemma}\label{lem:ener-part1-small-g}
Assume \eqref{eq:prop-S} that holds $\gamma\in[0,\tfrac12]$. For all $M>0$ there is $C=C(M)>0$ with the property: Whenever $\|c_0\|_{\Lo[\infty]}\leq M$, then for all $\epsi\in(0,1)$ the solution $(\nep,\cep,\uep)$ of \eqref{approxprob} satisfies
\begin{align*}
\frac{\intd}{\intd t}\intomega D_{2,\epsi}(\nep)+\frac{1}{2}\intomega\Dep^2(\nep)|\nabla\nep|^2\leq C\intomega\frac{\nep^2|\nabla\cep|^2}{\cep}
\end{align*}
holds on $(0,\infty)$.
\end{lemma}

\begin{bew}
Since $D_{2,\epsi}''(s)=\Dep(s)$, we find from integrating by parts that
\begin{align*}
\frac{\intd}{\intd t}\intomega D_{2,\epsi}(\nep)&=\intomega\nabla\cdot\big(\Dep(\nep)\nabla\nep\big) D_{2,\epsi}'(\nep)-\intomega\nabla\cdot\big(\nep\Sep(x,\nep,\cep)\cdot\nabla\cep\big)D_{2,\epsi}'(\nep)\\
&=-\intomega\Dep^2(\nep)|\nabla\nep|^2+\intomega\Dep(\nep)\nep\big(\Sep(x,\nep,\cep)\cdot\nabla\cep\big)\cdot\nabla\nep
\end{align*}
on $(0,\infty)$ for all $\epsi\in(0,1)$, where we used the divergence-free property of $\uep$ to remove the term containing the fluid-velocity. In light of Young's inequality, \eqref{eq:Sep}, \eqref{eq:prop-S}, the properties of $\rho_\epsi$ and $\chi_\epsi$, the non-decreasing property of $S_0$ and Lemma~\ref{lem:mass}, we conclude
\begin{align}\label{eq:ener-part1-eq1}
\frac{\intd}{\intd t}\intomega D_{2,\epsi}(\nep)+\frac{1}{2}\intomega\Dep^2(\nep)|\nabla\nep|^2\leq \frac{S_0^2(M+1)}{2}\intomega\frac{\nep^2|\nabla\cep|^2}{\cep^{2\gamma}}
\end{align}
on $(0,\infty)$ for all $\epsi\in(0,1)$. Here, the assumption $\gamma\in[0,\tfrac12]$ ensures $1-2\gamma\geq0$ and hence Lemma~\ref{lem:mass} and Young's inequality entail
\begin{align*}
\intomega\frac{\nep^2|\nabla\cep|^2}{\cep^{2\gamma}}\leq \|c_0\|_{\Lo[\infty]}^{1-2\gamma}\intomega\frac{\nep^2|\nabla\cep|^2}{\cep}\leq (M+1)\intomega\frac{\nep^2|\nabla\cep|^2}{\cep}
\end{align*}
on $(0,\infty)$ for all $\epsi\in(0,1)$, which upon combination with \eqref{eq:ener-part1-eq1} completes the proof.
\end{bew}

It is known from similar arguments within convex domains (\cite[Lemma 3.2]{winklerChemotaxisStokesInteractionVery2022}), that the an inequality for $\frac{\intd}{\intd t}\intomega\frac{\nep|\nabla\cep|^2}{\cep}$ can provide a favorably signed term controlling the mixed term on the right-hand side in Lemma~\ref{lem:ener-part1-small-g}. While, technically, this differential inequality can be established for any $\gamma\in[0,1]$, we can, however, only gain leverage from the left-hand side for small values of $\gamma$.
\begin{lemma}\label{lem:ener-part2}
Assume that \eqref{eq:prop-S} holds for some $\gamma\in[0,1]$. For all $M>0$ and $\eta>0$ there are $C_1=C_1(M,\eta)>0$ and $C_2>0$ with the property: Whenever $\|c_0\|_{\Lo[\infty]}\leq M$, then for all $\epsi\in(0,1)$
\begin{align}\label{eq:ener-part2}
&\frac{\intd}{\intd t}\intomega\frac{\nep|\nabla\cep|^2}{\cep}+\frac{1}{2}\intomega\frac{\nep^2|\nabla\cep|^2}{\cep}+2\intomega\frac{\nep|\nabla\cep|^4}{\cep^3}\nonumber\\
\leq\ &\eta\intomega\Dep^2(\nep)|\nabla\nep|^2+C_1\intomega|\nabla\nep|^2+C_1\intomega\frac{|D^2\cep|^2|\nabla\cep|^2}{\cep^3}+C_1\intomega\frac{|\nabla\cep|^6}{\cep^5}\nonumber\\
&\hspace*{1cm}+C_1\intomega|\nabla\uep|^3+C_1\intomega|\uep|^6+C_2\Big(\intomega n_0\Big)^2+C_2\quad\text{on }(0,\infty).
\end{align}
\end{lemma}

\begin{bew}
Straightforward computations using the first two equations of \eqref{approxprob} entail that
\begin{align*}
\frac{\intd}{\intd t}\intomega\frac{\nep|\nabla\cep|^2}{\cep}&=\intomega\frac{|\nabla\cep|^2}{\cep}\Big(\nabla\cdot\big(\Dep(\nep)\nabla\nep-\nep\Sep(x,\nep,\cep)\cdot\nabla\cep\big)-\uep\cdot\nabla\nep\Big)\\
&\quad+2\intomega\frac{\nep}{\cep}\nabla\cep\cdot\nabla\big(\Delta\cep-\nep\cep-\uep\nabla\cep\big)-\intomega\frac{\nep|\nabla\cep|^2}{\cep^2}\big(\Delta\cep-\nep\cep-\uep\nabla\cep\big)
\end{align*}
on $(0,\infty)$ for all $\epsi\in(0,1)$. Herein, we make use of the pointwise identity $\nabla\cep\cdot\nabla\Delta\cep=\frac{1}{2}\Delta|\nabla\cep|^2-|D^2\cep|^2$ to re-express the second term and thereby obtain
\begin{align*}
&\frac{\intd}{\intd t}\intomega\frac{\nep|\nabla\cep|^2}{\cep}\\
=\ &\intomega\frac{|\nabla\cep|^2}{\cep}\Big(\nabla\cdot\big(\Dep(\nep)\nabla\nep-\nep\Sep(x,\nep,\cep)\cdot\nabla\cep\big)-\uep\cdot\nabla\nep\Big)\\
&+\intomega\frac{\nep}{\cep}\Delta|\nabla\cep|^2-2\intomega\frac{\nep}{\cep}|D^2\cep|^2-2\intomega\frac{\nep^2}{\cep}|\nabla\cep|^2-2\intomega\nep(\nabla\cep\cdot\nabla\nep)\\
&\hspace*{0.5cm}-2\intomega\frac{\nep}{\cep}\nabla\cep\cdot(\nabla\uep\cdot\nabla\cep)-2\intomega\frac{\nep}{\cep}\nabla\cep\cdot(D^2\cep\cdot\uep)\\
&-\intomega\frac{\nep|\nabla\cep|^2}{\cep^2}\big(\Delta\cep-\nep\cep-\uep\nabla\cep\big)
\end{align*}
on $(0,\infty)$ for all $\epsi\in(0,1)$. Then, we find from multiple integrations by parts, additionally using $\nabla|\nabla\cep|^2=2 D^2\cep\cdot\nabla\cep$, that for all $\epsi\in(0,1)$ we have
\pagebreak
\begin{align}\label{eq:ener-part2_eq1}
&\frac{\intd}{\intd t}\intomega\frac{\nep|\nabla\cep|^2}{\cep}\nonumber\\
=\ &-2\intomega\frac{\Dep(\nep)}{\cep}\big(D^2\cep\cdot\nabla\cep\big)\cdot\nabla\nep+\intomega\frac{|\nabla\cep|^2}{\cep^2}\Dep(\nep)(\nabla\nep\cdot\nabla\cep)-\intomega\frac{|\nabla\cep|^2}{\cep}(\uep\cdot\nabla\nep)\nonumber\\
&\hspace*{0.5cm}+2\intomega\frac{\nep}{\cep}\big(D^2\cep\cdot\nabla\cep\big)\cdot\big(\Sep(x,\nep,\cep)\cdot\nabla\cep\big)-\intomega\frac{\nep}{\cep^2}|\nabla\cep|^2\nabla\cep\cdot\big(\Sep(x,\nep,\cep)\cdot\nabla\cep\big)
\nonumber\\
&-2\intomega\frac{1}{\cep}(D^2\cep\cdot\nabla\cep)\cdot\nabla\nep+2\intomega\frac{\nep}{\cep^2}\big((D^2\cep\cdot\nabla\cep)\cdot\nabla\cep\big)+\intromega\frac{\nep}{\cep}\frac{\partial|\nabla\cep|^2}{\partial\nu}-2\intomega\frac{\nep}{\cep}|D^2\cep|^2\nonumber\\
&\hspace*{0.5cm}-2\intomega\frac{\nep^2|\nabla\cep|^2}{\cep}-2\intomega\nep(\nabla\cep\cdot\nabla\nep)-2\intomega\frac{\nep}{\cep}\nabla\cep\cdot(\nabla\uep\cdot\nabla\cep)-2\intomega\frac{\nep}{\cep}\nabla\cep\cdot(D^2\cep\cdot\uep)\nonumber\\
&+\intomega\frac{|\nabla\cep|^2}{\cep^2}(\nabla\cep\cdot\nabla\nep)+2\intomega\frac{\nep}{\cep^2}\big((D^2\cep\cdot\nabla\cep)\cdot\nabla\cep\big)-2\intomega\frac{\nep|\nabla\cep|^4}{\cep^3}\nonumber\\
&\hspace*{0.5cm}+\intomega\frac{\nep^2|\nabla\cep|^2}{\cep}+\intomega\frac{\nep|\nabla\cep|^2}{\cep^2}(\uep\cdot\nabla\cep)=:I_{1,\epsi}+\dots+I_{18,\epsi}
\end{align}
on $(0,\infty)$. Note that $I_{9,\epsi}$, $I_{10,\epsi}$ and $I_{16,\epsi}$ in \eqref{eq:ener-part2_eq1} are already favorably signed terms and that $I_{17,\epsi}$ can be consumed by $I_{10,\epsi}$. To treat the remaining integrals, we will mostly rely on Young's inequality, one exception being the boundary integral $I_{8,\epsi}$, which requires additional work since we do not restrict this argument to convex domains. Applying Young's inequality to $I_{1,\epsi}$ and $I_{2,\epsi}$ entails that for all $\eta>0$ and $\epsi\in(0,1)$ we have
\begin{align*}
I_{1,\epsi}+I_{2,\epsi}\leq \eta\intomega\Dep^2(\nep)|\nabla\nep|^2+\frac{2M}{\eta}\intomega\frac{|D^2\cep|^2|\nabla\cep|^2}{\cep^3}+\frac{M}{2\eta}\intomega\frac{|\nabla\cep|^6}{\cep^5}
\end{align*}
on $(0,\infty)$, where we used Lemma~\ref{lem:mass} in both summands to increase the power of $\cep$ in the denominator by one at the cost of an additional factor $M$. Similarly, we find from two applications of Young's inequality in $I_{3,\epsi}$, that
\begin{align*}
I_{3,\epsi}\leq \intomega|\nabla\nep|^2+\frac14\intomega\frac{|\nabla\cep|^4}{\cep^2}|\uep|^2\leq\intomega|\nabla\nep|^2+\intomega|\uep|^6+M^2\intomega\frac{|\nabla\cep|^6}{\cep^5}
\end{align*}
on $(0,\infty)$ for all $\epsi\in(0,1)$. In $I_{4,\epsi}$ and $I_{5,\epsi}$ we estimate towards $I_{10}$ by employing Young's inequality, with some $\eta_1>0$ to be determined later, yielding that for all $\epsi\in(0,1)$
\begin{align*}
I_{4,\epsi}&\leq \eta_1\intomega\frac{\nep^2|\nabla\cep|^2}{\cep}+\frac{S_0^2(M+1)}{\eta_1}\intomega\frac{|D^2\cep|^2|\nabla\cep|^2}{\cep^{1+2\gamma}}\\&\leq\eta_1\intomega\frac{\nep^2|\nabla\cep|^2}{\cep}+\frac{M^{2(1-\gamma)}S_0^2(M+1)}{\eta_1}\intomega\frac{|D^2\cep|^2|\nabla\cep|^2}{\cep^{3}}\\
&\leq\eta_1\intomega\frac{\nep^2|\nabla\cep|^2}{\cep}+\frac{(M^2+1)S_0^2(M+1)}{\eta_1}\intomega\frac{|D^2\cep|^2|\nabla\cep|^2}{\cep^{3}}
\end{align*}
and
\begin{align*}
I_{5,\epsi}&\leq \eta_1\intomega\frac{\nep^2|\nabla\cep|^2}{\cep}+\frac{S_0^2(M+1)}{4\eta_1}\intomega\frac{|\nabla\cep|^6}{\cep^{3+2\gamma}}\leq \eta_1\intomega\frac{\nep^2|\nabla\cep|^2}{\cep}+\frac{(M^2+1)S_0^2(M+1)}{4\eta_1}\intomega\frac{|\nabla\cep|^6}{\cep^5}
\end{align*}
hold on $(0,\infty)$, where we also used \eqref{eq:Sep}, the boundedness of $\rho_\epsi$ and $\chi_\epsi$ together with \eqref{eq:prop-S}, the non-decreasing property of $S_0$ and $\gamma\leq 1$, as well as Lemma~\ref{lem:mass}. The same arguments also entail the estimates
\begin{align*}
I_{6,\epsi}&\leq \intomega|\nabla\nep|^2+M\intomega\frac{|D^2\cep|^2|\nabla\cep|^2}{\cep^3}, &I_{7,\epsi}&\leq\eta_1\intomega\frac{\nep^2|\nabla\cep|^2}{\cep}+\frac{1}{\eta_1}\intomega\frac{|D^2\cep|^2|\nabla\cep|^2}{\cep^3},\\
I_{11,\epsi}&\leq\eta_1\intomega\frac{\nep^2|\nabla\cep|^2}{\cep}+\frac{M}{\eta_1}\intomega|\nabla\nep|^2, &I_{12,\epsi}&\leq\eta_1\intomega\frac{\nep^2|\nabla\cep|^2}{\cep}+\frac{M^2}{\eta_1}\intomega\frac{|\nabla\cep|^6}{\cep^5}+\frac{1}{\eta_1}\intomega|\nabla\uep|^3,\\
I_{14,\epsi}&\leq \intomega|\nabla\nep|^2+\frac{M}{4}\intomega\frac{|\nabla\cep|^6}{\cep^5}, &I_{15,\epsi}&\leq \eta_1\intomega\frac{\nep^2|\nabla\cep|^2}{\cep}+\frac{1}{\eta_1}\intomega\frac{|D^2\cep|^2|\nabla\cep|^2}{\cep^3}
\end{align*}
on $(0,\infty)$ for all $\epsi\in(0,1)$. Analogously, multiple employments of Young's inequality to $I_{13,\epsi}$ and $I_{18,\epsi}$ show that also 
\begin{align*}
I_{13,\epsi}&\leq 2\intomega\frac{\nep}{\cep}|D^2\cep|^2+\frac12\intomega\frac{\nep|\nabla\cep|^2}{\cep}|\uep|^2\\
&\leq 2\intomega\frac{\nep}{\cep}|D^2\cep|^2+\eta_1\intomega\frac{\nep^2|\nabla\cep|^2}{\cep}+\frac{1}{16\eta_1}\intomega\frac{|\nabla\cep|^2}{\cep}|\uep|^4\\
&\leq 2\intomega\frac{\nep}{\cep}|D^2\cep|^2+\eta_1\intomega\frac{\nep^2|\nabla\cep|^2}{\cep}+\frac{1}{16\eta_1}\intomega|\uep|^6+\frac{M^2}{16\eta_1}\intomega\frac{|\nabla\cep|^6}{\cep^5}
\end{align*}
and
\begin{align*}
I_{18,\epsi}&\leq\eta_1\intomega\frac{\nep^2|\nabla\cep|^2}{\cep}+\frac{1}{4\eta_1}\intomega|\uep|^6+\frac{M^\nfrac12}{4\eta_1}\intomega\frac{|\nabla\cep|^6}{\cep^5}
\end{align*}
hold on $(0,\infty)$ for all $\epsi\in(0,1)$. We are only left with the treatment of the boundary integral $I_{8,\epsi}$. For this, we use arguments inspired by \cite[Proposition 3.2]{ISY-quasilin-pp-JDE14}. First, noting that there is $C_3=C_3(\Omega)>0$ such that whenever $w\in\CSp{2}{\bomega}$ satisfies $\frac{\partial w}{\partial\nu}=0$ on $\romega$, then $\frac{\partial|\nabla w|^2}{\partial\nu}\leq C_3|\nabla w|^2$ on $\romega$ (\cite[Lemma 4.2]{MS14}) and that by the trace inequality (see e.g. \cite[Theorem~1.6.6]{brennerMathematicalTheoryFinite1994} or \cite[Remark~52.9]{QS07}) for any $\delta>0$ there is $C_4=C_4(\delta,\Omega)>0$ such that for all $w\in\W[1,2]$
$\|w\|_{\LSp{2}{\romega}}\leq \delta\|\nabla w\|_{\Lo[2]}+C_4\|w\|_{\Lo[2]}$, we may estimate
\begin{align*}
I_{8,\epsi}&\leq C_3\intromega\frac{\nep|\nabla\cep|^2}{\cep}\\&\leq \|\nep\|_{\LSp{2}{\romega}}^2+C_3^2\Big\|\tfrac{|\nabla\cep|^2}{\cep}\Big\|_{\LSp{2}{\romega}}^2\\&
\leq\|\nabla\nep\|_{\Lo[2]}^2+C_4\|\nep\|_{\Lo[2]}^2+\Big\|\tfrac{\nabla|\nabla\cep|^2}{\cep}\Big\|_{\Lo[2]}^2+\Big\|\tfrac{|\nabla\cep|^2\nabla\cep}{\cep^2}\Big\|_{\Lo[2]}^2+C_4\Big\|\tfrac{|\nabla\cep|^2}{\cep}\Big\|_{\Lo[2]}^2
\end{align*}
on $(0,\infty)$ for all $\epsi\in(0,1)$. After which, drawing on a Poincaré-inequality, Lemma~\ref{lem:mass}, the identity $\nabla|\nabla\cep|^2=2D^2\cep\cdot\nabla\cep$ and Young's inequality, we find $C_5=C_5(\Omega)>0$ such that for all $\epsi\in(0,1)$
\begin{align*}
I_{8,\epsi}\leq C_5\intomega|\nabla\nep|^2+C_5\Big(\intomega n_0\Big)^2+4M\intomega\frac{|D^2\cep|^2|\nabla\cep|^2}{\cep^3}+\big(M+M^2)\intomega\frac{|\nabla\cep|^6}{\cep^5}+C_5
\end{align*}
holds on $(0,\infty)$. Gathering the estimates for $I_{k,\epsi}$, $k\in\{1,\dots,18\}$ and plugging them into \eqref{eq:ener-part2_eq1} we obtain that
\begin{align*}
&\frac{\intd}{\intd t}\intomega\frac{\nep|\nabla\cep|^2}{\cep}+(1-8\eta_1)\intomega\frac{\nep^2|\nabla\cep|^2}{\cep}+2\intomega\frac{\nep|\nabla\cep|^4}{\cep^3}\\
\leq\ &\eta\intomega\Dep^2(\nep)|\nabla\nep|^2+\Big(3+\frac{M}{\eta_1}+C_5\Big)\intomega|\nabla\nep|^2\\
&+\Big(\frac{2M}{\eta}+\frac{2}{\eta_1}+\frac{(M^2+1)S_0^2(M+1)}{\eta_1}+5M\Big)\!\intomega\frac{|D^2\cep|^2|\nabla\cep|^2}{\cep^3}\\
&\qquad+\Big(\frac{M}{2\eta}+\frac{5M^2+M^\nfrac12+(M^2+1)S_0^2(M+1)}{4\eta_1}+2M^2+\frac{5M}{4}\Big)\intomega\frac{|\nabla\cep|^6}{\cep^5}\\
&\qquad\qquad+\frac{1}{\eta_1}\intomega|\nabla\uep|^3+\Big(1+\frac{5}{16\eta_1}\Big)\intomega|\uep|^6+C_5\Big(\intomega n_0\Big)^2+C_5
\end{align*}
is valid on $(0,\infty)$ for all $\epsi\in(0,1)$. Choosing $\eta_1=\frac{1}{16}$ turns this into \eqref{eq:ener-part2} upon picking 
\begin{align*}
C_1=C_1(M,\eta):=\max\Big\{3+16M+C_5,\ &\tfrac{2M}{\eta}+16(M^2+1)S_0^2(M+1)+5M+32,\\& \tfrac{M}{2\eta}+4(M^2+1)S_0^2(M+1)+22M^2+\frac{5M}4+4M^\frac12\Big\}
\end{align*}
and $C_2:=C_5$.
\end{bew}

While, for small values of $\gamma$, the lemma above provides us with the necessary leverage to deal with $\intomega\frac{\nep^2|\nabla\cep|^2}{\cep}$ appearing in Lemma~\ref{lem:ener-part1-small-g}, it also introduces three new terms we currently lack control of. Most of these terms will be dealt with by considering the evolution of $\intomega\frac{|\nabla\cep|^4}{\cep^3}$. This is also a key element for the arguments in the setting with larger values of $\gamma$ discussed in Section~\ref{sec4:large-g}. In preparation we formulate the following lemma (cf. \cite[Lemma 3.2 and Lemma 3.4]{winklerApproachingLogarithmicSingularities2022}).

\begin{lemma}\label{lem:ln-ident}
Let $\phi\in\CSp{2}{\bomega}$ be such that $\phi>0$ in $\bomega$ and $\frac{\partial\phi}{\partial\nu}=0$ on $\romega$. Then
\begin{align}\label{eq:ln-ident}
|D^2\phi|^2=\phi^2|D^2\ln\phi|^2+\frac{1}{\phi}\nabla|\nabla\phi|^2\cdot\nabla\phi-\frac{1}{\phi^2}|\nabla\phi|^4\quad\text{in }\Omega
\end{align}
and
\begin{align}\label{eq:ln-ident-est1}
\intomega\frac{|\nabla\phi|^6}{\phi^5}\leq (4+\sqrt{2})^2\intomega\frac{|\nabla\phi|^2}{\phi}|D^2\ln\phi|^2
\end{align}
as well as
\begin{align}\label{eq:ln-ident-est2}
\intomega\frac{|D^2\phi|^2|\nabla\phi|^2}{\phi^3}\leq (5+\sqrt{2})^2\intomega\frac{|\nabla\phi|^2}{\phi}|D^2\ln\phi|^2.
\end{align}
\end{lemma}

These properties at hand, we follow the idea of \cite[Lemma 3.4]{winklerChemotaxisStokesInteractionVery2022} to obtain the following.

\begin{lemma}\label{lem:ener-part3}
There is $\delta>0$ such that for all $M>0$ one can find $C=C(M)>0$ with the property: Whenever $\|c_0\|_{\Lo[\infty]}\leq M$, for all $\epsi\in(0,1)$ the solution $(\nep,\cep,\uep)$ of \eqref{approxprob} satisfies
\begin{align*}
\frac{\intd}{\intd t}\intomega\frac{|\nabla\cep|^4}{\cep^3}&+\delta\intomega\frac{|D^2\cep|^2|\nabla\cep|^2}{\cep^3}+\delta\intomega\frac{|\nabla\cep|^6}{\cep^5}+\intomega\frac{\nep|\nabla\cep|^4}{\cep^3}\\
&\leq C\intomega|\nabla\nep|^2+C\intomega|\nabla\uep|^3+C\intomega|\uep|^6+C
\end{align*}
for all $t>0$.
\end{lemma}

\begin{bew}
In a quite similar approach to Lemma~\ref{lem:ener-part2}, we find from direct calculations involving integration by parts and the pointwise identity $\nabla\Delta\cep\cdot\nabla\cep=\frac{1}{2}\Delta|\nabla\cep|^2-|D^2\cep|^2$, that for all $\epsi\in(0,1)$
\begin{align}\label{eq:ener-part3-eq1}
&\frac{\intd}{\intd t}\intomega\frac{|\nabla\cep|^4}{\cep^3}\\=\ &4\intomega\frac{|\nabla\cep|^2}{\cep^3}\nabla\cep\cdot\nabla\big(\Delta\cep-\nep\cep-\uep\cdot\nabla\cep\big)-3\intomega\frac{|\nabla\cep|^4}{\cep^4}\big(\Delta\cep-\nep\cep-\uep\cdot\nabla\cep\big)\nonumber\\
=\ &-2\!\intomega\frac{\big|\nabla|\nabla\cep|^2\big|^2}{\cep^3}+6\!\intomega\frac{|\nabla\cep|^2}{\cep^4}\big(\nabla|\nabla\cep|^2\cdot\nabla\cep\big)+2\!\intromega\frac{|\nabla\cep|^2}{\cep^3}\frac{\partial|\nabla\cep|^2}{\partial\nu}-4\!\intomega\frac{|D^2\cep|^2|\nabla\cep|^2}{\cep^3}\nonumber\\
&-\!\intomega\frac{\nep|\nabla\cep|^4}{\cep^3}-4\!\intomega\frac{|\nabla\cep|^2}{\cep^2}(\nabla\cep\cdot\nabla\nep)-4\!\intomega\frac{|\nabla\cep|^2}{\cep^3}\big(\nabla\uep\cdot\nabla\cep)\cdot\nabla\cep-4\!\intomega\frac{|\nabla\cep|^2}{\cep^3}\big(D^2\cep\cdot\uep)\cdot\nabla\cep\nonumber\\
&+6\!\intomega\frac{|\nabla\cep|^2}{\cep^4}\big(\nabla|\nabla\cep|^2\cdot\nabla\cep\big)-12\!\intomega\frac{|\nabla\cep|^6}{\cep^5}+3\!\intomega\frac{|\nabla\cep|^4}{\cep^4}(\uep\cdot\nabla\cep)=:J_{1,\epsi}+\dots+J_{11,\epsi}
\end{align}
holds on $(0,\infty)$. Then, employing the identity \eqref{eq:ln-ident} from Lemma~\ref{lem:ln-ident} for $\phi=\cep$, we see that
\begin{align*}
&J_{1,\epsi}+J_{4,\epsi}+J_{2,\epsi}+J_{9,\epsi}+J_{10,\epsi}\\
=\ &-2\intomega\frac{\big|\nabla|\nabla\cep|^2\big|^2}{\cep^3}-4\intomega\frac{|\nabla\cep|^2}{\cep^3}\Big(\cep^2|D^2\ln\cep|^2+\frac{1}{\cep}\nabla|\nabla\cep|^2\cdot\nabla\cep-\frac{1}{\cep^2}|\nabla\cep|^4\Big)\\&+12\intomega\frac{|\nabla\cep|^2}{\cep^4}\big(\nabla|\nabla\cep|^2\cdot\nabla\cep\big)-12\intomega\frac{|\nabla\cep|^6}{\cep^5}\\
=\ &-4\intomega\frac{|\nabla\cep|^2}{\cep}|D^2\ln\cep|^2-2\intomega\frac{1}{\cep^3}\bigg|\nabla|\nabla\cep|^2-\frac{2|\nabla\cep|^2\nabla\cep}{\cep}\bigg|^2\leq-4\intomega\frac{|\nabla\cep|^2}{\cep}|D^2\ln\cep|^2
\end{align*}
on $(0,\infty)$ for all $\epsi\in(0,1)$, so that in \eqref{eq:ener-part3-eq1} we can accordingly estimate
\begin{align}\label{eq:ener-part3-eq3}
\frac{\intd}{\intd t}\intomega\frac{|\nabla\cep|^4}{\cep^3}&\leq -4\intomega\frac{|\nabla\cep|^2}{\cep}|D^2\ln\cep|^2+2\intromega\frac{|\nabla\cep|^2}{\cep^3}\frac{\partial|\nabla\cep|^2}{\partial\nu}-\intomega\frac{\nep|\nabla\cep|^4}{\cep^3}-4\intomega\frac{|\nabla\cep|^2}{\cep^2}(\nabla\cep\cdot\nabla\nep)\nonumber\\
&\quad-4\intomega\frac{|\nabla\cep|^2}{\cep^3}(\nabla\uep\cdot\nabla\cep)\cdot\nabla\cep-4\intomega\frac{|\nabla\cep|^2}{\cep^3}(D^2\cep\cdot\uep)\cdot\nabla\cep+3\intomega\frac{|\nabla\cep|^4}{\cep^4}(\uep\cdot\nabla\cep)\nonumber\\
&:=\tilde{J}_{1,\epsi}+\dots+\tilde{J}_{7,\epsi}
\end{align}
on $(0,\infty)$ for all $\epsi\in(0,1)$. Introducing $\delta:=\min\big\{\frac{1}{(5+\sqrt{2})^2},\frac{1}{(4+\sqrt{2})^2}\big\}=\frac{1}{(5+\sqrt{2})^2}$, we obtain from \eqref{eq:ln-ident-est1} and \eqref{eq:ln-ident-est2} of Lemma~\ref{lem:ln-ident} that
\begin{align}\label{eq:ener-part3-eq4}
\tilde{J}_{1,\epsi}&\leq -\frac{2}{(5+\sqrt{2})^2}\intomega\frac{|D^2\cep|^2|\nabla\cep|^2}{\cep^3}-\frac{2}{(4+\sqrt{2})^2}\intomega\frac{|\nabla\cep|^6}{\cep^5}\nonumber\\
&\leq -2\delta\intomega\frac{|D^2\cep|^2|\nabla\cep|^2}{\cep^3}-2\delta\intomega\frac{|\nabla\cep|^6}{\cep^5}\quad\text{ on }(0,\infty)
\end{align}
holds for all $\epsi\in(0,1)$, whereas Young's inequality employed in $\tilde{J}_{4,\epsi},\tilde{J}_{5,\epsi},\tilde{J}_{6,\epsi}$ and $\tilde{J}_{7,\epsi}$ entails that
\begin{align}\label{eq:ener-part3-eq5}
\tilde{J}_{4,\epsi}&\leq \frac{\delta}{5}\intomega\frac{|\nabla\cep|^6}{\cep^5}+\frac{20}{\delta}\intomega\cep|\nabla\nep|^2\leq \frac{\delta}{5}\intomega\frac{|\nabla\cep|^6}{\cep^5}+\frac{20M}{\delta}\intomega|\nabla\nep|^2,\\\label{eq:ener-part3-eq6}
\tilde{J}_{5,\epsi}&\leq 4\intomega\Big(\frac{\delta}{20}\frac{|\nabla\cep|^6}{\cep^5}\Big)^\nfrac23\Big(\frac{20^2}{\delta^2}\cep|\nabla\uep|^3\Big)^\nfrac13\leq\frac{\delta}{5}\intomega\frac{|\nabla\cep|^6}{\cep^5}+\frac{4\cdot20^2 M}{\delta^2}\intomega|\nabla\uep|^3,
\end{align}
and
\begin{align}\label{eq:ener-part3-eq7}
\tilde{J}_{6,\epsi}&\leq \frac{4\delta}{5}\intomega\frac{|D^2\cep|^2|\nabla\cep|^2}{\cep^3}+\frac{5}{\delta}\intomega\frac{|\nabla\cep|^4}{\cep^3}|\uep|^2\nonumber\\
&=\frac{4\delta}{5}\intomega\frac{|D^2\cep|^2|\nabla\cep|^2}{\cep^3}+\frac{5}{\delta}\intomega\Big(\frac{\delta^2}{5^2}\frac{|\nabla\cep|^6}{\cep^5}\Big)^\nfrac23\Big(\frac{5^4}{\delta^4}\cep|\uep|^6\Big)^\nfrac13\nonumber\\
&\leq \frac{4\delta}{5}\intomega\frac{|D^2\cep|^2|\nabla\cep|^2}{\cep^3}+\frac{\delta}{5}\intomega\frac{|\nabla\cep|^6}{\cep^5}+\frac{5^5M}{\delta^5}\intomega|\uep|^6 ,
\end{align}
as well as
\begin{align}\label{eq:ener-part3-eq8}
\tilde{J}_{7,\epsi}&\leq 3\intomega\Big(\frac{\delta}{15}\frac{|\nabla\cep|^6}{\cep^5}\Big)^\nfrac56\Big(\frac{15^5}{\delta^5}\cep|\uep|^6\Big)^\nfrac16\leq \frac{\delta}{5}\intomega\frac{|\nabla\cep|^6}{\cep^5}+\frac{3\cdot15^5 M}{\delta^5}\intomega|\uep|^6
\end{align}
on $(0,\infty)$ for all $\epsi\in(0,1)$. We are left with estimating the boundary integral $\tilde{J}_{2,\epsi}$. Similarly to the step in Lemma~\ref{lem:ener-part2}, we rely on the fact that there is $C_1=C_1(\Omega)>0$ such that $\frac{\partial|\nabla w|^2}{\partial\nu}\leq C_1|\nabla w|^2$ holds on $\romega$ for all $w\in\CSp{2}{\bomega}$ satisfying $\frac{\partial w}{\partial \nu}=0$ on $\romega$ (cf. \cite[Lemma~4.2]{ISY-quasilin-pp-JDE14}) and that by the trace inequality (\cite[Theorem~1.6.6]{brennerMathematicalTheoryFinite1994}, \cite[Remark~52.9]{QS07}) there is $C_2=C_2(\Omega)>0$ satisfying $2C_1\|w\|_{\LSp{2}{\romega}}^2\leq\frac{\delta}{20}\|\nabla w\|_{\Lo[2]}^2+C_2\|w\|_{\Lo[2]}^2$ and estimate
\begin{align*}
\tilde{J}_{2,\epsi}&\leq 2C_1\intromega\frac{|\nabla\cep|^4}{\cep^3}=2C_1\Big\|\frac{|\nabla\cep|^2}{\cep^\nfrac32}\Big\|_{\LSp{2}{\romega}}^2\leq \frac{\delta}{20}\Big(\intomega\frac{\big|\nabla|\nabla\cep|^2\big|^2}{\cep^3}+\frac{9}{4}\intomega\frac{|\nabla\cep|^6}{\cep^5}\Big)+C_2\intomega\frac{|\nabla\cep|^4}{\cep^3}
\end{align*}
on $(0,\infty)$ for all $\epsi\in(0,1)$. Drawing once more on the pointwise identity $\nabla|\nabla\cep|^2=2D^2\cep\nabla\cep$ and employing Young's inequality in the last term, we thus find $C_3=C_3(\Omega)>0$ such that
\begin{align}\label{eq:ener-part3-eq9}
\tilde{J}_{2,\epsi}\leq \frac{\delta}{5}\intomega\frac{|D^2\cep|^2|\nabla\cep|^2}{\cep^3}+\frac{\delta}{5}\intomega\frac{|\nabla\cep|^6}{\cep^5}+C_3 M\quad\text{on }(0,\infty)\text{ for all }\epsi\in(0,1).
\end{align}
Gathering \eqref{eq:ener-part3-eq3}--\eqref{eq:ener-part3-eq9} and letting $C_4=C_4(M):=\max\big\{\frac{20M}{\delta},\frac{4\cdot20^2 M}{\delta^2},\frac{(5^5+3\cdot15^5)M}{\delta^5},C_3 M\big\}>0$ yields
\begin{align*}
\frac{\intd}{\intd t}\intomega\frac{|\nabla\cep|^4}{\cep^3}&+\delta\intomega\frac{|D^2\cep|^2|\nabla\cep|^2}{\cep^3}+\delta\intomega\frac{|\nabla\cep|^6}{\cep^5}+\intomega\frac{\nep|\nabla\cep|^4}{\cep^3}\\
&\leq C_4\intomega|\nabla\nep|^2+C_4\intomega|\nabla\uep|^3+C_4\intomega|\uep|^6+C_4\quad\text{on }(0,\infty)\text{ for all }\epsi\in(0,1),
\end{align*}
completing the proof.
\end{bew}

Since the diffusion features a possible degeneracy at population density zero, we can only rely on the dissipative term $\intomega \Dep^2(\nep)|\nabla\nep|^2$ to control the troublesome term $\intomega|\nabla\nep|^2$ -- which is present in the inequalities of both Lemma~\ref{lem:ener-part2} and Lemma~\ref{lem:ener-part3} -- in the region $\{\nep>s_0\}$.
Hence, we need to augment the approach by an additional testing procedure, which captures the behavior of the first equation for small densities in a more detailed fashion. To this end, it was considered in \cite{winklerChemotaxisStokesInteractionVery2022} to investigate a suitably truncated version of the reciprocal diffusion, i.e. introduce

\begin{align}\label{eq:Psi}
\Psi_{0,\epsi}^{(s_0)}(s):=\begin{cases}
\frac{1}{\Dep(s)}\quad&\text{if }s\in(0,s_0),\\
\frac{2s_0-s}{s_0\Dep(s_0)}&\text{if }s\in[s_0,2s_0],\\
0&\text{if }s>2s_0.
\end{cases}
\end{align}
and the primitives
$\Psi_{1,\epsi}^{(s_0)}(s):=-\int_s^{2s_0}\Psi_{0,\epsi}^{(s_0)}(\sigma)\intd\sigma$ and $\Psi_{2,\epsi}^{(s_0)}(s):=-\int_s^{2s_0}\Psi_{1,\epsi}^{(s_0)}(\sigma)\intd\sigma$, and then consider the time-evolution of the second primitive. This approach provides us with an additional well-signed term for the gradient of $\nep$ in subdomains with small density. Since our first equation does not differ significantly from the system in the referenced work, the lemma below is more or less a verbatim copy of \cite[Lemma~3.5]{winklerChemotaxisStokesInteractionVery2022}. Note that the condition $\liminf_{n\searrow0}\frac{D(n)}{n}>0$ is crucial here to ensure that $\kappa>0$.

\begin{lemma}\label{lem:ener-part4}
For $s_0>0$ set $$\kappa=\kappa(s_0):=\inf_{n\in(0,2s_0)}\frac{D(n)}{n}>0.$$
Then $\Psi_{2,\epsi}^{(s_0)}\in\CSp{2}{(0,\infty)}$ with $$0\leq \Psi_{2,\epsi}^{(s_0)}(s)\leq \frac{3s_0}{\kappa}\quad\text{for all }s>0.$$ Moreover, for all $M>0$ and $\eta>0$ there is $C=C(M,\eta)>0$ with the property: Whenever $\gamma\in[0,\tfrac56]$ and $\|c_0\|_{\Lo[\infty]}\leq M$, then
\begin{align}\label{eq:ener-part4}
\frac{\intd}{\intd t}\intomega\Psi_{2,\epsi}^{(s_0)}(\nep)+\frac{1}{2}\int_{\{\nep<s_0\}}|\nabla\nep|^2\leq \eta\intomega\frac{|\nabla\cep|^6}{\cep^{5}}+\frac{C}{\kappa^3}
\end{align}
holds on $(0,\infty)$ for all $\epsi\in(0,1)$.
\end{lemma}

\begin{bew}
For the results concerning the regularity of $\Psi_{2,\epsi}^{(s_0)}$ and its boundedness, we refer the reader to \cite[Lemma~3.5]{winklerChemotaxisStokesInteractionVery2022}. Here,  for the sake of completeness regarding the quasi-energy functional. we only reproduce the necessary steps to verify the differential inequality.

Integrating by parts and drawing on Young's inequality as well as the nonnegativity of $\Dep$ and $\Psi_{0,\epsi}^{(s_0)}$ shows that
\begin{align}\label{eq:ener-part4-eq1}
\frac{\intd}{\intd t}\intomega\Psi_{2,\epsi}^{(s_0)}(\nep)\leq\ &-\frac{1}{2}\intomega\Psi_{0,\epsi}^{(s_0)}(\nep)\Dep(\nep)|\nabla\nep|^2+\frac{S_0^2(M+1)}{2}\intomega\Psi_{0,\epsi}^{(s_0)}(\nep)\frac{\nep^2|\nabla\cep|^2}{\Dep(\nep)\cep^{2\gamma}}\nonumber\\
\leq\ &-\frac{1}{2}\int\limits_{\{\nep<s_0\}}\!\!\Psi_{0,\epsi}^{(s_0)}(\nep)\Dep(\nep)|\nabla\nep|^2+\frac{S_0^2(M+1)}{2}\intomega\Psi_{0,\epsi}^{(s_0)}(\nep)\frac{\nep^2|\nabla\cep|^2}{\Dep(\nep)\cep^{2\gamma}}
\end{align}
on $(0,\infty)$ for all $\epsi\in(0,1)$. Making use of the explicit form of $\Psi_{0,\epsi}^{(s_0)}$ and the definition of $\kappa$ to estimate
\begin{align*}
\Psi_{0,\epsi}^{(s_0)}(s)\frac{s^2}{\Dep(s)}= \frac{(2s_0-s)s^2}{s_0\Dep(s_0)\Dep(s)}\leq \frac{2s_0}{\Dep(s_0)}\cdot\frac{s}{\Dep(s)}\leq \frac{2}{\kappa^2}\quad\text{for all }s\in[s_0,2s_0],
\end{align*}
and
\begin{align*}
\Psi_{0,\epsi}^{(s_0)}(s)\frac{s^2}{\Dep(s)}=\frac{s^2}{\Dep^2(s)}\leq \frac{1}{\kappa^2}\quad\text{for all }s\in(0,s_0),
\end{align*}
we can split the domain of integration in the last term of \eqref{eq:ener-part4-eq1} to obtain
\begin{align*}
\intomega\Psi_{0,\epsi}^{(s_0)}(\nep)\frac{\nep^2|\nabla\cep|^2}{\Dep(\nep)\cep^{2\gamma}}
\leq\ &\frac{1}{\kappa^2}\!\!\int\limits_{\{\nep<s_0\}}\!\!\frac{|\nabla\cep|^2}{\cep^{2\gamma}} +\frac{2}{\kappa^2}\!\!\!\int\limits_{\{\nep\in[s_0,2s_0]\}}\!\!\!\frac{|\nabla\cep|^2}{\cep^{2\gamma}}\leq\frac{2}{\kappa^2}\intomega\frac{|\nabla\cep|^2}{\cep^{2\gamma}}
\end{align*}
on $(0,\infty)$ for all $\epsi\in(0,1)$. Using $\gamma\leq \frac56$, a final employment of Young's inequality then entails the existence of $C_1=C_1(\eta)>0$ such that
\begin{align*}
&\frac{\intd}{\intd t}\intomega\Psi_{2,\epsi}^{(s_0)}(\nep)+\frac{1}{2}\!\int\limits_{\{\nep<s_0\}}\!\!\Psi_{0,\epsi}^{(s_0)}(\nep)\Dep(\nep)|\nabla\nep|^2\leq \eta\intomega\frac{|\nabla\cep|^6}{\cep^5}+\frac{C_1 S_0^3(M+1)(M^3+1)|\Omega|}{\kappa^3}
\end{align*}
on $(0,\infty)$ for all $\epsi\in(0,1)$. Plugging in the explicit form of $\Psi_{0,\epsi}^{(s_0)}$ yields \eqref{eq:ener-part4} as claimed.
\end{bew}

Amalgamating suitable multiples of the inequalities of Lemmas~\ref{lem:ener-part1-small-g} and \ref{lem:ener-part2} as well as Lemmas~\ref{lem:ener-part3} and \ref{lem:ener-part4}, we are now able to identify some $L=L(M)$ for which the following energy-like structure emerges.

\begin{lemma}\label{lem:ener-part5-small-g}
Assume that \eqref{eq:prop-S} holds with some $\gamma\in[0,\tfrac12]$. Let $\delta>0$ be given by Lemma~\ref{lem:ener-part3}. For all $M>0$ one can fix $L=L(M)>0$, $b_1=b_1(M)>0$, $b_2=b_2(M)>0$, $b_3=b_3(M)>0$ such that: Whenever $D$ satisfies \eqref{eq:liminf} and $\|c_0\|_{\Lo[\infty]}\leq M$, there are $s_0=s_0(L)>0$, $\Gamma_1=\Gamma_1(M)>0$ and $\Gamma_2=\Gamma_2(M,n_0)>0$ such that for all $\epsi\in(0,1)$
\begin{align}\label{eq:ener-F-def}
\Fep(t):=\intomega D_{2,\epsi}\big(\nep(\cdot,t)\big)+b_1\intomega\frac{\nep(\cdot,t)\big|\nabla\cep(\cdot,t)\big|^2}{\cep(\cdot,t)}+b_2\intomega\frac{\big|\nabla\cep(\cdot,t)\big|^4}{\cep^3(\cdot,t)}+b_3\intomega\Psi_{2,\epsi}^{(s_0)}\big(\nep(\cdot,t)\big)
\end{align}
satisfies
\begin{align}\label{eq:ener-part5-eq0}
\frac{\intd}{\intd t}\Fep(t)&+\frac{1}{8}\intomega\Dep^2(\nep)|\nabla\nep|^2+\frac{b_1}{4}\intomega\frac{\nep^2|\nabla\cep|^2}{\cep}+\frac{b_2\delta}{4}\intomega\frac{|\nabla\cep|^6}{\cep^5}+\intomega|\nabla\nep|^2+\frac{b_1}{4}\intomega\frac{|\nabla\cep|^2}{\cep}\nonumber\\
&\leq \Gamma_1\Big(\intomega|\nabla\uep|^3+\intomega|\uep|^6\Big) +\Gamma_2
\end{align}
for all $t>0$.
\end{lemma}

\begin{bew}
Let $C_1=C_1(M)>0$ be the constant obtained in Lemma~\ref{lem:ener-part1-small-g}. Pick $b_1=b_1(M)=4C_1$, $\eta_1=\eta_1(M)=\frac{1}{4b_1}$ and then denote by $C_2=C_2(M,\eta_1)>0$, $C_3>0$ the constants provided by Lemma~\ref{lem:ener-part2} and by $C_4=C_4(M)>0$ the constant obtained in Lemma~\ref{lem:ener-part3}. Then, choose $b_2=b_2(M)=\frac{2b_1 C_2}{\delta}$ and $b_3=b_3(M)=2C_2b_1+2C_4b_2+2$, let $\eta_2=\eta_2(M)=\frac{b_2\delta}{8b_3}$ and denote by $C_5=C_5(M,\eta_2)>0$ the corresponding constant provided by Lemma~\ref{lem:ener-part4}. Finally, set $L=L(M)=\sqrt{8(C_2 b_1+C_4 b_2+1)}$ and let $s_0=s_0(L)$ be given in accordance with \eqref{eq:s0}, i.e. such that for all $s\geq s_0$ we have $D(s)\geq L$.

With these coefficients fixed, we once more set $\kappa=\kappa(s_0):=\inf_{n\in(0,2s_0)}\frac{D(n)}{n}$ and then find from a combination of Lemma~\ref{lem:ener-part1-small-g}, Lemma~\ref{lem:ener-part2}, Lemma~\ref{lem:ener-part3} and Lemma~\ref{lem:ener-part4} that
\begin{align}\label{eq:ener-part5-eq1}
&\frac{\intd}{\intd t}\Fep(t)+\frac{1}{2}\intomega\Dep^2(\nep)|\nabla\nep|^2+\frac{b_1}{2}\intomega\frac{\nep^2|\nabla\cep|^2}{\cep}+2b_1\intomega\frac{\nep|\nabla\cep|^4}{\cep^3}\nonumber\\
&\hspace*{1.2cm}+b_2\delta\intomega\frac{|D^2\cep|^2|\nabla\cep|^2}{\cep^3}+b_2\delta\intomega\frac{|\nabla\cep|^6}{\cep^5}+b_2\intomega\frac{\nep|\nabla\cep|^4}{\cep^3}+\frac{b_3}{2}\!\!\int\limits_{\{\nep<s_0\}}\!\!|\nabla\nep|^2\nonumber\\
\leq\ &C_1\intomega\frac{\nep^2|\nabla\cep|^2}{\cep}+b_1\eta_1\intomega\Dep^2(\nep)|\nabla\nep|^2+C_2 b_1\intomega|\nabla\nep|^2+C_2 b_1\intomega\frac{|D^2\cep|^2|\nabla\cep|^2}{\cep^3}\\
&\hspace*{0.3cm}+C_2 b_1\intomega\frac{|\nabla\cep|^6}{\cep^5}+C_2 b_1\intomega|\nabla\uep|^3+ C_2 b_1\intomega|\uep|^6+C_4 b_2\intomega|\nabla\nep|^2+C_4 b_2\intomega|\nabla\uep|^3\nonumber\\
&\hspace*{0.6cm}+C_4 b_2\intomega|\uep|^6+\eta_2 b_3\intomega\frac{|\nabla\cep|^6}{\cep^5}+C_3b_1\Big(\intomega n_0\Big)^2+C_3b_1+C_4b_2+\frac{C_5b_3}{\kappa^3}\nonumber
\end{align}
on $(0,\infty)$ for all $\epsi\in(0,1)$. Since $\gamma\leq\frac12$, we may employ Lemma~\ref{lem:mass} and Young's inequality to find that
\begin{align*}
\frac{b_1}{4}\intomega\frac{|\nabla\cep|^2}{\cep}&=\intomega\Big(\frac{b_2\delta}{8}\frac{|\nabla\cep|^6}{\cep^5}\Big)^\nfrac13\cdot\Big(\Big(\frac{8}{b_2\delta}\Big)^\nfrac13 \frac{b_1}{4}\cep^{\frac23}\Big)\leq \frac{b_2\delta}{8}\intomega\frac{|\nabla\cep|^6}{\cep^5}+\Big(\frac{8}{b_2\delta}\Big)^\nfrac12\Big(\frac{b_1}{4}\Big)^\nfrac32|\Omega|M
\end{align*}
on $(0,\infty)$ for all $\epsi\in(0,1)$. Hence, adding $\frac{b_1}{4}\intomega\frac{|\nabla\cep|^2}{\cep}$ to both sides of \eqref{eq:ener-part5-eq1}, we obtain from rearranging, using the specific choices for $b_1,b_2,b_3$ and $\eta$ as well as putting $\Gamma_1=\Gamma_1(M):=C_2 b_1+C_4 b_2$ and $\Gamma_2=\Gamma_2(M,n_0):=\big(\frac{8}{b_2\delta}\big)^\nfrac12\big(\frac{b_1}{4}\big)^\nfrac32|\Omega|M+C_3b_1(\intomega n_0)^2+C_3b_1+C_4b_2+\frac{C_5b_3}{\kappa^3}$ that
\begin{align}\label{eq:ener-part5-eq2}
&\frac{\intd}{\intd t}\Fep(t)+\frac{1}{4}\intomega\Dep^2(\nep)|\nabla\nep|^2+\frac{b_1}{4}\intomega\frac{\nep^2|\nabla\cep|^2}{\cep}+\big(2b_1+b_2)\intomega\frac{\nep|\nabla\cep|^4}{\cep^3}\nonumber\\
&\hspace*{1.2cm}+\frac{b_2\delta}{2}\intomega\frac{|D^2\cep|^2|\nabla\cep|^2}{\cep^3}+\frac{b_2\delta}{4}\intomega\frac{|\nabla\cep|^6}{\cep^5}+\intomega|\nabla\nep|^2+\frac{b_1}{4}\intomega\frac{|\nabla\cep|^2}{\cep}\\
\leq\ &(\Gamma_1+1)\intomega|\nabla\nep|^2-(\Gamma_1+1)\!\!\int\limits_{\{\nep<s_0\}}\!\!|\nabla\nep|^2+\Gamma_1\Big(\intomega|\nabla\uep|^3+\intomega|\uep|^6\Big)+\Gamma_2\nonumber
\end{align}
on $(0,\infty)$ for all $\epsi\in(0,1)$. Then, recalling that by the choice for $s_0>0$ and $L=\sqrt{8(\Gamma_1+1)}$ we have $\Dep^2(s)\geq 8(\Gamma_1+1)$ for all $s\geq s_0$, we find that
\begin{align*}
(\Gamma_1+1)\Big(\intomega|\nabla\nep|^2-\!\!\int\limits_{\{\nep<s_0\}}\!\!|\nabla\nep|^2\Big)=(\Gamma_1+1)\!\!\int\limits_{\{\nep\geq s_0\}}\!\!\frac{\Dep^2(\nep)}{\Dep^2(\nep)}|\nabla\nep|^2\leq \frac{1}{8}\intomega\Dep^2(\nep)|\nabla\nep|^2
\end{align*}
on $(0,\infty)$ for all $\epsi\in(0,1)$. Accordingly, we arrive at \eqref{eq:ener-part5-eq0} from \eqref{eq:ener-part5-eq2} after dropping two nonnegative terms on the left.
\end{bew}

To derive bounds from this inequality, our next aim is to establish a linear dissipation term for $\Fep$. In a first step we will estimate $\intomega\Dep^2(\nep)|\nabla\nep|^2$ from below by $\intomega D_{2,\epsi}(\nep)$. In \cite{winklerChemotaxisStokesInteractionVery2022} this was achieved by relying on a Poincaré-type inequality requiring convexity of the domain (e.g. \cite[Corollary 9.1.4]{Jost-PDE07}). In order to get rid of this requirement we will instead draw on the following lemma taken from \cite[Lemma 9.1]{LanWin2017}, which does not have this limitation.
\begin{lemma}\label{lem:PC-Lan-Win}
Let $\Omega\subset\R^\dimN$ be a bounded domain with smooth boundary, and let $\varpi>0$ and $p\in[1,\infty)$. Then there exists $C(\Omega,\varpi,p)>0$ with the property that for all $\varphi\in\W[1,p]$,
\begin{align*}
\bigg(\intomega\Big|\varphi-\frac1{|B|}\int_B\varphi\Big|^p\bigg)^\frac{1}{p}\leq C(\Omega,\varpi,p)\Big(\intomega|\nabla\varphi|^p\Big)^\frac{1}{p}
\end{align*}
holds for any measurable set $B\subset\Omega$ with $|B|\geq\varpi$.
\end{lemma}
We briefly remark that the referenced work states the lemma only for subsets satisfying $|B|=\varpi$. Nevertheless, since $|B|$ enters the constant with a negative exponent this is easily extended to the form above. With this we can now prove the following. (Compare \cite[Lemma 3.6]{winklerChemotaxisStokesInteractionVery2022}.)

\begin{lemma}\label{lem:ener-part6}
Assume that \eqref{eq:liminf} holds with some $L>0$ and that $s_0=s_0(L)>0$ is as in \eqref{eq:s0}. Then there is $C=C(L)>0$ such that for all $\epsi\in(0,1)$
\begin{align*}
\intomega D_{2,\epsi}\big(\nep(\cdot,t)\big)\leq C\intomega\Dep^2\big(\nep(\cdot,t)\big)\big|\nabla\nep(\cdot,t)\big|^2+C\big(D_1(2\overline{n}_0)+4\overline{n}_0\big)^2+ C\quad\text{for all }t>0,
\end{align*}
where $D_1(n):=\int_0^{n}D(s)\intd s$ and $\overline{n}_0:=\frac{1}{|\Omega|}\intomega n_0$.
\end{lemma}

\begin{bew}
We introduce the number
\begin{align*}
\delta_\epsi=\delta_\epsi(\overline{n}_0):=D_{1,\epsi}\big(2\overline{n}_0\big)=\int_0^{2\overline{n}_0}\Dep(s)\intd s>0.
\end{align*}
Since Lemma~\ref{lem:mass} and the Chebyshev inequality imply that 
\begin{align*}
\intomega n_0=\intomega \nep\geq \Big(\frac{2}{|\Omega|}\intomega n_0\Big) \Big|\big\{\nep\geq\frac{2}{|\Omega|}\intomega n_0\big\}\Big|\quad\text{on }(0,\infty)\text{ for all }\epsi\in(0,1),
\end{align*}
we find from the fact that $D_{1,\epsi}'(s)=\Dep(s)>0$ for all $s>0$ and $\epsi\in(0,1)$ that
\begin{align*}
\Big|\big\{D_{1,\epsi}(\nep)<\delta_\epsi\big\}\Big|=\Big|\big\{\nep<\frac{2}{|\Omega|}\intomega n_0\big\}\Big|=|\Omega|-\Big|\big\{\nep\geq\frac{2}{|\Omega|}\intomega n_0\big\}\Big|\geq \frac{|\Omega|}{2}
\end{align*}
is valid on $(0,\infty)$ for all $\epsi\in(0,1)$. Accordingly, an application of Lemma~\ref{lem:PC-Lan-Win} for $\varphi=\big(D_{1,\epsi}(\nep)-\delta_\epsi\big)_+$ and $B=\big\{D_{1,\epsi}(\nep)<\delta_\epsi\big\}$ entails the existence of $C_1=C_1(\Omega)>0$ satisfying
\begin{align}\label{eq:ener-part-6-eq0}
\intomega\big(D_{1,\epsi}(\nep)-\delta_\epsi\big)_+^2\leq C_1\intomega \big|\nabla\big(D_{1,\epsi}(\nep)-\delta_\epsi\big)_+\big|^2\leq C_1\intomega\Dep^2(\nep)|\nabla\nep|^2
\end{align}
on $(0,\infty)$ for all $\epsi\in(0,1)$, because $\varphi=0$ on $B$ and $D_{1,\epsi}'=\Dep$ for all $\epsi\in(0,1)$. In view of \eqref{eq:s0} we have
\begin{align*}
D_{1,\epsi}(n)\geq \int_{s_0}^{n}\Dep(s)\intd s\geq L(n-s_0)\quad\text{for all }n\geq s_0\text{ and }\epsi\in(0,1),
\end{align*}
so that
\begin{align}\label{eq:ener-part-6-eq1}
\intomega \nep^2=\!\!\int\limits_{\{\nep\leq s_0\}}\!\!\nep^2+\!\!\int\limits_{\{\nep>s_0\}}\!\!\nep^2\leq s_0^2|\Omega|+\intomega\Big(\frac{D_{1,\epsi}(\nep)}{L}+s_0\Big)^2\leq 3s_0^2|\Omega|+\frac{2}{L^2}\intomega D_{1,\epsi}^2(\nep)
\end{align}
on $(0,\infty)$ for all $\epsi\in(0,1)$. Thus, noticing that by the positivity of $\Dep$ we have $D_{2,\epsi}(n)=\int_0^n D_{1,\epsi}(s)\intd s\leq n D_{1,\epsi}(n)$ for all $n>0$ and $\epsi\in(0,1)$, we can draw on Young's inequality and \eqref{eq:ener-part-6-eq1} to estimate
\begin{align*}
\intomega D_{2,\epsi}(\nep)&\leq \intomega \nep^2+\frac{1}{4}\intomega D_{1,\epsi}^2(\nep)\leq 3s_0^2|\Omega|+\Big(\frac{2}{L^2}+\frac{1}{4}\Big)\intomega D_{1,\epsi}^2(\nep)\\
&\leq 3s_0^2|\Omega|+\Big(\frac{2}{L^2}+\frac{1}{4}\Big)\delta_\epsi^2|\Omega|+\Big(\frac{2}{L^2}+\frac{1}{4}\Big)\intomega\big(D_{1,\epsi}(\nep)-\delta_\epsi\big)^2_+\quad\text{on }(0,\infty)\
\end{align*}
for all $\epsi\in(0,1)$. Here, plugging in \eqref{eq:ener-part-6-eq0} we arrive at
\begin{align*}
\intomega D_{2,\epsi}(\nep)\leq 3s_0^2|\Omega|+\Big(\frac{2}{L^2}+\frac{1}{4}\Big)\delta_\epsi^2|\Omega|+\Big(\frac{2}{L^2}+\frac{1}{4}\Big)C_1\intomega\Dep^2(\nep)|\nabla\nep|^2
\end{align*}
on $(0,\infty)$ for all $\epsi\in(0,1)$. Finally, noticing that due to \eqref{eq:Dep} we have $\Dep(n)\leq D(n)+2\epsi\leq D(n)+2$ for $n\geq0$ and hence also
\begin{align*}
\delta_\epsi=\int_0^{2\overline{n}_0}\Dep(s)\intd s\leq \int_0^{2\overline{n}_0} D(s)\intd s+4\overline{n}_0=D_1(2\overline{n}_0)+4\overline{n}_0,
\end{align*}
we conclude the lemma.
\end{bew}

With this we can now refine the inequality from Lemma~\ref{lem:ener-part5-small-g} into the following form.

\begin{lemma}\label{lem:ener-like-func-s-g}
Assume that \eqref{eq:prop-S} holds with some $\gamma\in[0,\tfrac12]$. For $M>0$ let $L=L(M)>0$ be provided by Lemma~\ref{lem:ener-part5-small-g}. Whenever $D$ satisfies \eqref{eq:liminf} and $\|c_0\|_{\Lo[\infty]}\leq M$, there are $\mu=\mu(M,n_0)>0$ and $\Gamma=\Gamma(M,n_0)>0$ such that for all $\epsi\in(0,1)$ the functional $\Fep$ as defined in Lemma~\ref{lem:ener-part5-small-g} satisfies
\begin{align*}
\frac{\intd}{\intd t}\Fep(t)+\mu\Fep(t)+\mu\intomega\frac{|\nabla\cep|^6}{\cep^5}+\intomega|\nabla\nep|^2\leq\Gamma\Big(\intomega|\nabla\uep|^3+\intomega|\uep|^6+1\Big)
\end{align*}
for all $t>0$.
\end{lemma}

\begin{bew}
From an application of Lemma~\ref{lem:ener-part6} we then find that there is $C_1=C_1(M,n_0)>0$ satisfying
\begin{align*}
\frac{1}{8C_1}\intomega D_{2,\epsi}(\nep)\leq \frac{1}{8}+\frac{1}{8}\intomega\Dep^2(\nep)|\nabla\nep|^2\quad\text{on }(0,\infty)\text{ for all }\epsi\in(0,1).
\end{align*}
On the other hand, with $b_i=b_i(M)>0$, $i\in\{1,2,3\}$, provided by Lemma~\ref{lem:ener-part5-small-g}, we infer from Young's inequality that
\begin{align*}
\frac{b_1}{4}\intomega\frac{\nep|\nabla\cep|^2}{\cep}\leq \frac{b_1}{4}\intomega\frac{\nep^2|\nabla\cep|^2}{\cep}+\frac{b_1}{4}\intomega\frac{|\nabla\cep|^2}{\cep}
\end{align*}
and
\begin{align*}
\frac{b_2}{4}\intomega\frac{|\nabla\cep|^4}{\cep^3}\leq \frac{b_2\delta}{8}\intomega\frac{|\nabla\cep|^6}{\cep^5}+\frac{b_2}{\delta^2}|\Omega|M
\end{align*}
on $(0,\infty)$ for all $\epsi\in(0,1)$, where $\delta>0$ is the constant from Lemma~\ref{lem:ener-part3}. Moreover, we conclude from Lemma~\ref{lem:ener-part4}, that
\begin{align*}
\frac{b_3}{4}\intomega\Psi_{2,\epsi}^{(s_0)}(\nep)\leq \frac{3s_0 b_3}{4\kappa}|\Omega|\quad\text{on }(0,\infty)\text{ for all }\epsi\in(0,1).
\end{align*}
Hence, recalling the definition of $\Fep(t)$ in \eqref{eq:ener-F-def} and the differential inequality presented in Lemma~\ref{lem:ener-part5-small-g}, we find that by adding $\frac{1}{8}+\frac{b_2}{\delta^2}|\Omega|M+\frac{3s_0 b_3}{4\kappa}|\Omega|$ to both sides of \eqref{eq:ener-part5-eq0} and estimating the terms on the left from below as presented above, we arrive at
\begin{align*}
\frac{\intd}{\intd t}\Fep(t)&+\mu\Fep(t)+\mu\intomega\frac{|\nabla\cep|^6}{\cep^5}+\intomega|\nabla\nep|^2\leq\Gamma\Big(\intomega|\nabla\uep|^3+\intomega|\uep|^6\Big)+\Gamma
\end{align*}
for all $t>0$ and $\epsi\in(0,1)$, with
\begin{align*}
\mu:=\min\Big\{\frac{1}{8C_1},\frac{1}{4},\frac{b_2\delta}{8}\Big\},\qquad\Gamma:=\max\Big\{\Gamma_1,\Gamma_2+\frac{1}{8}+\frac{b_2}{\delta^2}|\Omega|M+\frac{3s_0 b_3}{4\kappa}|\Omega|\Big\}
\end{align*}
and $\Gamma_1=\Gamma_1(M)>0$ as well as $\Gamma_2=\Gamma_2(M,n_0)>0$ taken from Lemma~\ref{lem:ener-part5-small-g}.
\end{bew}

\setcounter{equation}{0} 
\section{Dealing with \texorpdfstring{$\pmb{\gamma>\tfrac12}$}{larger values of gamma}: A second energy-like functional}\label{sec4:large-g}
If the singularity inside the chemotactic sensitivity function is stronger, i.e. $\gamma>\tfrac12$, we cannot rely on Lemma~\ref{lem:ener-part2} to directly consume the term $\intomega\frac{\nep^2|\nabla\cep|^2}{\cep^{2\gamma}}$ appearing in the proof of Lemma~\ref{lem:ener-part1-small-g}. Instead, we will have to invest the additionally assumption $\|n_0\|_{\Lo[1]}\leq M$ and fall back to splitting this term by an additional application of Young's inequality. The resulting terms on the right-hand side can then be treated by Lemma~\ref{lem:ener-part3} directly, without the necessity of considering the mixed term as in Lemma~\ref{lem:ener-part2}.

\begin{lemma}\label{lem:ener-part1-large-g}
Assume \eqref{eq:prop-S} that holds $\gamma\in[0,\tfrac56]$. For all $M>0$ there is $C=C(M)>0$ with the property: Whenever $\|c_0\|_{\Lo[\infty]}\leq M$ and $\|n_0\|_{\Lo[1]}\leq M$, then for all $\epsi\in(0,1)$ the solution $(\nep,\cep,\uep)$ of \eqref{approxprob} satisfies
\begin{align}\label{eq:ener-part1-large-g-eq0}
\frac{\intd}{\intd t}\intomega D_{2,\epsi}(\nep)+\frac{1}{2}\intomega\Dep^2(\nep)|\nabla\nep|^2\leq C\intomega|\nabla\nep|^2+C\intomega\frac{|\nabla\cep|^6}{\cep^5}+C\quad\text{on }(0,\infty).
\end{align}
\end{lemma}

\begin{bew}
Following the steps of Lemma~\ref{lem:ener-part1-small-g} we find that
\begin{align*}
\frac{\intd}{\intd t}\intomega D_{2,\epsi}(\nep)+\frac{1}{2}\intomega\Dep^2(\nep)|\nabla\nep|^2\leq \frac{S_0^2(M+1)}{2}\intomega\frac{\nep^2|\nabla\cep|^2}{\cep^{2\gamma}}
\end{align*}
on $(0,\infty)$ for all $\epsi\in(0,1)$. Then, since $\gamma\leq\tfrac56$ implies that $\frac{5}{3}-2\gamma\geq0$, we have $\|\cep\|_{\Lo[\infty]}^{\frac{5}{3}-2\gamma}\leq M^2+1$ and Young's inequality entails the existence of $C_1=C_1(M)>0$ satisfying
\begin{align}\label{eq:est-n2-nabc2-c2g-part1-eq1}
\frac{S_0^2(M+1)}{2}\intomega\frac{\nep^2|\nabla\cep|^2}{\cep^{2\gamma}}&\leq \frac{S_0^2(M+1)(M^2+1)}{2}\intomega\frac{\nep^2|\nabla\cep|^2}{\cep^\nfrac53}\leq C_1\intomega\nep^3+C_1\intomega\frac{|\nabla\cep|^6}{\cep^5}
\end{align}
on $(0,\infty)$ for all $\epsi\in(0,1)$. Here, we draw on the \GNI\ and the assumption $\|n_0\|_{\Lo[1]}\leq M$ to find $C_2>0$ and $C_3=C_3(M)>0$ such that for all $\epsi\in(0,1)$ we have
\begin{align*}
\|\nep\|_{\Lo[3]}^3\leq C_2\|\nabla\nep\|_{\Lo[2]}^2\|n_0\|_{\Lo[1]}+C_2\|n_0\|_{\Lo[1]}^3\leq C_3\|\nabla\nep\|_{\Lo[2]}^2+C_3\quad\text{on }(0,\infty).
\end{align*}
Plugging this back into \eqref{eq:est-n2-nabc2-c2g-part1-eq1} we obtain \eqref{eq:ener-part1-large-g-eq0}.
\end{bew}

Noticing that in this case the ill-signed terms can immediately be handled by Lemma~\ref{lem:ener-part3} and Lemma~\ref{lem:ener-part4}, we can establish the following.

\begin{lemma}\label{lem:ener-like-func-l-g}
Assume that \eqref{eq:prop-S} holds with some $\gamma\in[0,\tfrac56]$. Let $\delta>0$ be given by Lemma~\ref{lem:ener-part3}. For all $M>0$ one can fix $L=L(M)>0$, $\hb_2=\hb_2(M)>0$, $\hb_3=\hb_3(M)>0$ such that: Whenever $D$ satisfies \eqref{eq:liminf}, $\|c_0\|_{\Lo[\infty]}\leq M$ and $\|n_0\|_{\Lo[1]}\leq M$,  there are $s_0=s_0(L)>0$, $\hat{\mu}=\hat{\mu}(M)>0$ and $\hat{\Gamma}=\hat{\Gamma}(M)>0$ such that for all $\epsi\in(0,1)$
\begin{align}\label{eq:ener-G-def}
\Gep(t):=\intomega D_{2,\epsi}\big(\nep(\cdot,t)\big)+\hb_2\intomega\frac{\big|\nabla\cep(\cdot,t)\big|^4}{\cep^3(\cdot,t)}+\hb_3\intomega\Psi_{2,\epsi}^{(s_0)}\big(\nep(\cdot,t)\big)
\end{align}
satisfies
\begin{align*}
\frac{\intd}{\intd t}\Gep(t)+\hat{\mu}\Gep(t)+\hat{\mu}\intomega\frac{|\nabla\cep|^6}{\cep^5}+\intomega|\nabla\nep|^2\leq\hat{\Gamma}\Big(\intomega|\nabla\uep|^3+\intomega|\uep|^6+1\Big)
\end{align*}
for all $t>0$.
\end{lemma}

\begin{bew}
With $C_1=C_1(M)>0$ provided by Lemma~\ref{lem:ener-part1-large-g} and $\delta>0$ as well as $C_2=C_2(M)>0$ taken from Lemma~\ref{lem:ener-part3}, we pick $\hb_2=\hb_2(M)=\frac{4C_1}{\delta}+\frac{C_1}{C_2}+\frac{1}{C_2}$ and $\hb_3=\hb_3(M)=2C_1+2C_2\hb_2+2$ and $\eta=\eta(M)=\frac{\delta}{8C_2}>0$. Then, we let $C_3=C_3(M,\eta)>0$ denote the constant from Lemma~\ref{lem:ener-part4}, set $L=L(M)=2(C_1+C_2\hb_2+1)^\nfrac12$, pick $s_0=s_0(L)$ in accordance with \eqref{eq:s0}, i.e. such that for all $s\geq s_0$ we have $D(s)\geq L$, and denote $\kappa=\kappa(s_0):=\inf_{n\in(0,2s_0)}\frac{D(n)}{n}$. Combining Lemma~\ref{lem:ener-part1-large-g} with Lemma~\ref{lem:ener-part3} and Lemma~\ref{lem:ener-part4}, we find that
\begin{align*}
&\frac{\intd}{\intd t}\Gep(t)+\frac{1}{2}\intomega\Dep^2(\nep)|\nabla\nep|^2+\big(\hb_2\delta-C_1-\eta\hb_3\big)\intomega\frac{|\nabla\cep|^6}{\cep^5}+\intomega|\nabla\nep|^2\\
\leq\ &\big(\hat{\Gamma}_1+1\big)\intomega|\nabla\nep|^2-\big(\hat{\Gamma}_1+1\big)\!\!\int\limits_{\{\nep<s_0\}}\!\!|\nabla\nep|^2+C_2\hb_2\intomega|\nabla\uep|^3+C_2\hb_2\intomega|\uep|^6+\hat{\Gamma}_2
\end{align*}
on $(0,\infty)$ for all $\epsi\in(0,1)$, where we dropped two non-negative terms already and set $\hat{\Gamma}_1=\hat{\Gamma}_1(M):=C_1+C_2\hb_2$ and $\hat{\Gamma}_2=\hat{\Gamma}_2(M):=C_1+C_2\hb_2+\frac{C_3\hb_3}{\kappa^3}$. Noting that by the choices for $\hb_2$, $\hb_3$ and $\eta$ we have $\hb_2\delta-C_1-\eta\hb_3=\frac{\hb_2\delta}{2}$ and that $\Dep^2(s)\geq 4(\hat{\Gamma}_1+1)$ for all $s\geq s_0$ entails
\begin{align*}
\big(\hat{\Gamma}_1+1\big)\intomega|\nabla\nep|^2-\big(\hat{\Gamma}_1+1\big)\!\!\int\limits_{\{\nep<s_0\}}\!\!|\nabla\nep|^2\leq\frac{1}{4}\intomega\Dep^2(\nep)|\nabla\nep|^2
\end{align*}
on $(0,\infty)$ for all $\epsi\in(0,1)$, we have
\begin{align}\label{eq:ener-like-func-l-g-eq1}
&\frac{\intd}{\intd t}\Gep(t)+\frac{1}{4}\intomega\Dep^2(\nep)|\nabla\nep|^2+\frac{\hb_2\delta}{2}\intomega\frac{|\nabla\cep|^6}{\cep^5}+\intomega|\nabla\nep|^2
\leq\hat{\Gamma}_1\Big(\intomega|\nabla\uep|^3+\intomega|\uep|^6\Big)+\hat{\Gamma}_2
\end{align}
on $(0,\infty)$ for all $\epsi\in(0,1)$. Here, since $\overline{n}_0\leq|\Omega| M$, an application of Lemma~\ref{lem:ener-part6} provides $C_4=C_4(M)>0$ such that
\begin{align}\label{eq:ener-like-func-l-g-eq2}
\frac{1}{4C_4}\intomega D_{2,\epsi}(\nep)\leq\frac{1}{4}+\frac{1}{4}\intomega\Dep^2(\nep)|\nabla\nep|^2\quad\text{on }(0,\infty)\text{ for all }\epsi\in(0,1).
\end{align}
Moreover, Young's inequality yields
\begin{align}\label{eq:ener-like-func-l-g-eq3}
\frac{\hb_2}{4}\intomega\frac{|\nabla\cep|^4}{\cep^3}\leq\frac{\hb_2\delta}{4}\intomega\frac{|\nabla\cep|^6}{\cep^5}+\frac{\hb_2}{4\delta^2}M|\Omega|\quad\text{on }(0,\infty)\text{ for all }\epsi\in(0,1)
\end{align}
and Lemma~\ref{lem:ener-part4} entails
\begin{align}\label{eq:ener-like-func-l-g-eq4}
\frac{\hb_3}{4}\intomega\Psi_{2,\epsi}^{(s_0)}(\nep)\leq\frac{3s_0\hb_3}{4\kappa}|\Omega|\quad\text{on }(0,\infty)\text{ for all }\epsi\in(0,1).
\end{align}
Accordingly, setting $$\hat{\mu}:=\min\Big\{\frac{1}{4C_4},\frac14,\frac{\hb_2\delta}{4}\Big\}\quad\text{and}\quad\hat{\Gamma}:=\max\Big\{\hat{\Gamma}_1,\hat{\Gamma}_2+\frac14+\frac{\hb_2}{4\delta^2}M|\Omega|+\frac{3s_0\hb_3}{4\kappa}|\Omega|\Big\}$$
we conclude the claimed differential inequality from \eqref{eq:ener-like-func-l-g-eq1}--\eqref{eq:ener-like-func-l-g-eq4}.
\end{bew}

\setcounter{equation}{0} 
\section{Improved boundedness properties}\label{sec5:improved-bounds}
Let us now turn our attention to deriving useful boundedness information from the energy-like inequalities provided by Lemmas~\ref{lem:ener-like-func-s-g} and \ref{lem:ener-like-func-l-g}. Most crucially, in both cases we will establish $\epsi$-independent bounds on $\intomega\frac{|\nabla\cep|^4}{\cep^3}$ and $\int_t^{t+1}\!\intomega\frac{|\nabla\cep|^6}{\cep^5}$, which are key for further testing procedures regarding $\nep$. First, however, we need to discuss the treatment of the integrals containing the fluid-component in the energy-like inequalities.

\begin{lemma}\label{lem:st-bound-uep}
Assume that \eqref{eq:prop-S} holds with some $\gamma\in[0,\frac56]$ and let $L>0$, then there exists $C=C(L,n_0,c_0,u_0)>0$ such that if $D$ satisfies \eqref{eq:liminf}, then
\begin{align*}
\sup_{t\geq0}\intomega\big|\nabla\uep(\cdot,t)\big|^2+\sup_{t\geq0}\int_t^{t+1}\!\Big(\intomega|\nabla\uep|^3+\intomega|\uep|^6\Big)\leq C+C\sup_{t\geq0}\int_t^{t+1}\!\intomega \nep^2
\end{align*}
for all $\epsi\in(0,1)$.
\end{lemma}

\begin{bew}
We first pick $\gamma_0=\frac{5}{6}$ and note that Lemma~\ref{lem:uep-bound} provides $C_0=C_0(L,n_0,c_0,u_0)>0$ satisfying
$$\int_{t}^{t+1}\!\intomega\big|\nabla\uep(\cdot,t)\big|^2\leq C_0\quad\text{for all }t>0\text{ and }\epsi\in(0,1)$$ and, accordingly, for any $t\geq0$ and each $\epsi\in(0,1)$ we can find $t_\epsi\in[t,t+1]$ such that
\begin{align}\label{eq:upe-nabcep-linfty-eq1}
\intomega\big|\nabla\uep(\cdot,t_\epsi)\big|^2\leq C_0.
\end{align}
Then, applying the Helmholtz projection to the third equation of \eqref{approxprob} and testing by $A\uep$ we find from standard arguments (see e.g. \cite[p.21]{win_fluid_CPDE12} and \cite[Lemma 3.10]{heihoffTwoNewFunctional2023}) that there is some $C_1>0$ fulfilling
\begin{align}\label{eq:st-bound-uep-eq1}
\frac{\intd}{\intd t}\intomega\big|\nabla\uep(\cdot,t)\big|^2+\frac{1}{2}\intomega\big|A\uep(\cdot,t)\big|^2\leq C_1\intomega\nep^2(\cdot,t)+C_1\Big(\intomega\big|\nabla\uep(\cdot,t)\big|^2\Big)^2
\end{align}
holds on $(0,\infty)$ for all $\epsi\in(0,1)$. Abbreviating $y_\epsi(t):=\gep(t):=\intomega|\nabla\uep(\cdot,t)|^2$ we may rewrite the above as
\begin{align*}
y_\epsi'(t)\leq C_1\intomega\nep^2(\cdot,t)+C_1\gep(t)y_\epsi(t),\quad\text{for all }t>0\text{ and }\epsi\in(0,1),
\end{align*}
with $\int_t^{t+1}y_\epsi(s)\intd s\leq C_0$ as well as $\int_t^{t+1}\gep(s)\intd s\leq C_0$ for all $t>0$ and $\epsi\in(0,1)$. An ODE-comparison argument first entails that with $C_2=C_2(L,n_0,u_0):=\max\{\intomega|\nabla u_0|^2,C_0e^{C_1C_0}\}$ we have 
\begin{align}\label{eq:st-bound-uep-eq15}
y_\eps(t)\leq C_2+e^{C_1C_0}C_1\sup_{t\geq0}\int_t^{t+1}\!\intomega\nep^2(\cdot,s)\intd s\quad\text{for all }t>0\text{ and }\epsi\in(0,1).
\end{align}
Thereafter, we conclude from \eqref{eq:st-bound-uep-eq1} that there is $C_3=C_3(L,n_0,c_0,u_0)>0$ such that
\begin{align}\label{eq:st-bound-uep-eq2}
\int_t^{t+1}\!\intomega\big|A\uep(\cdot,t)\big|^2\leq C_3+C_3\sup_{t\geq0}\int_t^{t+1}\!\intomega\nep^2\quad\text{for all }t>0\text{ and }\epsi\in(0,1).
\end{align}
Drawing on the \GNI\ we then find $C_4>0$ such that
\begin{align*}
\intomega|\nabla\uep|^3\leq \|\uep\|_{\W[1,3]}^3\leq C_4\|\uep\|_{\W[2,2]}^2\|\uep\|_{\Lo[2]}
\end{align*}
and
\begin{align*}	
\intomega|\uep|^6\leq C_4\|\uep\|_{\W[2,2]}^2\|\uep\|_{\Lo[2]}^4
\end{align*}
hold on $(0,\infty)$ for all $\epsi\in(0,1)$ and conclude from \cite[Theorem IV.6.1]{galdiIntroductionMathematicalTheory2011} and Lemma~\ref{lem:uep-bound} that
\begin{align*}
\sup_{t\geq 0}\int_t^{t+1}\!\Big(\intomega |\nabla\uep|^3+\intomega |\uep|^6\Big)\leq C_5\sup_{t\geq 0}\int_t^{t+1}\!\intomega\big| A\uep\big|^2
\end{align*}
for all $\epsi\in(0,1)$ with some $C_5=C_5(L,n_0,c_0,u_0)>0$, which together with with \eqref{eq:st-bound-uep-eq15} and \eqref{eq:st-bound-uep-eq2} directly entails the assertion.
\end{bew}

This lemma at hand, we return to our energy-like functionals from Lemma~\ref{lem:ener-like-func-s-g} and Lemma~\ref{lem:ener-like-func-l-g} to improve upon our regularity information for $\nabla\cep$ as mentioned before.

\begin{lemma}\label{lem:nab-cep-l4-bound}
Assume that \eqref{eq:prop-S} holds with some $\gamma\in[0,\frac56]$ and let $M>0$. If $\gamma\in[0,\tfrac12)$ assume that $D$ satisfies \eqref{eq:liminf} with $L(M)>0$ provided by Lemma~\ref{lem:ener-part1-small-g} and that $\|c_0\|_{\Lo[\infty]}\leq M$. If instead $\gamma\in(\tfrac12,\tfrac56]$ assume that $D$ satisfies \eqref{eq:liminf} with $L(M)>0$ provided by Lemma~\ref{lem:ener-like-func-l-g} and that $\|c_0\|_{\Lo[\infty]}\leq M$ and  $\|n_0\|_{\Lo[1]}\leq M$. Then, there is $C=C(n_0,c_0,u_0)>0$ such that
\begin{align*}
\intomega\frac{\big|\nabla\cep(\cdot,t)\big|^4}{\cep^3(\cdot,t)}\leq C\quad\text{and}\quad\intomega\big|\nabla\cep(\cdot,t)\big|^4\leq C
\end{align*}
for all $t>0$ and all $\epsi\in(0,1)$ and such that moreover
\begin{align*}
\sup_{t\geq0}\left(\int_t^{t+1}\!\!\intomega\frac{|\nabla\cep|^6}{\cep^5}+\int_t^{t+1}\!\!\intomega|\nabla\nep|^2\right)\leq C\quad\text{for all }\epsi\in(0,1).
\end{align*}
\end{lemma}

\begin{bew}
If $\gamma\in[0,\tfrac12]$ and $\|c_0\|_{\Lo[\infty]}\leq M$, we let $L=L(M)>0$, $b_1=b_1(M)>0$, $b_2=b_2(M)>0$, $b_3=b_3(M)>0$ and $s_0=s_0(M)>0$ and $\Fep$ as in Lemma~\ref{lem:ener-part5-small-g} and conclude from Lemma~\ref{lem:ener-like-func-s-g} that there are $\mu=\mu(M,n_0)>0$ and $\Gamma=\Gamma(M,n_0)>0$ such that 
\begin{align*}
\frac{\intd}{\intd t}\Fep(t)+\mu\Fep(t)+\mu\intomega\frac{|\nabla\cep|^6}{\cep^5}+\intomega|\nabla\nep|^2\leq \Gamma\gep(t)\quad\text{for all }\epsi\in(0,1)\text{ and }t>0,
\end{align*}
where we abbreviated $\gep(t):=\intomega|\nabla\uep(\cdot,t)|^3+\intomega|\uep(\cdot,t)|^6+1$. Similarly, if $\gamma\in(\tfrac12,\tfrac56]$, $\|c_0\|_{\Lo[\infty]}\leq M$ and additionally $\|n_0\|_{\Lo[1]}\leq M$, we let $L=L(M)>0$, $\hb_1=\hb_1(M)>0$, $\hb_2=\hb_2(M)>0$, $s_0=s_0(M)$ and $\Gep$ as in Lemma~\ref{lem:ener-like-func-l-g} and obtain $\hat{\mu}=\hat{\mu}(M)>0$ and $\hat{\Gamma}=\hat{\Gamma}(M)>0$ satisfying
\begin{align*}
\frac{\intd}{\intd t}\Gep(t)+\hat{\mu}\Gep(t)+\hat{\mu}\intomega\frac{|\nabla\cep|^6}{\cep^5}+\intomega|\nabla\nep|^2\leq \hat{\Gamma}\gep(t)\quad\text{for all }\epsi\in(0,1)\text{ and }t>0.
\end{align*}
Fixing $\gamma_0=\frac{5}{6}$, in both cases for $L$, we can draw on Lemma~\ref{lem:st-bound-uep} to find $C_1=C_1(M,n_0,c_0,u_0)>0$ such that we have
\begin{align*}
\sup_{t\geq 0}\int_t^{t+1}\!\gep(s)\intd s\leq C_1+C_1\sup_{t\geq0}\int_t^{t+1}\!\!\intomega\nep^2\quad\text{for all }\epsi\in(0,1).
\end{align*}
Therefore, an ODE-comparison argument (see e.g. \cite[Lemma 3.4]{winklerThreedimensionalKellerSegel2019} and \cite[Lemma 2.1]{winklerLpboundLocalSensing}) entails that there is $C_2=C_2(M,n_0,c_0,u_0)>0$ such that if $\gamma\in[0,\tfrac12]$, we have
\begin{align}\label{eq:nab-cep-l4-bound-eq1}
&\sup_{t\geq0}\Fep(t)+\sup_{t\geq 0}\left(\mu\int_t^{t+1}\!\!\intomega\frac{|\nabla\cep|^6}{\cep^5}+\int_t^{t+1}\!\!\intomega|\nabla\nep|^2\right)\leq C_2+C_2\sup_{t\geq0}\int_t^{t+1}\!\!\intomega\nep^2
\intertext{for all $\epsi\in(0,1)$ and such that if $\gamma\in(\tfrac12,\tfrac56]$, we have}
\label{eq:nab-cep-l4-bound-eq2}
&\sup_{t\geq0}\Gep(t)+\sup_{t\geq 0}\left(\hat{\mu}\int_t^{t+1}\!\!\intomega\frac{|\nabla\cep|^6}{\cep^5}+\int_t^{t+1}\!\!\intomega|\nabla\nep|^2\right)\leq C_2+C_2\sup_{t\geq0}\int_t^{t+1}\!\!\intomega\nep^2
\end{align}
for all $\epsi\in(0,1)$. Here, we note that according to the \GNI\ there is $C_3=C_3(M,n_0,c_0,u_0)>0$ such that for all $\epsi\in(0,1)$
\begin{align*}
C_2\|\nep\|_{\Lo[2]}^2&\leq C_3\|\nabla\nep\|_{\Lo[2]}\|\nep\|_{\Lo[1]}+C_3\|\nep\|_{\Lo[1]}^2
\end{align*}
holds on $(0,\infty)$. Hence, by means of Young's inequality and Lemma~\ref{lem:mass}, we see that
\begin{align*}
C_2\|\nep\|_{\Lo[2]}^2&\leq \frac12\|\nabla\nep\|_{\Lo[2]}^2+\Big(C_3+\frac{C_3^2}{2}\Big)\|\nep\|_{\Lo[1]}^2\leq \frac12\|\nabla\nep\|_{\Lo[2]}^2+\Big(C_3+\frac{C_3^2}{2}\Big)\|n_0\|_{\Lo[1]}^2
\end{align*}
for all $\epsi\in(0,1)$ on $(0,\infty)$ and plugging this back into \eqref{eq:nab-cep-l4-bound-eq1} and \eqref{eq:nab-cep-l4-bound-eq2}, respectively, yields $C_4=C_4(M,n_0,c_0,u_0)>0$ such that for all $\epsi\in(0,1)$
\begin{align*}
&\sup_{t\geq 0}\Fep(t)+\sup_{t\geq 0}\left(\mu\int_t^{t+1}\!\!\intomega\frac{|\nabla\cep|^6}{\cep^5}+\frac12\int_t^{t+1}\!\!\intomega|\nabla\nep|^2\right)\leq C_4,\\
\text{and}\ &\sup_{t\geq0}\Gep(t)+\sup_{t\geq 0}\left(\hat{\mu}\int_t^{t+1}\!\!\intomega\frac{|\nabla\cep|^6}{\cep^5}+\frac12\int_t^{t+1}\!\!\intomega|\nabla\nep|^2\right)\leq C_4\quad\text{on }(0,\infty).
\end{align*}
Recalling the definition of $\Fep$ and $\Gep$ from \eqref{eq:ener-F-def} and \eqref{eq:ener-G-def}, respectively, this readily entails that for some $C_5=C_5(M,n_0,c_0,u_0)>0$
\begin{align*}
\sup_{t\geq0}\intomega\frac{|\nabla\cep|^4}{\cep^3}\leq C_5\quad\text{and}\quad\sup_{t\geq0}\left(\int_t^{t+1}\!\!\intomega\frac{|\nabla\cep|^6}{\cep^5}+\int_t^{t+1}\!\!\intomega|\nabla\nep|^2\right)\leq C_5\quad\text{for all }\epsi\in(0,1),
\end{align*}
where the former, by rewriting and relying on Lemma~\ref{lem:mass} once more, also shows that
\begin{align*}
\intomega|\nabla\cep|^4=\intomega\frac{|\nabla\cep|^4\cep^3}{\cep^3}\leq \|c_0\|_{\Lo[\infty]}^3\intomega\frac{|\nabla\cep|^4}{\cep^3}\leq C_5\|c_0\|_{\Lo[\infty]}^3
\end{align*}
on $(0,\infty)$ for all $\epsi\in(0,1)$, confirming the last open claim.
\end{bew}

In \cite[Lemma 4.2]{winklerChemotaxisStokesInteractionVery2022} the time-uniform control on $\intomega\frac{|\nabla\cep|^4}{\cep^3}$ was exploited to derive $\Lo[p]$-bounds for $\nep$. If we were to mimic the proof in our setting, we would, however, have to further restrict the range of permissible $\gamma$. Instead, deviating slightly from the steps in the reference, we will rely on the leeway provided by the bound on $\int_t^{t+1}\!\!\intomega\frac{|\nabla\cep|^6}{\cep^5}$ to build an appropriate bootstrap argument. Let us prepare the reasoning with the following lemma (cf. \cite[Lemma 3.3]{Win24_ev-reg-2d-doubly-deg-nutrient}).

\begin{lemma}\label{lem:Mk}
Let $a\geq 1$, $b\geq 1$ and $(\M_k)_{k\geq0}\subset[1,\infty)$ be such that
\begin{align*}
\M_k\leq a^k \M_{k-1}^2+b^{2^k}\quad\text{for all }k\geq 1.
\end{align*}
Then,
\begin{align*}
\liminf_{k\to\infty} \M_k^{\frac{1}{2^k}}\leq 2\sqrt{2}a^3 b \M_0.
\end{align*}
\end{lemma}

\begin{lemma}\label{lem:moser-it}
Assume that \eqref{eq:prop-S} holds with some $\gamma\in[0,\frac56]$ and let $M>0$. If $\gamma\in[0,\tfrac12)$ assume that $D$ satisfies \eqref{eq:liminf} with $L(M)>0$ provided by Lemma~\ref{lem:ener-part1-small-g} and that $\|c_0\|_{\Lo[\infty]}\leq M$. If instead $\gamma\in(\tfrac12,\tfrac56]$ assume that $D$ satisfies \eqref{eq:liminf} with $L(M)>0$ provided by Lemma~\ref{lem:ener-like-func-l-g} and that $\|c_0\|_{\Lo[\infty]}\leq M$ and  $\|n_0\|_{\Lo[1]}\leq M$. Then, there is $K=K(M,n_0,c_0,u_0)>0$ such that
\begin{align*}
\big\|\nep(\cdot,t)\big\|_{\Lo[\infty]}\leq K\quad\text{for all }t>0\text{ and }\epsi\in(0,1).
\end{align*}
\end{lemma}

\begin{bew}
For $k\geq0$ we set $p_k:=2^k$ and test the first equation of \eqref{approxprob} against $(\nep-1)_+^{p_k-1}$ to find that
\begin{align*}
\frac{\intd}{\intd t}\intomega(\nep-1)_+^{p_k}=&-p_k(p_k-1)\intomega\Dep(\nep)(\nep-1)_+^{p_k-2}|\nabla\nep|^2\\
&+p_k(p_k-1)\intomega(\nep-1)_+^{p_k-2}\nep\big(\Sep(x,\nep,\cep)\nabla\cep\cdot\nabla\nep\big)\quad\text{on }(0,\infty)
\end{align*}
holds for all $\epsi\in(0,1)$ and $k\geq1$, where we used the prescribed boundary conditions and the divergence-free property of $\uep$. As \eqref{eq:liminf} implies that $d:=\inf_{n\geq 1} D(n)$ is strictly positive, we conclude from \eqref{eq:Dep} that $\Dep(\nep)\geq d>0$ on $\{\nep\geq 1\}$, so that by Young's inequality and the fact that $p_k\geq 2$ for $k\geq1$, we have
\begin{align}\label{eq:nep-lp-eq0}
&\frac{\intd}{\intd t}\intomega(\nep-1)_+^{p_k}\leq-\frac{p_k^2 d}{4}\intomega(\nep-1)_+^{p_k-2}|\nabla\nep|^2+\frac{p_k^2}{d}\intomega (\nep-1)_+^{p_k-2}\nep^2\big|\Sep(x,\nep,\cep)\nabla\cep\big|^2
\end{align}
on $(0,\infty)$ for all $\epsi\in(0,1)$ and $k\geq1$. Abbreviating $\fep:=\Sep(x,\nep,\cep)\nabla\cep$, we make use of the fact that $\sigma^2(\sigma-1)_+^{p_k-2}\leq 4(\sigma-1)_+^{p_k}+2$ for all $\sigma\geq 0$ and $k\geq1$ and Hölder's inequality to obtain that
\begin{align*}
\intomega(\nep-1)_+^{p_k-2}\nep^2\big|\fep\big|^2&\leq 4\big\|(\nep-1)_+^{\frac{p_k}{2}}\fep\big\|_{\Lo[2]}^2+2\big\|\fep\big\|_{\Lo[2]}^2\\
&\leq 4\big\|\fep\big\|_{\Lo[6]}^2\big\|(\nep-1)_+^\frac{p_k}{2}\big\|_{\Lo[3]}^2+2|\Omega|^\frac23\big\|\fep\big\|_{\Lo[6]}^2
\end{align*}
on $(0,\infty)$ for all $\epsi\in(0,1)$ and $k\geq1$. Adding $\intomega(\nep-1)_+^{p_k}$ to both sides of \eqref{eq:nep-lp-eq0} and estimating $$\intomega(\nep-1)_+^{p_k}\leq |\Omega|^\frac{1}{3}\big\|(\nep-1)_+^\frac{p_k}{2}\big\|_{\Lo[3]}^2\quad\text{on }(0,\infty)\text{ for all }\epsi\in(0,1)\text{ and }k\geq 1,$$
we hence find that
\begin{align*}
\frac{\intd}{\intd t}\intomega(\nep-1)_+^{p_k}&+d\intomega\big|\nabla(\nep-1)_+^\frac{p_k}{2}\big|^2+\intomega(\nep-1)_+^{p_k}\nonumber\\
&\leq \big(4\big\|\fep\big\|_{\Lo[6]}^2+|\Omega|^\frac13\big)\frac{p_k^2}{2d}\big\|(\nep-1)_+^\frac{p_k}{2}\big\|_{\Lo[3]}^2+2|\Omega|^\frac{4}{3}\frac{p_k^2}{2d}\big\|\fep\big\|_{\Lo[6]}^2
\end{align*}
on $(0,\infty)$ for all $\epsi\in(0,1)$ and $k\geq1$. In light of the \GNI, there is $C_1>0$ such that 
\begin{align*}
\big\|(\nep-1)_+^\frac{p_k}{2}\big\|_{\Lo[3]}^2\leq C_1\big\|\nabla(\nep-1)_+^\frac{p_k}{2}\big\|_{\Lo[2]}^\frac{4}{3}\big\|(\nep-1)_+^\frac{p_k}{2}\big\|_{\Lo[1]}^\frac{2}{3}+C_1\big\|(\nep-1)_+^\frac{p_k}{2}\big\|_{\Lo[1]}^2
\end{align*}
on $(0,\infty)$ for all $\epsi\in(0,1)$ and $k\geq1$, so that an application of Young's inequality entails that there is $C_2
>0$ satisfying
\begin{align}\label{eq:nep-lp-eq1}
\frac{\intd}{\intd t}\intomega(\nep&-1)_+^{p_k}+\intomega(\nep-1)_+^{p_k}\nonumber\\&\leq C_2p_k^6\big\|\fep\big\|_{\Lo[6]}^6\big\|(\nep-1)_+^\frac{p_k}{2}\big\|_{\Lo[1]}^2+C_2 p_k^2\big\|\fep\big\|_{\Lo[6]}^2\big\|(\nep-1)_+^\frac{p_k}{2}\big\|_{\Lo[1]}^2\\&\qquad+C_2 p_k^2\big\|(\nep-1)_+^\frac{p_k}{2}\big\|_{\Lo[1]}^2+C_2 p_k^6\big\|\fep\big\|_{\Lo[6]}^2\quad\text{on }(0,\infty)\text{ for all }\epsi\in(0,1)\text{ and }k\geq 1.\nonumber
\end{align}
Introducing $$\M_k:=\M_{\epsi,k}:=1+\sup_{t>0}\intomega(\nep-1)_+^{p_k},\quad \epsi\in(0,1), k\geq 0,$$ by definition of $p_k$ we have $\big\|(\nep-1)_+^\frac{p_k}{2}\big\|_{\Lo[1]}\leq \M_{k-1}$ on $(0,\infty)$ for all $\epsi\in(0,1)$ and $k\geq1$, so that upon letting $y_k(t):=y_{\epsi,k}(t):=\intomega(\nep-1)_+^{p_k}$ we infer from \eqref{eq:nep-lp-eq1} that there is $C_3
:=2C_2>0$ such that
\begin{align*}
y'_{k}(t)+y_{k}(t)\leq C_3 p_k^6\big\|\fep\big\|_{\Lo[6]}^6 \M_{k-1}^2+C_3 p_k^2\big\|\fep\big\|_{\Lo[6]}^2 \M_{k-1}^2+C_3 p_k^2 \M_{k-1}^2
\end{align*}
holds on $(0,\infty)$ for all $\epsi\in(0,1)$ and $k\geq1$. Here, since \eqref{eq:prop-S} implies that $$\big\|\fep\big\|_{\Lo[6]}^6\leq S_0^6(M+1)\intomega\frac{|\nabla\cep|^6}{\cep^{6\gamma}}\leq S_0^6(M+1)\|\cep\|_{\Lo[\infty]}^{5-6\gamma}\intomega\frac{|\nabla\cep|^6}{\cep^5}\quad\text{on }(0,\infty)\text{ for all }\epsi\in(0,1),$$ we infer from $\gamma\leq\frac{5}{6}$, Lemma~\ref{lem:mass} and Lemma~\ref{lem:nab-cep-l4-bound} that there is $C_4=C_4(M,n_0,c_0,u_0)>0$ satisfying
\begin{align*}
\int_t^{t+1}\big\|\fep(\cdot,s)\big\|_{\Lo[6]}^6\intd s+\int_t^{t+1}\big\|\fep(\cdot,s)\big\|_{\Lo[6]}^2\intd s\leq C_4\quad\text{for all }t>0\text{ and }\epsi\in(0,1).
\end{align*}
Letting $$h_k(s):=h_{\epsi,k}(s):=C_3 p_k^6 \big\|\fep(\cdot,s)\big\|_{\Lo[6]}^6\M_{k-1}^2+C_3 p_k^2 \big\|\fep(\cdot,s)\big\|_{\Lo[6]}^2\M_{k-1}^2+C_3 p_k^2\M_{k-1}^2,$$
we conclude from the above and $p_k\geq 1$, that with $C_5=C_5(
M,n_0,c_0,u_0):=C_3(C_4+1)>0$ we have $$\int_t^{t+1} h_k(s)\intd s\leq C_5 p_k^6\M_{k-1}^2\quad\text{on }(0,\infty)\text{ for all }\epsi\in(0,1)\text{ and }k\geq 1.$$ Accordingly, we may draw on Lemma~\ref{lem:ode} to find that $$y_k\leq e^{-t}\intomega (n_0-1)_+^{p_k}+\frac{C_5}{1-e^{-1}} p_k^6\M_{k-1}^2\quad\text{for all }t\geq0,\ \epsi\in(0,1)\text{ and }k\geq1,$$
from which, using precise form of $p_k=2^k$, we derive that $$\M_k\leq b^{2^k}+a^k \M_{k-1}^2\quad\text{for all }\epsi\in(0,1)\text{ and }k\geq 1,$$
with $$a=a(M,n_0,c_0,u_0):=\max\{2C_5,1\}2^6\quad\text{and}\quad b=b(n_0):=1+\|(n_0-1)_+\|_{\Lo[\infty]}\max\{|\Omega|,1\}.$$
Here, Lemma~\ref{lem:Mk} applies to say that actually
$$\liminf_{k\to\infty}\M_{\epsi,k}^\frac{1}{2^k}\leq 2\sqrt{2}a^3 b\M_0\quad\text{for all }\epsi\in(0,1).$$
In view of Lemma~\ref{lem:mass} we find that $\M_{0}\leq 1+\intomega n_0$, so that indeed
\begin{align*}
\big\|\big(n(\cdot,t)-1\big)_+\big\|_{\Lo[\infty]}&=\liminf_{k\to\infty}\Big(\intomega(n-1)_+^{2^k}\Big)^{\frac{1}{2^k}}\\
&\leq \liminf_{k\to\infty}\M_{\epsi,k}^\frac{1}{2^k} \leq 2\sqrt{2}a^3 b\Big(1+\intomega n_0\Big)\quad\text{on }(0,\infty)\text{ for all }\epsi\in(0,1).
\end{align*}
The claim is an evident consequence of this.
\end{bew}

With $\nep$ being bounded, we can now also obtain an $\epsi$-independent positive lower bound for $\cep$.

\begin{lemma}\label{lem:lower-c}
Assume that \eqref{eq:prop-S} holds with some $\gamma\in[0,\frac56]$ and let $M>0$. If $\gamma\in[0,\tfrac12)$ assume that $D$ satisfies \eqref{eq:liminf} with $L(M)>0$ provided by Lemma~\ref{lem:ener-part1-small-g} and that $\|c_0\|_{\Lo[\infty]}\leq M$. If instead $\gamma\in(\tfrac12,\tfrac56]$ assume that $D$ satisfies \eqref{eq:liminf} with $L(M)>0$ provided by Lemma~\ref{lem:ener-like-func-l-g} and that $\|c_0\|_{\Lo[\infty]}\leq M$ and  $\|n_0\|_{\Lo[1]}\leq M$. Denote by $K=K(M,n_0,c_0,u_0)>0$ the constant obtained in Lemma~\ref{lem:moser-it}, then
\begin{align*}
\cep(x,t)\geq \big(\inf_{x\in\Omega}c_0\big)e^{-K t}\quad\text{for all }(x,t)\in\Omega\times(0,\infty)\text{ and all }\epsi\in(0,1).
\end{align*}
\end{lemma}

\begin{bew}
With $K=K(M,n_0,c_0,u_0)>0$ given by Lemma~\ref{lem:moser-it}, we have
\begin{align*}
\nep \leq K\quad\text{in }\Omega\times(0,\infty)\text{ for all }\epsi\in(0,1).
\end{align*}
Letting $\underline{c}(t):=\big(\inf_{x\in\Omega} c_0\big) e^{-K t}$ we have 
$$\underline{c}_t-\Delta\underline{c}+\nep\underline{c}+\uep\nabla\underline{c}\leq0\quad
\text{and}\quad\underline{c}(0)=\inf_{x\in\Omega} c_0\leq \cep(\cdot,0).$$ The assertion is therefore an evident consequence of parabolic comparison.
\end{bew}

Combining Lemma~\ref{lem:moser-it} with Lemma~\ref{lem:st-bound-uep} also readily entails the following.
\begin{corollary}\label{cor:nab-uep-l2}
Assume that \eqref{eq:prop-S} holds with some $\gamma\in[0,\frac56]$ and let $M>0$. If $\gamma\in[0,\tfrac12)$ assume that $D$ satisfies \eqref{eq:liminf} with $L(M)>0$ provided by Lemma~\ref{lem:ener-part1-small-g} and that $\|c_0\|_{\Lo[\infty]}\leq M$. If instead $\gamma\in(\tfrac12,\tfrac56]$ assume that $D$ satisfies \eqref{eq:liminf} with $L(M)>0$ provided by Lemma~\ref{lem:ener-like-func-l-g} and that $\|c_0\|_{\Lo[\infty]}\leq M$ and  $\|n_0\|_{\Lo[1]}\leq M$. Then, there is $C=C(M,n_0,c_0,u_0)>0$ such that for all $\epsi\in(0,1)$ the solution $(\nep,\cep,\uep)$ of \eqref{approxprob} satisfies
\begin{align*}
\big\|\nabla\uep(\cdot,t)\big\|_{\Lo[2]}\leq C\quad\text{for all }t>0.
\end{align*}
\end{corollary}

Next, we will make use of the gradient information on the fluid-flow, to derive $L^\infty$-boundedness of $\uep$ itself and afterwards also of $\nabla\cep$. These bounds will later provide a way to transfer the implied regularity to the limit functions constructed in Section~\ref{sec7:solution}.

\begin{lemma}\label{lem:A-uep-Lr}
Assume that \eqref{eq:prop-S} holds with some $\gamma\in[0,\frac56]$ and let $M>0$. If $\gamma\in[0,\tfrac12)$ assume that $D$ satisfies \eqref{eq:liminf} with $L(M)>0$ provided by Lemma~\ref{lem:ener-part1-small-g} and that $\|c_0\|_{\Lo[\infty]}\leq M$. If instead $\gamma\in(\tfrac12,\tfrac56]$ assume that $D$ satisfies \eqref{eq:liminf} with $L(M)>0$ provided by Lemma~\ref{lem:ener-like-func-l-g} and that $\|c_0\|_{\Lo[\infty]}\leq M$ and  $\|n_0\|_{\Lo[1]}\leq M$. Let $\fracpow\in(\tfrac12,1)$ be given by \eqref{IR}, then for all $\beta\in(\frac12,\fracpow]$ and $r\in(2,\infty)$ there is $C=C(M,n_0,c_0,u_0,\beta,r)>0$ such that for all $\epsi\in(0,1)$ the solution $(\nep,\cep,\uep)$ of \eqref{approxprob} satisfies
\begin{align*}
\big\|A^\beta\uep(\cdot,t)\big\|_{\Lo[r]}\leq C\quad\text{for all }t>0.
\end{align*}
\end{lemma}

\begin{bew}
In view of the Sobolev embedding theorem, we conclude from Corollary~\ref{cor:nab-uep-l2} that for any $\vartheta\geq 1$
\begin{align}\label{eq:A-uep-Lr-eq1}
\intomega\big|\uep(\cdot,t)\big|^\vartheta\leq C_1\quad\text{for all }t>0\text{ and }\epsi\in(0,1),
\end{align}
with some $C_1=C_1(M,n_0,c_0,u_0,\vartheta)>0$. With $\fracpow\in(\tfrac12,1)$ provided by \eqref{IR}, $\beta\in(\tfrac12,\fracpow]$ and $r\in(2,\infty)$ we pick $q=q(\beta,r)\in(2,r)$ close to $r$ such that still $\beta+\frac{1}{q}-\frac{1}{r}<1$ holds. Since $\beta>\tfrac12$ we have $1-\frac{2}{r}<2\beta-\frac{2}{r}$ and known embedding properties on the domains of fractional powers of the Stokes operator (see e.g. \cite[Theorem 1.6.1]{hen81} and \cite[Lemma 2.3]{caolan16_smalldatasol3dnavstokes}) provide $C_2=C_2(\beta,r)>0$ such that
\begin{align*}
\|\varphi\|_{\W[1,r]}\leq C_2\|A^\beta\varphi\|_{\Lo[r]}\quad\text{for all }\varphi\in D\big(A_r^\beta\big).
\end{align*}
Therefore, fixing some $q\in(2,r)$, we obtain from the interpolating by means of Hölder's inequality and Corollary~\ref{cor:nab-uep-l2} that there is $C_3=C_3(M,n_0,c_0,u_0,\beta,r)>0$ satisfying
\begin{align}\label{eq:A-uep-Lr-eq2}
\|\nabla\uep\|_{\Lo[\frac{r+q}{2}]}\leq\|\nabla\uep\|_{\Lo[2]}^{1-a}\|\nabla\uep\|_{\Lo[r]}^a\leq C_3\big\|A^\beta\uep\|_{\Lo[r]}^a\quad\text{for all }t>0\text{ and }\epsi\in(0,1),
\end{align}
where $a=\frac{\frac{1}{2}-\frac{2}{r+q}}{\frac12-\frac1r}\in(0,1)$. With this, we denote by $\Po=\Po_r$ the Helmholtz projection from $\LSp{r}{\Omega;\R^2}$ to $L^r_\sigma(\Omega)$ and draw on semigroup estimates for the Stokes semigroup (cf. \cite{gig86} and \cite[Lemma 2.3]{caolan16_smalldatasol3dnavstokes}) and Hölder's inequality to obtain $C_4=C_4(\beta,r)>0$ and $C_5(\beta,r)>0$ such that
\begin{align*}
\big\|A^\beta\uep(\cdot,t)\big\|_{\Lo[r]}&\leq \big\|A^\beta e^{-tA} u_0\big\|_{\Lo[r]}+C_4\int_0^t (t-s)^{-\fracpow}e^{-\mu t}\|\Po(\nep(\cdot,s)\nabla\phi)\|_{\Lo[r]}\intd s\\
&\hspace*{1.8cm}+C_4\int_0^t(t-s)^{-\beta-\frac1q+\frac1r}e^{-\mu t}\big\|\Po\big((\uep\cdot\nabla)\uep\big)(\cdot,s)\big\|_{\Lo[q]}\intd s\\
&\leq C_5\big\|A^\fracpow u_0\big\|_{\Lo[r]}+C_5\int_0^t (t-s)^{-\fracpow}e^{-\mu t}\|\nabla\Phi\|_{\Lo[\infty]}\|\nep\|_{\Lo[r]}\intd s\\
&\hspace*{1.8cm}+C_5\int_0^t(t-s)^{-\beta-\frac1q+\frac1r}e^{-\mu t}\|\uep(\cdot,s)\|_{\Lo[\frac{q(r+q)}{r-q}]}\|\nabla\uep(\cdot,s)\|_{\Lo[\frac{r+q}{2}]}\intd s
\end{align*}
is valid for all $t>0$ and any $\epsi\in(0,1)$, where we used the contractivity of the Stokes semigroup, that $D(A^\fracpow_r)\hookrightarrow D(A^\beta_r)$ due to $\beta<\fracpow$ and that $\Po$ is bounded on $\LSp{r}{\Omega;\R^2}$. Accordingly, relying on Lemma~\ref{lem:moser-it}, \eqref{IR}, \eqref{eq:A-uep-Lr-eq1} and \eqref{eq:A-uep-Lr-eq2}, we infer that there is $C_6=C_6(M,n_0,c_0,u_0,\beta,r)>0$ such that $M_\epsi(t):=\sup_{t>0}\big\|A^\beta\uep(\cdot,t)\big\|_{\Lo[r]}$ satisfies
\begin{align*}
M_\epsi(t)\leq C_6+C_6 M_\epsi^a(t)\quad\text{for all }\epsi\in(0,1).
\end{align*}
Due to $a<1$ this readily leads to the claimed bound.
\end{bew}

\begin{corollary}\label{cor:uep-nabcep-linfty}
Assume that \eqref{eq:prop-S} holds with some $\gamma\in[0,\frac56]$ and let $M>0$. If $\gamma\in[0,\tfrac12)$ assume that $D$ satisfies \eqref{eq:liminf} with $L(M)>0$ provided by Lemma~\ref{lem:ener-part1-small-g} and that $\|c_0\|_{\Lo[\infty]}\leq M$. If instead $\gamma\in(\tfrac12,\tfrac56]$ assume that $D$ satisfies \eqref{eq:liminf} with $L(M)>0$ provided by Lemma~\ref{lem:ener-like-func-l-g} and that $\|c_0\|_{\Lo[\infty]}\leq M$ and  $\|n_0\|_{\Lo[1]}\leq M$. Then, there is $C=K(M,n_0,c_0,u_0)>0$ such that
\begin{align*}
\big\|\uep(\cdot,t)\big\|_{\Lo[\infty]}\leq C\quad\text{and}\quad \big\|\nabla\cep(\cdot,t)\big\|_{\Lo[\infty]}\leq C\quad\text{for all }t>0\text{ and }\epsi\in(0,1).
\end{align*}
\end{corollary}

\begin{bew}
Employing Lemma~\ref{lem:A-uep-Lr} for $\beta=\fracpow$ and $r>2$ entails the existence of $C_1=C_1(M,n_0,c_0,u_0)>0$ such that $\|A^\fracpow\uep(\cdot,t)\|_{\Lo[r]}\leq C_1$ for all $t>0$ and all $\epsi\in(0,1)$. By these choices we moreover have $1-\frac{2}{r}<2\beta-\frac{2}{r}$ and hence $D(A_r^\fracpow)\hookrightarrow \W[1,r]\hookrightarrow \Lo[\infty]$ (\cite[Thm. 1.6.1]{hen81}). This evidently implies the desired boundedness of $\uep$. To verify the bound on $\nabla\cep$ we employ semigroup estimates for the Neumann heat semigroup to find $C_{2}>0$ satisfying
\begin{align*}
\|\nabla\cep(\cdot,t)\|_{\Lo[\infty]}\leq \|\nabla c_0\|_{\Lo[\infty]}&+C_{2}\int_0^t \big(1+(t-s)^{-\frac{3}{4}}\big)e^{-\lambda(t-s)}\big\|\cep(\cdot,s)\big\|_{\Lo[\infty]}\big\|\nep(\cdot,s)\big\|_{\Lo[4]}\intd s \\&+C_{2}\int_0^t\big(1+(t-s)^{-\frac34}\big)e^{-\lambda(t-s)}\big\|\uep(\cdot,s)\big\|_{\Lo[\infty]}\big\|\nabla\cep(\cdot,s)\big\|_{\Lo[4]}\intd s
\end{align*}
for all $t>0$ and all $\epsi\in(0,1)$. Seeing that we have bounds for $\|\cep\|_{\Lo[\infty]}$, $\|\nep\|_{\Lo[4]}$, $\|\uep\|_{\Lo[\infty]}$ and $\|\nabla\cep\|_{\Lo[4]}$ by Lemma~\ref{lem:mass}, Lemma~\ref{lem:moser-it}, the previous part of this lemma and Lemma~\ref{lem:nab-cep-l4-bound}, respectively, we may conclude with standard reasoning.
\end{bew}

\setcounter{equation}{0} 
\section{Estimates for the time derivatives}\label{sec6:time}
For the construction of a limit solution via an Aubin--Lions type argument, we also require some boundedness information on the time derivatives of $\nep,\cep$ and $\uep$, which we establish in the following lemma.
\begin{lemma}\label{lem:time-reg}
Assume that \eqref{eq:prop-S} holds with some $\gamma\in[0,\frac56]$ and let $M>0$ and $s>2$. If $\gamma\in[0,\tfrac12)$ assume that $D$ satisfies \eqref{eq:liminf} with $L(M)>0$ provided by Lemma~\ref{lem:ener-part1-small-g} and that $\|c_0\|_{\Lo[\infty]}\leq M$. If instead $\gamma\in(\tfrac12,\tfrac56]$ assume that $D$ satisfies \eqref{eq:liminf} with $L(M)>0$ provided by Lemma~\ref{lem:ener-like-func-l-g} and that $\|c_0\|_{\Lo[\infty]}\leq M$ and  $\|n_0\|_{\Lo[1]}\leq M$. For all $T>0$ there exists $C=C(M,n_0,c_0,u_0,T)>0$ such that for all $\epsi\in(0,1)$ the solution $(\nep,\cep,\uep)$ of \eqref{approxprob} satisfies
\begin{align*}
\int_0^T\big\|\partial_t\nep\big\|_{(\W[1,2])^*}+\int_0^T\big\|\partial_t\cep\big\|_{(\W[1,2])^*}+\int_0^T\big\|\partial_t\uep\big\|_{(W^{1,s}_{0,\sigma}(\Omega))^*}\leq C.
\end{align*}
\end{lemma}

\begin{bew}
Let $K=K(M,n_0,c_0,u_0)>0$ denote the constant from Lemma~\ref{lem:moser-it}. Then, according to Lemma~\ref{lem:moser-it} and \eqref{eq:Sep}, for any $\varphi\in\W[1,2]$ we may estimate
\begin{align*}
\Big|\intomega n_{\epsi t}\varphi\Big|&=\Big|\intomega\nabla\cdot\big(\Dep(\nep)\nabla\nep-\nep\Sep(x,\nep,\cep)\cdot\nabla\cep-\nep\uep\big)\varphi\Big|\\
&\leq\!\sup_{\epsi\in(0,1)}\!\|\Dep\|_{\LSp{\infty}{(0,K)}}\|\nabla\nep\|_{\Lo[2]}\|\nabla\varphi\|_{\Lo[2]}^2\\&\qquad+K S_0\big(M+1\big)\Big\|\frac{\nabla\cep}{\cep^\nfrac56}\Big\|_{\Lo[6]}\|\cep^{\nfrac56-\gamma}\|_{\Lo[3]}\|\nabla\varphi\|_{\Lo[2]}\\
&\qquad\qquad +K\|\uep\|_{\Lo[2]}\|\nabla\varphi\|_{\Lo[2]}\quad\text{on }(0,\infty)\ \text{for all }\epsi\in(0,1).
\end{align*}
Similarly, for any $\varphi\in\W[1,2]$
\begin{align*}
\Big|\intomega c_{\epsi t}\varphi\Big|&=\Big|\intomega\big(\Delta\cep-\nep\cep-\nabla\cdot(\cep\uep)\big)\varphi\Big|\\
&\leq\!\|\nabla\cep\|_{\Lo[2]}\|\nabla\varphi\|_{\Lo[2]}+K\|\cep\|_{\Lo[2]}\|\varphi\|_{\Lo[2]}+\|\cep\|_{\Lo[\infty]}\|\uep\|_{\Lo[2]}\|\nabla\varphi\|_{\Lo[2]}
\end{align*}
holds on $(0,\infty)$ for all $\epsi\in(0,1)$. Finally, for $\psi\in W_{0,\sigma}^{1,s}(\Omega)$ we find from likewise arguments that
\begin{align*}
\Big|\intomega u_{\epsi t}\cdot\psi\Big|&=\Big|-\intomega\nabla\uep\cdot\nabla\psi-\intomega(\uep\cdot\nabla)\uep\cdot\psi+\intomega\nep\nabla\Phi\cdot\psi\Big|\\&\leq\|\nabla\uep\|_{\Lo[2]}\|\nabla\psi\|_{\Lo[2]}+\!\|\uep\|_{\Lo[2]}\|\nabla\uep\|_{\Lo[2]}\|\psi\|_{\Lo[\infty]}+\!\|\nabla\Phi\|_{\Lo[\infty]}\|\nep\|_{\Lo[1]}\|\psi\|_{\Lo[\infty]}
\end{align*}
is valid on $(0,\infty)$ for all $\epsi\in(0,1)$. Integrating from $0$ to $T$, the claim is then a consequence of the bounds established in Lemma~\ref{lem:moser-it}, Lemma~\ref{lem:nab-cep-l4-bound}, Corollary~\ref{cor:uep-nabcep-linfty}, Lemma~\ref{lem:mass} and Lemma~\ref{lem:uep-bound} with $\gamma_0=\frac56$.
\end{bew}

\setcounter{equation}{0} 
\section{Existence of a global bounded solution. Proof of Theorem~\ref{theo1}}\label{sec7:solution}
The precompactness properties entailed from the bounds presented in the previous sections at hand, we can now pass to the limit $\epsi\searrow0$ and obtain limit functions, which, as illustrated in Lemma~\ref{lem:sol-prop}, solve \eqref{problem} in the standard weak sense.
\begin{lemma}\label{lem:aubin}
Assume that \eqref{eq:prop-S} holds with some $\gamma\in[0,\frac56]$ and let $M>0$. If $\gamma\in[0,\tfrac12)$ assume that $D$ satisfies \eqref{eq:liminf} with $L(M)>0$ provided by Lemma~\ref{lem:ener-part1-small-g} and that $\|c_0\|_{\Lo[\infty]}\leq M$. If instead $\gamma\in(\tfrac12,\tfrac56]$ assume that $D$ satisfies \eqref{eq:liminf} with $L(M)>0$ provided by Lemma~\ref{lem:ener-like-func-l-g} and that $\|c_0\|_{\Lo[\infty]}\leq M$ and  $\|n_0\|_{\Lo[1]}\leq M$. Then there exist $(\epsi_j)_{j\in\N}\subset(0,1)$ and functions
\begin{align*}
n&\in\LSp{\infty}{\Omega\times(0,\infty)}\cap\LSploc{2}{[0,\infty);\W[1,2]},\\
c&\in\LSp{\infty}{(0,\infty);\W[1,\infty]}\\
u&\in\LSpn{\infty}{(0,\infty);W^{1,\infty}_{0,\sigma}(\Omega)}
\end{align*}
such that $\epsi_j\to0$ as $j\to\infty$, $n\geq 0$ and $c>0$ a.e. in $\Omega\times(0,\infty)$, and such that 
\begin{alignat}{2}
\label{eq:nep-strong-conv}
\nep&\to n \qquad &&\text{in }\LSploc{p}{\bomega\times[0,\infty)}\text{ for all }p>1\text{ and a.e. in }\Omega\times(0,\infty),\\
\label{eq:nab-nep-w-conv}
\nabla\nep&\wto\nabla n &&\text{in }\LSploc{2}{\bomega\times[0,\infty);\R^2},\\
\label{eq:cep-strong-conv}
\cep&\to c \qquad &&\text{in }\LSploc{p}{\bomega\times[0,\infty)}\text{ for all }p>1\text{ and a.e. in }\Omega\times(0,\infty),\\
\label{eq:nab-cep-wstar-conv}
\nabla\cep&\wsto \nabla c&& \text{in }\LSp{\infty}{\Omega\times(0,\infty);\R^2},\\
\label{eq:uep-strong-conv}
\uep&\to u \qquad &&\text{in }\LSploc{p}{\bomega\times[0,\infty);\R^2}\text{ for all }p>1\text{ and a.e. in }\Omega\times(0,\infty),\\
\label{eq:nab-uep-wstar-conv}
\nabla\uep&\wsto \nabla u&&\text{in }\LSp{\infty}{\Omega\times(0,\infty);\R^{2\times 2}}
\end{alignat}
as $\epsi=\epsi_j\to0$.
\end{lemma}

\begin{bew}
Note that as a consequence of Lemmas~\ref{lem:nab-cep-l4-bound} and \ref{lem:moser-it} $\{\nep\}_{\epsi\in(0,1)}$ is bounded in $\LSp{2}{(0,T);\W[1,2]}$ and by Lemma~\ref{lem:time-reg} $\{n_{\epsi t}\}_{\epsi\in(0,1)}$ is bounded in $\LSp{1}{(0,T);(\W[1,2])^*}$ for any $T>0$. Thus, we may employ an Aubin--Lions type lemma (e.g. \cite[Corollary 8.4]{sim87}) to obtain that $\{\nep\}_{\epsi\in(0,1)}$ is precompact in $\LSp{2}{(0,T);\Lo[2]}$ for any $T>0$, which ensures the existence of $n\in\LSploc{2}{[0,\infty);\Lo[2]}$ and a subsequence $(\epsi_j)_{j\in\N}\subset(0,1)$ satisfying $\nep\to n$ a.e. in $\Omega\times(0,\infty)$ and in $\LSploc{2}{\bomega\times[0,\infty)}$. Extracting an additional subsequence (still denoted by $\epsi_j$) we conclude that \eqref{eq:nab-nep-w-conv} holds from the boundedness of $\{\nabla\nep\}_{\epsi\in(0,1)}$ in $\LSploc{2}{\bomega\times[0,\infty)}$. Moreover, combining the equi-integrability of $\{n_{\epsi_j}^p\}_{j\in\N}$ for any $p>1$ as entailed by Lemma~\ref{lem:moser-it} with Vitali's convergence theorem yields \eqref{eq:nep-strong-conv}. Arguing similarly for $\cep$ and $\uep$, while drawing on the bounds presented in Lemma~\ref{lem:time-reg}, Lemma~\ref{lem:mass}, Corollary~\ref{cor:uep-nabcep-linfty}, and Lemma~\ref{lem:A-uep-Lr} when combined with the embedding $D(A^\beta_r)\hookrightarrow \CSp{1+\theta}{\bomega}\hookrightarrow \W[1,\infty]$ for $r>2$, $\beta\in(\tfrac12,\fracpow)$ and $\theta\in(0,1)$ satisfying $1+\theta<2\beta-\frac2r$, we can extract further subsequences along which also \eqref{eq:cep-strong-conv}--\eqref{eq:nab-uep-wstar-conv} hold. The nonnegativity of $n$ is a consequence of the nonnegativity of $\nep$, while the positivity of $c$ follows from Lemma~\ref{lem:lower-c}.
\end{bew}

\begin{definition}\label{def:weak_sol}
Let $S\in\CSp{2}{\bomega\times[0,\infty)\times(0,\infty);\R^{2\times2}}$. Assume $D\in\bigcup_{\theta\in(0,1)}\CSp{\theta}{[0,\infty)}\cap\CSp{2}{(0,\infty)}$ is positive on $(0,\infty)$ and suppose that the triple $(n,c,u)$ of functions
\begin{align*}
n&\in\LSploc{1}{\bomega\times[0,\infty)},\\
c&\in\LSploc{\infty}{\bomega\times[0,\infty)}\cap\LSploc{1}{[0,\infty);\W[1,1]},\\
u&\in L_{loc}^1\big([0,\infty); W_{0}^{1,1}\!\left(\Omega;\R^2\right)\!\big),
\end{align*}
are such that $n\geq0$ and $c>0$ a.e. in $\Omega\times(0,\infty)$ and that $$D_1(n),\quad n|S(x,n,c)||\nabla c|,\quad n|u|,\quad\text{and}\quad|u\otimes u|\quad \text{belong to }\LSploc{1}{\bomega\times[0,\infty)},$$ where $D_1(s)=\int_0^s D(\sigma)\intd\sigma$ for $s\geq0$. Then $(n,c,u)$ will be called a global weak solution of \eqref{problem}, if $\nabla\cdot u=0$ a.e. in $\Omega\times(0,\infty)$, if
\begin{align}\label{eq:weak-sol-n}
-\intinfomega n\varphi_t&-\intomega n_0\,\varphi(\cdot,0)=\intinfomega D_1(n)\Delta\varphi+\intinfomega n\big(S(x,n,c)\nabla c\cdot\nabla\varphi\big)+\intinfomega n(u\cdot\nabla\varphi)
\end{align}
for all $\varphi\in C_0^\infty\left(\bomega\times[0,\infty)\right)$ satisfying $\frac{\partial\varphi}{\partial\nu}=0$ on $\romega\times(0,\infty)$, if
\begin{align}\label{eq:weak-sol-c}
-\intinfomega c\varphi_t&-\intomega c_0\,\varphi(\cdot,0)=-\intinfomega \nabla c\cdot\nabla\varphi-\intinfomega nc\varphi+\intinfomega c(u\cdot\nabla\varphi)
\end{align}
for all $\varphi\in C_0^\infty\left(\bomega\times[0,\infty)\right)$, and if
\begin{align}\label{eq:weak-sol-u}
-\intinfomega u\cdot\psi_t-\intomega u_0\cdot\psi(\cdot,0)=-\intinfomega\nabla u\cdot\nabla\psi+\intinfomega(u\otimes u)\cdot\nabla\psi+\intinfomega n\nabla\Phi\cdot\psi
\end{align}
for all $\psi\in C_0^\infty\left(\Omega\times[0,\infty);\R^2\right)$ such that $\nabla\cdot\psi=0$ in $\Omega\times(0,\infty)$.
\end{definition}

\begin{lemma}\label{lem:sol-prop}
Assume that \eqref{eq:prop-S} holds with some $\gamma\in[0,\frac56]$ and let $M>0$. If $\gamma\in[0,\tfrac12)$ assume that $D$ satisfies \eqref{eq:liminf} with $L(M)>0$ provided by Lemma~\ref{lem:ener-part1-small-g} and that $\|c_0\|_{\Lo[\infty]}\leq M$. If instead $\gamma\in(\tfrac12,\tfrac56]$ assume that $D$ satisfies \eqref{eq:liminf} with $L(M)>0$ provided by Lemma~\ref{lem:ener-like-func-l-g} and that $\|c_0\|_{\Lo[\infty]}\leq M$ and  $\|n_0\|_{\Lo[1]}\leq M$. Then the functions $(n,c,u)$ obtained in Lemma~\ref{lem:aubin} constitute a global weak solution in the sense of Definition~\ref{def:weak_sol}.
\end{lemma}

\begin{bew}
Let us first note that for any $T>0$ $$\intoTomega \big|\Sep(x,\nep,\cep)\nabla\cep\big|^6\leq S_0^6(M+1)\intoTomega\frac{|\nabla\cep|^6}{\cep^{6\gamma}}\leq S_0^6(M+1)M^{5-6\gamma}\intoTomega\frac{|\nabla\cep|^6}{\cep^5},$$ so that by Lemma~\ref{lem:nab-cep-l4-bound}, Lemma~\ref{lem:aubin} and \eqref{eq:Sep} $\Sep(x,\nep,\cep)\nabla\cep\wto S(x,n,c)\nabla c$ in $\LSploc{6}{\bomega\times[0,\infty);\R^2}$. Then, starting from the weak formulations corresponding to \eqref{approxprob}, we can easily verify that the convergence properties established in Lemma~\ref{lem:aubin} and the properties in \eqref{eq:Dep} and \eqref{eq:Sep} suffice to let $\epsi=\epsi_j\searrow0$ in each of the integrals to arrive at \eqref{eq:weak-sol-n}, \eqref{eq:weak-sol-c} and \eqref{eq:weak-sol-u}, respectively. 
\end{bew}

With the lemma above we are basically finished for the case where the diffusion features a degeneracy at zero. For the part of Theorem~\ref{theo1} requiring the additional hypothesis that $D(0)>0$, we can attain the improved regularity properties by following the arguments presented in the second part of \cite[Lemma 5.1]{winklerChemotaxisStokesInteractionVery2022}.

\begin{lemma}\label{lem:more-reg-D-positive}
Assume that \eqref{eq:prop-S} holds with some $\gamma\in[0,\frac56]$ and let $M>0$. If $\gamma\in[0,\tfrac12)$ assume that $D$ satisfies \eqref{eq:liminf} with $L(M)>0$ provided by Lemma~\ref{lem:ener-part1-small-g} and that $\|c_0\|_{\Lo[\infty]}\leq M$. If instead $\gamma\in(\tfrac12,\tfrac56]$ assume that $D$ satisfies \eqref{eq:liminf} with $L(M)>0$ provided by Lemma~\ref{lem:ener-like-func-l-g} and that $\|c_0\|_{\Lo[\infty]}\leq M$ and  $\|n_0\|_{\Lo[1]}\leq M$.  Suppose moreover that $D(0)>0$, then the functions $n,c,u$ obtained in Lemma~\ref{lem:aubin} satisfy the enhanced regularity conditions expressed in \eqref{classical} and, together with some $P\in\CSp{1,0}{\Omega\times(0,\infty)}$, solve \eqref{problem} in the classical sense.
\end{lemma}

\begin{bew}
By assumption $D$ is uniformly positive in $[0,\infty)$ and by Lemma~\ref{lem:nab-cep-l4-bound} and Corollary~\ref{cor:uep-nabcep-linfty} there is $C_1=C_1(M,n_0,c_0,u_0)>0$ such that for each $T>0$ \begin{align*}\big\|\Sep(x,\nep,\cep)\nabla\cep+\uep\big\|_{\LSp{6}{\Omega\times(0,T)}}\leq C_1\quad\text{for all }t>0.\end{align*} Hence, well-known results on the Hölder regularity in parabolic equations (e.g. \cite[Theorem 1.3]{PorzVesp93}) entail the existence of $\theta_1=\theta_1(T,D,n_0,c_0,u_0)\in(0,1)$ such that $\{\nep\}_{\epsi\in(0,1)}$ is bounded in $\CSp{\theta_1,\frac{\theta_1}{2}}{\bomega\times[0,T]}$. Furthermore, adjusting the arguments presented in \cite[Lemma 4.4]{winklerDoesLerayStructure2022} and \cite[Lemma 3.16]{heihoffTwoNewFunctional2023} to our setting, we first find that for all $\beta_1\in(\tfrac12,1)$ and all $r_1\in(2,\infty)$ there is and $C_2=C_2(M,n_0,c_0,u_0,r_1,\beta_1)>0$ such that
\begin{align}\label{eq:more-reg-D-pos-eq1}
\big\|A^{\beta_1}\big(\uep(\cdot,t_2)-\uep(\cdot,t_1)\big)\big\|_{\Lo[r_1]}\leq C_2 (t_2-t_1)^{1-\beta_1}
\end{align}
holds for all $t_2>t_1>0$ and all $\epsi\in(0,1)$. Picking $\beta_2\in(\frac12,\fracpow)$ and $r_2\in(2,\infty)$ such that $2\beta_2-\frac2{r_2}>1$ and $r_2>\frac{4}{\beta_1}$ hold, the latter ensuring that with $\theta_2:=\frac{1-\beta_1}{2}$ we have $1+\theta_2<2\beta_2-\frac2{r_2}$ and hence $D(A^{\beta_2}_{r_2})\hookrightarrow \CSp{1+\theta_2}{\bomega;\R^2}$ (cf. \cite{gig85} and \cite[p.39]{hen81}), we can conclude from \eqref{eq:more-reg-D-pos-eq1} and Lemma~\ref{lem:A-uep-Lr} that there is $C_3=C_3(M,n_0,c_0,u_0)>0$ such that $$\big\|\uep\big\|_{\CSp{\theta_2}{[t,t+1];\CSpnl{1+\theta_2}{\bomega;\R^2}}}\leq C_3\quad\text{for all }t>0\text{ and }\epsi\in(0,1).$$ Thereafter, Schauder theory for the Stokes evolution equations (\cite[Theorem 1.1]{Solonnikov2007} and \cite[V.1.5.1]{sohr}) shows that for any $\tau>0$ and $T>\tau$ there is $\theta_3=\theta_3(\tau,T,D,n_0,c_0,u_0)>0\in(0,1)$ such that $\{\uep\}_{\epsi\in(0,1)}$ is moreover bounded in $\CSp{2+\theta_3,1+\frac{\theta_3}{2}}{\bomega\times[\tau,T];\R^2}$. In similar fashion, parabolic regularity theory (e.g. \cite[IV.5.3 and III.5.1]{LSU}) yields $\theta_4=\theta_4(T,D,n_0,c_0,u_0)\in(0,1)$ and $\theta_5=\theta_5(\tau,T,D,n_0,c_0,u_0)>0$ such that $\{\cep\}$ is bounded in $\CSp{\theta_4,\tfrac{\theta_4}{2}}{\bomega\times[0,T]}$ and in $\CSp{2+\theta_5,1+\frac{\theta_5}{2}}{\bomega\times[\tau,T]}$. In view of \eqref{eq:nep-strong-conv}, \eqref{eq:cep-strong-conv}, \eqref{eq:uep-strong-conv} and the Arzelà--Ascoli theorem we hence conclude that $\nep\to n$ in $\CSploc{0}{\bomega\times[0,\infty)}$, that $\cep\to c$ in $\CSploc{0}{\bomega\times[0,\infty)}\cap\CSp{2,1}{\bomega\times(0,\infty)}$ and that $\uep\to u$ in $\CSploc{0}{\bomega\times[0,\infty);\R^2}\cap\CSp{2,1}{\bomega\times(0,\infty);\R^2}$ as $\epsi_j\to 0$. Constructing the pressure $P\in\CSp{1,0}{\Omega\times(0,\infty)}$ by standard arguments (\cite{sohr}) we hence find, that the third and second equations of \eqref{problem} are satisfied in the classical sense. With these regularity properties established, arguing similar to \cite[Lemma 5.7]{caolan16_smalldatasol3dnavstokes}, we find that since $n$ is a bounded generalized solution of $n_t=\nabla\cdot(D(n)\nabla n-n S(x,n,c)\nabla c - nu)$ in $\Omega\times(0,\infty)$ with $(D(n)\nabla n-n S(x,n,c)\nabla c-nu)\cdot\nu\vert_{\romega}=0$ and $n(\cdot,0)=n_0\in\W[1,\infty]$ as considered in \cite[Theorem 1.1]{liebermanHolderContinuityGradient1987} and \cite[IV.5.3]{LSU} we also conclude that $n$ actually belongs to $\CSp{0}{\bomega\times[0,\infty)}\cap\CSp{2,1}{\bomega\times(0,\infty)}$ and solves its respective equation classically in $\Omega\times(0,\infty)$.
\end{bew}

We just have to combine the results from the last two lemmas above in order to see that Theorem~\ref{theo1} holds.

\begin{proof}[\textbf{Proof of Theorem~\ref{theo1}}]
The theorem is an evident consequence of Lemma~\ref{lem:sol-prop} and Lemma~\ref{lem:more-reg-D-positive}.
\end{proof}

Finally, let us briefly cover the arguments for the proof of Corollary~\ref{cor:2}.
\begin{proof}[\textbf{Proof of Corollary~\ref{cor:2}}]
Even without the uniform positivity of $D$ -- crucially exploited in Lemma~\ref{lem:more-reg-D-positive} -- we can still rely on quite general Hölder regularity properties for bounded weak solutions in the porous medium setting $D(s)=s^{m-1}$, $m\in(1,2]$. In fact in this context, we can explicitly pick $D_\epsi(s)=(s+\epsi^\frac{1}{m-1})^{m-1}$ for the proof of Theorem~\ref{theo1} noting that this way the conditions in \eqref{eq:Dep} are clearly satisfied. Then, denoting by $(\nep,\cep,\uep)_{\epsi\in(0,1)}$ the solution family of \eqref{approxprob} and by $(n,c,u)$ the bounded weak solution obtained from the limiting procedure in Lemma~\ref{lem:aubin}, we start by mimicking the argument of Lemma~\ref{lem:more-reg-D-positive} for $\cep$ and $\uep$, to conclude that actually $(c,u)\in\CSp{0}{\bomega\times[0,\infty)}\times\CSp{0}{\bomega\times[0,\infty);\R^2}$ by means of the Arzelà--Ascoli theorem. To verify the continuity of $n$, we introduce $\wep:=\nep+\epsi^{\frac{1}{m-1}}$ and $\Phi(s):=\int_0^s \sigma^{m-1}\intd\sigma=\frac{1}{m}s^{m}$ as well as $$a_\epsi(x,t):=\nep\big(\nep+\epsi^{\frac{1}{m-1}}\big)^{-1}\Sep(x,\nep,\cep)\nabla\cep-\nep\big(\nep+\epsi^{\frac{1}{m-1}}\big)^{-1}\uep$$ and observe that then $\wep$ solves $$w_{\epsi t}=\Delta\Phi(\wep)+\nabla\cdot \big(a_\epsi(x,t) \wep\big)\quad\text{in }\Omega\times(0,\infty).$$ Here, \eqref{eq:prop-S} and $\gamma\leq\frac{5}{6}$ together with Lemma~\ref{lem:mass} and the bounds established in Lemma~\ref{lem:nab-cep-l4-bound} and Corollary~\ref{cor:uep-nabcep-linfty} entail that for any $T>0$ there is $C=C(M,n_0,c_0,u_0,T)>0$ such that $$\|a_\epsi\|_{\LSpn{6}{\Omega\times(0,T);\R^2}}^6\leq S_0^6(M+1)\|c_0\|^{5-6\gamma}_{\Lo[\infty]}\intoTomega\frac{|\nabla\cep|^6}{\cep^5}+\intoTomega|\uep|^6\leq C\quad\text{for all }\epsi\in(0,1).$$ Moreover, $\epsi\in(0,1)$ and Lemma~\ref{lem:moser-it} ensure that with some $K=K(M,n_0,c_0,u_0)>0$ we have $$\|\wep(\cdot,t)\|_{\Lo[\infty]}\leq\|\nep(\cdot,t)+1\|_{\Lo[\infty]}\leq K\quad\text{for all }t\in(0,T)\quad\text{ and }\epsi\in(0,1).$$ Since these bounds are independent of $\epsi\in(0,1)$, we obtain Hölder bounds uniform in $\epsi\in(0,1)$ for $\wep$ from \cite[Theorem 1.8 and Corollary 1.9]{TB23_hoeldertaxis}. Drawing once more on the Arzelà--Ascoli theorem we may hence conclude that indeed also $n$ is continuous in $\bomega\times[0,\infty)$.
\end{proof}

\section*{Acknowledgements}
The author acknowledges support of the {\em Deutsche Forschungsgemeinschaft} (Project No.~462888149) and the {\em Research Institute for Mathematical Sciences}, an International Joint Usage/Research Center located in Kyoto University.

\footnotesize{
\setlength{\bibsep}{3pt plus 0.5ex}

}


\begin{thebibliography}{56}
\providecommand{\natexlab}[1]{#1}

\bibitem[Bellomo et~al.(2015)Bellomo, Bellouquid, Tao, and
  Winkler]{BBWT15}{https://doi.org/10.1142/S021820251550044X}
N.~Bellomo, A.~Bellouquid, Y.~Tao, and M.~Winkler.
\newblock Toward a mathematical theory of {{Keller-Segel}} models of pattern
  formation in biological tissues.
\newblock \emph{Math. Models Methods Appl. Sci.}, 25\penalty0 (9):\penalty0
  1663--1763, 2015.

\bibitem[Black(2023)]{TB23_hoeldertaxis}{}
T.~Black.
\newblock Refining h\"older regularity theory in degenerate drift-diffusion
  equations.
\newblock 2023.
\newblock Preprint.

\bibitem[Brenner and
  Scott(1994)]{brennerMathematicalTheoryFinite1994}{https://doi.org/10.1007/978-1-4757-4338-8}
S.~C. Brenner and L.~R. Scott.
\newblock \emph{The {{Mathematical Theory}} of {{Finite Element Methods}}},
  Volume~15 of \emph{Texts in {{Applied Mathematics}}}.
\newblock {Springer New York}, 1994.

\bibitem[Cao and
  Ishida(2014)]{caoGlobalintimeBoundedWeak2014}{https://doi.org/10.1088/0951-7715/27/8/1899}
X.~Cao and S.~Ishida.
\newblock Global-in-time bounded weak solutions to a degenerate quasilinear
  {{Keller}}--{{Segel}} system with rotation.
\newblock \emph{Nonlinearity}, 27\penalty0 (8):\penalty0 1899--1913, 2014.

\bibitem[Cao and
  Lankeit(2016)]{caolan16_smalldatasol3dnavstokes}{https://doi.org/10.1007/s00526-016-1027-2}
X.~Cao and J.~Lankeit.
\newblock Global classical small-data solutions for a three-dimensional
  chemotaxis {{Navier-Stokes}} system involving matrix-valued sensitivities.
\newblock \emph{Calc. Var. Partial Differential Equations}, 55\penalty0
  (4):\penalty0 Paper No. 107, 2016.

\bibitem[Cao and
  Wang(2015)]{CaoWang-GlobClass-DCDSB15}{https://doi.org/10.3934/dcdsb.2015.20.3235}
X.~Cao and Y.~Wang.
\newblock Global classical solutions of a {{3D}} chemotaxis-{{Stokes}} system
  with rotation.
\newblock \emph{Discret. Contin. Dyn. Syst. - Ser. B}, 20\penalty0
  (9):\penalty0 3235--3254, 2015.

\bibitem[Cie{\'s}lak and
  Stinner(2012)]{CieslakStinner-JDE12}{https://doi.org/10.1016/j.jde.2012.01.045}
T.~Cie{\'s}lak and C.~Stinner.
\newblock Finite-time blowup and global-in-time unbounded solutions to a
  parabolic-parabolic quasilinear {{Keller-Segel}} system in higher dimensions.
\newblock \emph{J. Differential Equations}, 252\penalty0 (10):\penalty0
  5832--5851, 2012.

\bibitem[Duan and
  Xiang(2014)]{duanNoteGlobalExistence2014}{https://doi.org/10.1093/imrn/rns270}
R.~Duan and Z.~Xiang.
\newblock A {{Note}} on {{Global Existence}} for the {{Chemotaxis}}--{{Stokes
  Model}} with {{Nonlinear Diffusion}}.
\newblock \emph{International Mathematics Research Notices}, 2014\penalty0
  (7):\penalty0 1833--1852, 2014.

\bibitem[Francesco et~al.(2010)Francesco, Lorz, and
  Markowich]{FrancescoLorzMarkowich-DCDS10}{https://doi.org/10.3934/dcds.2010.28.1437}
M.~Francesco, A.~Lorz, and P.~A. Markowich.
\newblock Chemotaxis-fluid coupled model for swimming bacteria with nonlinear
  diffusion: {{Global}} existence and asymptotic behavior.
\newblock \emph{Discret. Contin. Dyn. Syst.}, 28\penalty0 (4):\penalty0
  1437--1453, 2010.

\bibitem[Fuest(2023)]{fuestChemotaxisFluidSystems2023}{https://doi.org/10.3934/dcdsb.2022232}
M.~Fuest.
\newblock Chemotaxis(-fluid) systems with logarithmic sensitivity and slow
  consumption: {{Global}} generalized solutions and eventual smoothness.
\newblock \emph{Discrete Contin. Dyn. Syst. Ser. B}, 28\penalty0 (10):\penalty0
  5177--5202, 2023.

\bibitem[Galdi(2011)]{galdiIntroductionMathematicalTheory2011}{https://doi.org/10.1007/978-0-387-09620-9}
G.~Galdi.
\newblock \emph{An {{Introduction}} to the {{Mathematical Theory}} of the
  {{Navier-Stokes Equations}}: {{Steady-State Problems}}}.
\newblock Springer {{Monographs}} in {{Mathematics}}. {Springer New York},
  2011.
\newblock ISBN 978-0-387-09619-3 978-0-387-09620-9.

\bibitem[Giga(1985)]{gig85}{https://doi.org/10.1007/BF00276874}
Y.~Giga.
\newblock Domains of fractional powers of the {{Stokes}} operator in
  {{L}}{\textsubscript{r}} spaces.
\newblock \emph{Arch. Rational Mech. Anal.}, 89\penalty0 (3):\penalty0
  251--265, 1985.

\bibitem[Giga(1986)]{gig86}{https://doi.org/10.1016/0022-0396(86)90096-3}
Y.~Giga.
\newblock Solutions for semilinear parabolic equations in
  {{L\textsuperscript{p}}} and regularity of weak solutions of the
  {{Navier--Stokes}} system.
\newblock \emph{J. Differential Equations}, 62\penalty0 (2):\penalty0 186--212,
  1986.

\bibitem[Heihoff(2023)]{heihoffTwoNewFunctional2023}{https://doi.org/10.1137/22M1531178}
F.~Heihoff.
\newblock Two {{New Functional Inequalities}} and {{Their Application}} to the
  {{Eventual Smoothness}} of {{Solutions}} to a
  {{Chemotaxis-Navier}}\textendash{{Stokes System}} with {{Rotational Flux}}.
\newblock \emph{SIAM J. Math. Anal.}, 55\penalty0 (6):\penalty0 7113--7154,
  2023.

\bibitem[Henry(1981)]{hen81}{https://doi.org/10.1007/BFb0089647}
D.~Henry.
\newblock \emph{Geometric {{Theory}} of {{Semilinear Parabolic Equations}}},
  Volume 840 of \emph{Lecture {{Notes}} in {{Mathematics}}}.
\newblock {Springer Berlin Heidelberg}, 1981.

\bibitem[Hillen and
  Painter(2009)]{HP09}{https://doi.org/10.1007/s00285-008-0201-3}
T.~Hillen and K.~J. Painter.
\newblock A user's guide to {{PDE}} models for chemotaxis.
\newblock \emph{J. Math. Biol.}, 58\penalty0 (1-2):\penalty0 183--217, 2009.

\bibitem[Ishida et~al.(2014)Ishida, Seki, and
  Yokota]{ISY-quasilin-pp-JDE14}{https://doi.org/10.1016/j.jde.2014.01.028}
S.~Ishida, K.~Seki, and T.~Yokota.
\newblock Boundedness in quasilinear {{Keller-Segel}} systems of
  parabolic-parabolic type on non-convex bounded domains.
\newblock \emph{J. Differential Equations}, 256\penalty0 (8):\penalty0
  2993--3010, 2014.

\bibitem[Jost(2007)]{Jost-PDE07}{https://doi.org/10.1007/978-0-387-49319-0}
J.~Jost.
\newblock \emph{Partial {{Differential Equations}}}, Volume 214 of
  \emph{Graduate {{Texts}} in {{Mathematics}}}.
\newblock {Springer New York}, 2007.
\newblock ISBN 978-0-387-49318-3.

\bibitem[Keller and
  Segel(1970)]{KS70}{https://doi.org/10.1016/0022-5193(70)90092-5}
E.~F. Keller and L.~A. Segel.
\newblock Initiation of slime mold aggregation viewed as an instability.
\newblock \emph{J. Theor. Biol.}, 26\penalty0 (3):\penalty0 399--415, 1970.

\bibitem[Kim(2023)]{kimGlobalSolutionsChemotaxisfluid2023}{https://doi.org/10.3934/dcdsb.2023040}
D.~Kim.
\newblock Global solutions for chemotaxis-fluid systems with singular
  chemotactic sensitivity.
\newblock \emph{Discrete Contin. Dyn. Syst. Ser. B}, 28\penalty0 (10):\penalty0
  5380--5395, 2023.

\bibitem[Lady{\v z}enskaja et~al.(1968)Lady{\v z}enskaja, Solonnikov, and
  Ural'ceva]{LSU}{https://doi.org/10.1090/mmono/023}
O.~A. Lady{\v z}enskaja, V.~A. Solonnikov, and N.~N. Ural'ceva.
\newblock \emph{Linear and Quasilinear Equations of Parabolic Type}.
\newblock Translations of Mathematical Monographs. {American Mathematical
  Society}, 1968.
\newblock ISBN 978-0-8218-1573-1.

\bibitem[Lankeit(2017)]{Lan17-LocBddGlobSolNonlinDiff-JDE}{https://doi.org/10.1016/j.jde.2016.12.007}
J.~Lankeit.
\newblock Locally bounded global solutions to a chemotaxis consumption model
  with singular sensitivity and nonlinear diffusion.
\newblock \emph{J. Differential Equations}, 262\penalty0 (7):\penalty0
  4052--4084, 2017.

\bibitem[Lankeit and
  Winkler(2017)]{LanWin2017}{https://doi.org/10.1007/s00030-017-0472-8}
J.~Lankeit and M.~Winkler.
\newblock A generalized solution concept for the {{Keller-Segel}} system with
  logarithmic sensitivity: Global solvability for large nonradial data.
\newblock \emph{NoDEA Nonlinear Differ. Equ. Appl.}, 24\penalty0 (4):\penalty0
  Art. 49, 33, 2017.

\bibitem[Lankeit and
  Winkler(2023)]{lankeitDepletingSignalAnalysis2023}{https://doi.org/10.1111/sapm.12625}
J.~Lankeit and M.~Winkler.
\newblock Depleting the signal: {{Analysis}} of chemotaxis-consumption
  models---{{A}} survey.
\newblock \emph{Stud. Appl. Math.}, 151\penalty0 (4):\penalty0 1197--1229,
  2023.

\bibitem[Lieberman(1987)]{liebermanHolderContinuityGradient1987}{https://doi.org/10.1007/BF01774284}
G.~M. Lieberman.
\newblock H{\"o}lder continuity of the gradient of solutions of uniformly
  parabolic equations with conormal boundary conditions.
\newblock \emph{Ann. Mat. Pura Appl.}, 148\penalty0 (1):\penalty0 77--99, 1987.

\bibitem[Liu(2023)]{liuGlobalClassicalSolvability2023}{https://doi.org/10.1142/S0218202523400031}
J.~Liu.
\newblock Global classical solvability and stabilization in a two-dimensional
  chemotaxis--fluid system with sub-logarithmic sensitivity.
\newblock \emph{Math. Models Methods Appl. Sci.}, 33\penalty0 (11):\penalty0
  2271--2311, 2023.

\bibitem[Mizoguchi and
  Souplet(2014)]{MS14}{https://doi.org/10.1016/j.anihpc.2013.07.007}
N.~Mizoguchi and P.~Souplet.
\newblock Nondegeneracy of blow-up points for the parabolic {{Keller-Segel}}
  system.
\newblock \emph{Ann. Inst. H. Poincar\'e Anal. Non Lin\'eaire}, 31\penalty0
  (4):\penalty0 851--875, 2014.

\bibitem[Nagai and Senba(1998)]{NagSen_JAMA98}{}
T.~Nagai and T.~Senba.
\newblock Global existence and blow-up of radial solutions to a
  parabolic-elliptic system of chemotaxis.
\newblock \emph{Adv. Math. Sci. Appl.}, 8\penalty0 (1):\penalty0 145--156,
  1998.

\bibitem[Nagai et~al.(1997)Nagai, Senba, and Yoshida]{NSY97}{}
T.~Nagai, T.~Senba, and K.~Yoshida.
\newblock Application of the {{Trudinger-Moser}} inequality to a parabolic
  system of chemotaxis.
\newblock \emph{Funkcial. Ekvac.}, 40\penalty0 (3):\penalty0 411--433, 1997.

\bibitem[Porzio and
  Vespri(1993)]{PorzVesp93}{https://doi.org/10.1006/jdeq.1993.1045}
M.~M. Porzio and V.~Vespri.
\newblock H{\"o}lder estimates for local solutions of some doubly nonlinear
  degenerate parabolic equations.
\newblock \emph{J. Differential Equations}, 103\penalty0 (1):\penalty0
  146--178, 1993.

\bibitem[Quittner and Souplet(2007)]{QS07}{}
P.~Quittner and P.~Souplet.
\newblock \emph{Superlinear Parabolic Problems}.
\newblock Birkh\"auser {{Advanced Texts}}: {{Basler Lehrb\"ucher}}.
  {Birkh\"auser Verlag, Basel}, 2007.
\newblock ISBN 978-3-7643-8441-8.

\bibitem[Rosen(1978)]{ROSEN1978}{https://doi.org/10.1007/BF02460738}
G.~Rosen.
\newblock Steady-state distribution of bacteria chemotactic toward oxygen.
\newblock \emph{Bull. Math. Biol.}, 40\penalty0 (5):\penalty0 671--674, 1978.

\bibitem[Simon(1987)]{sim87}{https://doi.org/10.1007/BF01762360}
J.~Simon.
\newblock Compact sets in the space {{L\textsuperscript{p}}}(0,{{T}};{{B}}).
\newblock \emph{Ann. Mat. Pura Appl. (4)}, 146:\penalty0 65--96, 1987.

\bibitem[Sohr(2001)]{sohr}{https://doi.org/10.1007/978-3-0348-8255-2}
H.~Sohr.
\newblock \emph{The {{Navier-Stokes}} Equations}.
\newblock Birkh{\"a}user {{Advanced Texts}}: {{Basler Lehrb{\"u}cher}}.
  {Birkh{\"a}user Verlag, Basel}, 2001.
\newblock ISBN 3-7643-6545-5.

\bibitem[Solonnikov(2007)]{Solonnikov2007}{https://doi.org/10.1090/trans2/220/08}
V.~A. Solonnikov.
\newblock Schauder estimates for the evolutionary generalized {{Stokes}}
  problem.
\newblock In \emph{Nonlinear Equations and Spectral Theory}, Volume 220 of
  \emph{Amer. {{Math}}. {{Soc}}. {{Transl}}. {{Ser}}. 2}, pp. 165--200. {Amer.
  Math. Soc., Providence, RI}, 2007.

\bibitem[Tao(2011)]{Tao-consumption_JMAA11}{https://doi.org/10.1016/j.jmaa.2011.02.041}
Y.~Tao.
\newblock Boundedness in a chemotaxis model with oxygen consumption by
  bacteria.
\newblock \emph{J. Math. Anal. Appl.}, 381\penalty0 (2):\penalty0 521--529,
  2011.

\bibitem[Tao and
  Winkler(2011)]{tao_winkler_chemohapto11-siam11}{https://doi.org/10.1137/100802943}
Y.~Tao and M.~Winkler.
\newblock A chemotaxis-haptotaxis model: The roles of nonlinear diffusion and
  logistic source.
\newblock \emph{SIAM J. Math. Anal.}, 43\penalty0 (2):\penalty0 685--704, 2011.

\bibitem[Tao and
  Winkler(2012)]{TaoWin-GlobExAndBdd-DCDS12}{https://doi.org/10.3934/dcds.2012.32.1901}
Y.~Tao and M.~Winkler.
\newblock Global existence and boundedness in a {{Keller-Segel-Stokes}} model
  with arbitrary porous medium diffusion.
\newblock \emph{Discrete Contin. Dyn. Syst.}, 32\penalty0 (5):\penalty0
  1901--1914, 2012.

\bibitem[Tao and
  Winkler(2013)]{TaoWin-LocBddGlobSol-AnnInstHP13}{https://doi.org/10.1016/j.anihpc.2012.07.002}
Y.~Tao and M.~Winkler.
\newblock Locally bounded global solutions in a three-dimensional
  chemotaxis-{{Stokes}} system with nonlinear diffusion.
\newblock \emph{Ann. Inst. H. Poincar{\'e} Anal. Non Lin{\'e}aire}, 30\penalty0
  (1):\penalty0 157--178, 2013.

\bibitem[Tuval et~al.(2005)Tuval, Cisneros, Dombrowski, Wolgemuth, Kessler, and
  Goldstein]{tuval2005bacterial}{https://doi.org/10.1073/pnas.0406724102}
I.~Tuval, L.~Cisneros, C.~Dombrowski, C.~W. Wolgemuth, J.~O. Kessler, and R.~E.
  Goldstein.
\newblock Bacterial swimming and oxygen transport near contact lines.
\newblock \emph{Proc. Natl. Acad. Sci. U.S.A.}, 102\penalty0 (7):\penalty0
  2277--2282, 2005.

\bibitem[Wang(2023)]{wangGlobalBoundedSolution2023}{https://doi.org/10.1007/s10440-023-00599-x}
J.~Wang.
\newblock Global {{Bounded Solution}} in a {{Chemotaxis-Stokes Model}} with
  {{Porous Medium Diffusion}} and {{Singular Sensitivity}}.
\newblock \emph{Acta Appl Math}, 187\penalty0 (1):\penalty0 7, 2023.

\bibitem[Winkler()]{winklerLpboundLocalSensing}{}
M.~Winkler.
\newblock {$L^p$} bounds in the two-dimensional {N}avier--{S}tokes system and
  application to blow-up suppression in chemotaxis-fluid systems accounting for
  local sensing.
\newblock preprint.

\bibitem[Winkler(2012)]{win_fluid_CPDE12}{https://doi.org/10.1080/03605302.2011.591865}
M.~Winkler.
\newblock Global large-data solutions in a chemotaxis-({{Navier-}}){{Stokes}}
  system modeling cellular swimming in fluid drops.
\newblock \emph{Commun. Partial. Differ. Equ.}, 37\penalty0 (2):\penalty0
  319--351, 2012.

\bibitem[Winkler(2013)]{Win13pure}{https://doi.org/10.1016/j.matpur.2013.01.020}
M.~Winkler.
\newblock Finite-time blow-up in the higher-dimensional parabolic--parabolic
  {{Keller}}--{{Segel}} system.
\newblock \emph{J. Math. Pures Appl.}, 100\penalty0 (5):\penalty0 748--767,
  2013.

\bibitem[Winkler(2014)]{win-stab2d-ArchRatMechAna12}{https://doi.org/10.1007/s00205-013-0678-9}
M.~Winkler.
\newblock Stabilization in a two-dimensional chemotaxis-{{Navier-Stokes}}
  system.
\newblock \emph{Arch. Ration. Mech. Anal.}, 211\penalty0 (2):\penalty0
  455--487, 2014.

\bibitem[Winkler(2015{\natexlab{a}})]{Win-ct_fluid_3d-CPDE15}{https://doi.org/10.1007/s00526-015-0922-2}
M.~Winkler.
\newblock Boundedness and large time behavior in a three-dimensional
  chemotaxis-{{Stokes}} system with nonlinear diffusion and general
  sensitivity.
\newblock \emph{Calc. Var. Partial Differential Equations}, 54\penalty0
  (4):\penalty0 3789--3828, 2015{\natexlab{a}}.

\bibitem[Winkler(2015{\natexlab{b}})]{win15_chemorot}{https://doi.org/10.1137/140979708}
M.~Winkler.
\newblock Large-data global generalized solutions in a chemotaxis system with
  tensor-valued sensitivities.
\newblock \emph{SIAM J. Math. Anal.}, 47\penalty0 (4):\penalty0 3092--3115,
  2015{\natexlab{b}}.

\bibitem[Winkler(2018)]{Winkler18_stokesrot}{https://doi.org/10.1007/s00028-018-0440-8}
M.~Winkler.
\newblock Global mass-preserving solutions in a two-dimensional
  chemotaxis-{{Stokes}} system with rotational flux components.
\newblock \emph{Journal of Evolution Equations}, 2018.

\bibitem[Winkler(2019)]{winklerThreedimensionalKellerSegel2019}{https://doi.org/10.1016/j.jfa.2018.12.009}
M.~Winkler.
\newblock A three-dimensional
  {{Keller}}\textendash{{Segel}}\textendash{{Navier}}\textendash{{Stokes}}
  system with logistic source: {{Global}} weak solutions and asymptotic
  stabilization.
\newblock \emph{J. Funct. Anal.}, 276\penalty0 (5):\penalty0 1339--1401, 2019.

\bibitem[Winkler(2020)]{winklerSmallMassSolutionsTwoDimensional2020}{https://doi.org/10.1137/19M1264199}
M.~Winkler.
\newblock Small-{{Mass Solutions}} in the {{Two-Dimensional Keller--Segel
  System Coupled}} to the {{Navier--Stokes Equations}}.
\newblock \emph{SIAM J. Math. Anal.}, 52\penalty0 (2):\penalty0 2041--2080,
  2020.

\bibitem[Winkler(2022{\natexlab{a}})]{winklerApproachingLogarithmicSingularities2022}{https://doi.org/10.3934/dcdsb.2022009}
M.~Winkler.
\newblock Approaching logarithmic singularities in quasilinear
  chemotaxis-consumption systems with signal-dependent sensitivities.
\newblock \emph{Discrete Continuous Dyn. Syst. Ser. B.}, 27\penalty0
  (11):\penalty0 6565, 2022{\natexlab{a}}.

\bibitem[Winkler(2022{\natexlab{b}})]{winklerChemotaxisStokesInteractionVery2022}{https://doi.org/10.1515/ans-2022-0004}
M.~Winkler.
\newblock Chemotaxis-{{Stokes}} interaction with very weak diffusion
  enhancement: {{Blow-up}} exclusion via detection of absorption-induced
  entropy structures involving multiplicative couplings.
\newblock \emph{Adv. Nonlinear Stud.}, 22\penalty0 (1):\penalty0 88--117,
  2022{\natexlab{b}}.

\bibitem[Winkler(2022{\natexlab{c}})]{winklerDoesLerayStructure2022}{https://doi.org/10.4171/JEMS/1226}
M.~Winkler.
\newblock Does {{Leray}}'s structure theorem withstand buoyancy-driven
  chemotaxis-fluid interaction?
\newblock \emph{J. Eur. Math. Soc.}, 2022{\natexlab{c}}.

\bibitem[Winkler(2024)]{Win24_ev-reg-2d-doubly-deg-nutrient}{}
M.~Winkler.
\newblock Eventual regularity in a two-dimensional doubly degenerate nutrient
  taxis system.
\newblock 2024.
\newblock Preprint.

\bibitem[Xue(2015)]{xueMacroscopicEquationsBacterial2015}{https://doi.org/10.1007/s00285-013-0748-5}
C.~Xue.
\newblock Macroscopic equations for bacterial chemotaxis: Integration of
  detailed biochemistry of cell signaling.
\newblock \emph{J. Math. Biology}, 70\penalty0 (1-2):\penalty0 1--44, 2015.

\bibitem[Xue and
  Othmer(2009)]{XueOthmer-Multiscale-SIAM09}{https://doi.org/10.1137/070711505}
C.~Xue and H.~G. Othmer.
\newblock Multiscale models of taxis-driven patterning in bacterial
  populations.
\newblock \emph{SIAM J. Appl. Math.}, 70\penalty0 (1):\penalty0 133--167, 2009.

\end{thebibliography}
\end{document}